\documentclass{article}

\usepackage[a4paper, total={6in, 8in}]{geometry}
\usepackage[utf8]{inputenc}

\usepackage{amsmath}
\usepackage{amssymb}
\usepackage{authblk}
\usepackage{booktabs}
\usepackage{cleveref}
\usepackage{placeins}
\usepackage{todonotes}
\usepackage{tikz}
\usepackage{pgfplots}
\pgfplotsset{compat=1.17}
\usepackage[acronym,toc,nogroupskip]{glossaries}
\newacronym{bvp}{BVP}{boundary value problem}
\newacronym{csst}{C-SST}{time-con\-tin\-u\-ous simplex space-time}
\newacronym[firstplural=degrees of freedom (DOFs)]{dof}{DOF}{degree of freedom}
\newacronym{dpst}{D-PST}{time-Dis\-con\-tin\-u\-ous prismatic space-time}
\newacronym{eim}{EIM}{empirical interpolation method}
\newacronym{emum}{EMUM}{elastic mesh update method}
\newacronym{fem}{FEM}{finite element method}
\newacronym{fom}{FOM}{full-order model}
\newacronym{gls}{GLS}{Galerkin-least-squares}
\newacronym{lbb}{LBB}{Ladyzhenskaya–Babuška–Brezzi}
\newacronym{lhs}{LHS}{left-hand side}
\newacronym{mor}{MOR}{model order reduction}
\newacronym{pc}{PC}{polycarbonate}
\newacronym{pde}{PDE}{partial differential equation}
\newacronym{pgd}{PGD}{proper generalized decomposition}
\newacronym{pod}{POD}{proper orthogonal decomposition}
\newacronym{rb}{RB}{Reduced Basis}
\newacronym{rhs}{RHS}{right-hand side}
\newacronym{rom}{ROM}{reduced-order model}
\newacronym{uq}{UQ}{uncertainty quantification}
\usepackage[url=false,eprint=false,giveninits=true]{biblatex} 
\usepackage[font=footnotesize,labelfont=bf,position=bottom, skip=0.5ex]{subcaption}
\usepackage{siunitx}
\usepackage{tikz-dimline}

\title{Model Order Reduction for Deforming Domain Problems in a Time-Continuous Space-Time Setting}
\author[1,*]{Fabian Key}
\author[2]{Max von Danwitz}
\author[3]{Francesco Ballarin}
\author[4]{Gianluigi Rozza}
\affil[1]{Institute of Lightweight Design and Structural Biomechanics (ILSB), TU Wien, Vienna, Austria}
\affil[2]{Institute for Mathematics and Computer-Based Simulation (IMCS), Universität der Bundeswehr München, Neubiberg, Germany}
\affil[3]{Department of Mathematics and Physics, Università Cattolica del Sacro Cuore, Brescia, Italy}
\affil[4]{Mathematics Area, mathLab, Scuola Internazionale Superiore di Studi Avanzati (SISSA), Trieste, Italy}

\date{}

\newcommand{\vek}[1]{\boldsymbol{#1}}
\newcommand{\mat}[1]{\boldsymbol{#1}}

\newcommand{\lp}{\left(}
\newcommand{\rp}{\right)}
\newcommand{\lb}{\left[}
\newcommand{\rb}{\right]}
\newcommand{\dd}[2]{\frac{\partial#1}{\partial#2}}
\newcommand{\ddt}[1]{#1_t}
\newcommand{\dv}[1]{\vek{\nabla} \cdot #1}
\newcommand{\gr}[1]{\vek{\nabla} #1}

\newcommand{\domainSpace}[1][]{\Omega_{#1}}
\newcommand{\x}{\vek{x}}

\newcommand{\domainSpaceTime}{Q}
\newcommand{\boundarySpaceTime}{P}
\newcommand{\boundarySpaceTimeDirichlet}{\boundarySpaceTime_{\text{D}}}
\newcommand{\boundarySpaceTimeNeumann}{\boundarySpaceTime_{\text{N}}}

\newcommand{\intDomainSpaceTime}[1]{\int_{\domainSpaceTime} #1 \ d\domainSpaceTime}
\newcommand{\intElementSpaceTime}[1]{\int_{\domainSpaceTime^e}#1 \ d\domainSpaceTime}

\newcommand{\intBoundarySpaceTimeNeumann}[1]{\int_{\boundarySpaceTimeNeumann} #1 \ d\boundarySpaceTime}

\newcommand{\sobolevSpace}[1][]
{
    \ifthenelse{\equal{#1}{}}
        {\mathcal{H}^{1}}
        {\mathcal{H}^{1}\lp#1\rp}
}

\newcommand{\nBasisFOM}{N{^{h}}}

\newcommand{\bilinearFormNeumannSpaceTimeParam}[3]{\left(#1,#2;#3\right)_{\boundarySpaceTimeNeumann}}

\newcommand{\density}{\rho}
\newcommand{\visc}{\eta}

\newcommand{\param}{\mu}
\newcommand{\paramVec}{\vek{\param}}

\newcommand{\vel}{\vek{u}}
\newcommand{\velDirichlet}{\vel_{\text{D}}}
\newcommand{\velX}{u}
\newcommand{\velY}{v}
\newcommand{\velZ}{w}

\newcommand{\trialSpaceVelocity}{\vek{\mathcal{S}}_{\vel}}

\newcommand{\trialSpaceVelocityDiscrete}{\trialSpaceVelocity^h}

\newcommand{\trialVelocity}{\vek{u}}
\newcommand{\trialVelocityParam}{\vek{u}\lp\paramVec\rp}
\newcommand{\trialVelocityDiscrete}{\trialVelocity^h}
\newcommand{\trialVelocityHom}{\vek{v}}
\newcommand{\trialVelocityHomDiscrete}{\trialVelocityHom^h}
\newcommand{\liftingVelocityDiscrete}{\vek{l}^h}
\newcommand{\weightingSpaceVelocity}{\vek{\mathcal{V}}_{\vel}}

\newcommand{\weightingSpaceVelocityDiscrete}{\weightingSpaceVelocity^h}

\newcommand{\weightVelocity}{\vek{w}}
\newcommand{\weightVelocityDiscrete}{\weightVelocity^h}
\newcommand{\sobolevSpaceVector}[1][]
{
    \ifthenelse{\equal{#1}{}}
        {\vek{\mathcal{H}}^1}
        {\vek{\mathcal{H}}^1\lp#1\rp}
}
\newcommand{\nBasisVelocityFOM}{N^{h}_{\vel}}
\newcommand{\coeffVelocityFOM}[1]{v_{#1}}
\newcommand{\basisFunctionVelocityFOM}[1]{\vek{N}^{\vel}_{#1}}
\newcommand{\p}{p}

\newcommand{\trialSpacePressure}{\mathcal{S}_{\p}}

\newcommand{\trialSpacePressureDiscrete}{\trialSpacePressure^h}

\newcommand{\trialPressure}{p}
\newcommand{\trialPressureParam}{p\lp\paramVec\rp}
\newcommand{\trialPressureDiscrete}{\trialPressure^h}
\newcommand{\weightingSpacePressure}{\trialSpacePressure}

\newcommand{\weightingSpacePressureDiscrete}{\weightingSpacePressure^h}

\newcommand{\weightPressure}{q}
\newcommand{\weightPressureDiscrete}{\weightPressure^h}
\newcommand{\nBasisPressureFOM}{N^{h}_{\p}}
\newcommand{\coeffPressureFOM}[1]{p_{#1}}
\newcommand{\basisFunctionPressureFOM}[1]{N^{\p}_{#1}}
\newcommand{\nval}{h}
\newcommand{\bdf}{\vek{f}}
\newcommand{\cstressSymbol}{\mat{\sigma}}

\newcommand{\cstressParam}[2]{\cstressSymbol\lp #1,#2 ; \paramVec\rp}

\newcommand{\tractionNeumann}{\vek{h}}
\newcommand{\rosSymbol}{\mat{\epsilon}}
\newcommand{\ros}[1]{\rosSymbol\lp #1 \rp}
\newcommand{\shearrate}{\dot{\gamma}}

\newcommand{\bilinearFormPressure}[2]{b\lp #1, #2 \rp}
\newcommand{\bilinearFormPressureParam}[3]{b\lp #1, #2; #3 \rp}

\newcommand{\tMom}{\tau_{\textnormal{MOM}}\lp\trialVelocityDiscrete, \paramVec\rp}
\newcommand{\tMomPlain}{\tau_{\textnormal{MOM}}}

\newcommand{\bilinearFormStabMomStokesTemporal}[3]
{
s_{\text{MOM}}\lp#1,#2,#3\rp
}
\newcommand{\superscriptROM}{N}
\newcommand{\trialVelocityReduced}{\trialVelocity^{\superscriptROM}}
\newcommand{\trialVelocityHomDiscreteParam}{\trialVelocityHomDiscrete\lp\paramVec\rp}
\newcommand{\trialVelocityHomReduced}{\trialVelocityHom^{\superscriptROM}}

\newcommand{\trialPressureDiscreteParam}{\trialPressureDiscrete\lp\paramVec\rp}
\newcommand{\trialPressureReduced}{\trialPressure^{\superscriptROM}}

\newcommand{\bilinearFormViscousStressParam}[4]{a\lp #1 ,#2;#3, #4 \rp}

\newcommand{\matrixViscousStressSymbol}{A}
\newcommand{\matrixViscousStress}{\mat{\matrixViscousStressSymbol}}
\newcommand{\matrixViscousStressParam}{\matrixViscousStress\lp\trialVelocityDiscrete,\paramVec\rp}
\newcommand{\matrixPressureSymbol}{B}
\newcommand{\matrixPressure}{\mat{\matrixPressureSymbol}}
\newcommand{\matrixPressureTrans}{\matrixPressure^{\text{T}}}
\newcommand{\matrixStabMomStokesSymbol}{S}
\newcommand{\matrixStabMomStokes}{\mat{\matrixStabMomStokesSymbol}}
\newcommand{\matrixStabMomStokesParam}{\matrixStabMomStokes\lp\trialVelocityDiscrete,\paramVec\rp}
\newcommand{\vectorRHSVelocitySymbol}{F}
\newcommand{\vectorRHSVelocity}{\vek{\vectorRHSVelocitySymbol}}

\newcommand{\vectorRHSVelocityNonlinearSymbol}{L}
\newcommand{\vectorRHSVelocityNonlinear}{\vek{\vectorRHSVelocityNonlinearSymbol}}
\newcommand{\vectorRHSVelocityNonlinearParam}{\vek{\vectorRHSVelocityNonlinearSymbol}\lp\trialVelocityDiscrete, \paramVec\rp}
\newcommand{\vectorRHSPressureSymbol}{G}
\newcommand{\vectorRHSPressure}{\vek{\vectorRHSPressureSymbol}}
\newcommand{\velocityDOFVectorSymbol}{U}
\newcommand{\velocityDOFVector}{\vek{\velocityDOFVectorSymbol}}
\newcommand{\velocityDOFVectorParam}{\vek{\velocityDOFVectorSymbol}\lp\paramVec\rp}
\newcommand{\pressureDOFVectorSymbol}{P}
\newcommand{\pressureDOFVector}{\vek{\pressureDOFVectorSymbol}}
\newcommand{\pressureDOFVectorParam}{\vek{\pressureDOFVectorSymbol}\lp\paramVec\rp}

\newcommand{\basisFunctionMatrix}{\mat{Z}}
\newcommand{\nBasisROM}{N}
\newcommand{\nBasisVelocityROM}{N_{\trialVelocity}}
\newcommand{\nBasisVelocityHomROM}{N_{\trialVelocityHom}}
\newcommand{\weightingSpaceVelocityReduced}{\weightingSpaceVelocity^{\superscriptROM}}
\newcommand{\trialSpaceVelocityReduced}{\trialSpaceVelocity^{\superscriptROM}}

\newcommand{\basisFunctionMatrixVelocity}{\basisFunctionMatrix_{\trialVelocity}}
\newcommand{\basisFunctionMatrixVelocityTrans}{\basisFunctionMatrixVelocity^{\text{T}}}

\newcommand{\functionEIM}{g_{\text{EIM}}}
\newcommand{\functionEIMParam}[1]
{
    \ifthenelse{\equal{#1}{}}
    {\functionEIM^{\paramVec}}
    {\functionEIM^{\paramVec_{#1}}}
}

\newcommand{\coeffEIMSymbol}{c}
\newcommand{\coeffEIM}[1]{\coeffEIMSymbol_{#1}}

\newcommand{\basisFunctionEIM}[2]
{
    \ifthenelse{\equal{#2}{}}
    {h_{#1}}
    {h_{#1}\lp#2\rp}
}

\newcommand{\nBasisPressureROM}{N_{\trialPressure}}
\newcommand{\weightingSpacePressureReduced}{\weightingSpacePressure^{\superscriptROM}}

\newcommand{\basisFunctionMatrixPressure}{\basisFunctionMatrix_{\trialPressure}}
\newcommand{\basisFunctionMatrixPressureTrans}{\basisFunctionMatrixPressure^{\text{T}}}

\newcommand{\matrixViscousStressReduced}{\matrixViscousStress_{\superscriptROM}}
\newcommand{\matrixViscousStressReducedParam}{\matrixViscousStressReduced\lp\trialVelocityDiscrete,\paramVec\rp}
\newcommand{\matrixPressureReduced}{\mat{\matrixPressureSymbol}_{\superscriptROM}}
\newcommand{\matrixPressureReducedTrans}{\matrixPressureReduced^{\text{T}}}
\newcommand{\matrixStabMomStokesReduced}{\matrixStabMomStokes_{\superscriptROM}}
\newcommand{\matrixStabMomStokesReducedParam}{\matrixStabMomStokesReduced\lp\trialVelocityDiscrete,\paramVec\rp}
\newcommand{\vectorRHSVelocityReduced}{\vek{\vectorRHSVelocitySymbol}_{\superscriptROM}}

\newcommand{\vectorRHSVelocityNonlinearReduced}{\vek{\vectorRHSVelocityNonlinearSymbol}_{\superscriptROM}}
\newcommand{\vectorRHSVelocityNonlinearReducedParam}{\vek{\vectorRHSVelocityNonlinearSymbol}_{\superscriptROM}\lp\trialVelocityDiscrete,\paramVec\rp}
\newcommand{\vectorRHSPressureReduced}{\vek{\vectorRHSPressureSymbol}_{\superscriptROM}}
\newcommand{\velocityDOFVectorReduced}{\vek{\velocityDOFVectorSymbol}_{\superscriptROM}}
\newcommand{\velocityDOFVectorReducedParam}{\vek{\velocityDOFVectorSymbol}_{\superscriptROM}\lp\paramVec\rp}
\newcommand{\pressureDOFVectorReduced}{\vek{\pressureDOFVectorSymbol}_{\superscriptROM}}
\newcommand{\pressureDOFVectorReducedParam}{\vek{\pressureDOFVectorSymbol}_{\superscriptROM}\lp\paramVec\rp}

\newcommand{\matrixTemporalSymbol}{E}
\newcommand{\matrixTemporal}{\mat{\matrixTemporalSymbol}}
\newcommand{\matrixTemporalReduced}{\matrixTemporal_{\superscriptROM}}

\newcommand{\vectorRHSTemporalSymbol}{H}
\newcommand{\vectorRHSTemporal}{\vek{\vectorRHSTemporalSymbol}}
\newcommand{\vectorRHSTemporalReduced}{\vectorRHSTemporal_{\superscriptROM}}

\newcommand{\bilinearFormStabMomStokesTemporalParam}[5]
{s_{\text{MOM}}\lp#1,#2,#3;#4,#5\rp
}
\newcommand{\matrixStabMomStokesTemporalSymbol}{C}
\newcommand{\matrixStabMomStokesTemporal}{\mat{\matrixStabMomStokesTemporalSymbol}}
\newcommand{\matrixStabMomStokesTemporalParam}{\mat{\matrixStabMomStokesTemporalSymbol}\lp\trialVelocityDiscrete, \paramVec\rp}
\newcommand{\matrixStabMomStokesTemporalReduced}{\matrixStabMomStokesTemporal_{\superscriptROM}}
\newcommand{\matrixStabMomStokesTemporalReducedParam}{\matrixStabMomStokesTemporal_{\superscriptROM}\lp\trialVelocityDiscrete, \paramVec\rp}

\newcommand{\vectorRHSStabMomStokesTemporalSymbol}{D}
\newcommand{\vectorRHSStabMomStokesTemporal}{\vek{\vectorRHSStabMomStokesTemporalSymbol}}
\newcommand{\vectorRHSStabMomStokesTemporalParam}{\vek{\vectorRHSStabMomStokesTemporalSymbol}\lp\trialVelocityDiscrete, \paramVec\rp}
\newcommand{\vectorRHSStabMomStokesTemporalReduced}{\vectorRHSStabMomStokesTemporal_{\superscriptROM}}
\newcommand{\vectorRHSStabMomStokesTemporalReducedParam}{\vectorRHSStabMomStokesTemporal_{\superscriptROM}\lp\trialVelocityDiscrete, \paramVec\rp}

\newcommand{\nTrain}{N_{\text{train}}}

\newcommand{\nTest}{N_{\text{test}}}

\newcommand{\relErrorVelocity}{\varepsilon_{\trialVelocity}}
\newcommand{\relErrorPressure}{\varepsilon_{\trialPressure}}

\begin{document}

\maketitle


\section*{Abstract}
In the context of simulation-based methods, multiple challenges arise, two of which are considered in this work. As a first challenge, problems including time-dependent phenomena with complex domain deformations, potentially even with changes in the domain topology, need to be tackled appropriately. The second challenge arises when computational resources and the time for evaluating the model become critical in so-called many query scenarios for parametric problems. For example, these problems occur in optimization, \acrfull{uq}, or automatic control and using highly resolved \acrfullpl{fom} may become impractical.
To address both types of complexity, we present a novel projection-based \acrfull{mor} approach for deforming domain problems that takes advantage of the time-continuous space-time formulation. We apply it to two examples that are relevant for engineering or biomedical applications and conduct an error and performance analysis.
In both cases, we are able to drastically reduce the computational expense for a model evaluation and, at the same time, to maintain an adequate accuracy level. 
All in all, this work indicates the effectiveness of the presented \acrshort{mor} approach for deforming domain problems taking advantage of a time-continuous space-time setting.

\section{Introduction}
\label{sec:introdcution}
The use of simulation-based methods is nowadays widespread in scientific and engineering applications. Commonly, these methods are intended to serve purposes like enhanced insight and prediction capabilities with respect to the processes under investigation. In this way, they can be useful tools in the context of product design or optimization procedures. Furthermore, they may also help to quantify or enhance the reliability of a given process by techniques such as \gls{uq}. Lastly, they further allow to be integrated into ongoing operations, e.g., as digital twins or to set up optimal control.
\bigskip\par
Arising in the context of simulation-based methods, we will focus on two specific types of complexity in this work. First, we consider applications that show a \textit{transient} behavior involving a \textit{deforming domain} and, potentially, even \textit{topology changes}. These aspects need to be treated appropriately by the computational approach to cover all the relevant effects.
In the following, we restrict ourselves to methods that make use of a computational mesh although there exists a variety of further simulation methods. In the case discussed here, the deformation of the domain requires methods that can appropriately handle both the computational mesh and the unsteady solution field. Existing methods will be discussed in more detail in \Cref{subsec:methodologicalBackgroundFOM}.
The second kind of complexity is related to the \textit{expense} --- in terms of computational resources and time --- needed for an evaluation of the computational model. Especially the design and development phase of a process or a product may entail the assessment of various operating points or configurations, the optimization of process settings or component features, as well as uncertainty regarding involved quantities. In all cases, one ends up with a so-called many query scenario in which a great number of model evaluations needs to be performed. Furthermore, the integration into an automatic control environment may demand fast feedback from the simulation model. Typically, one can formulate the involved problems in a parametric manner where each of the model evaluations is characterized by a certain sample of parameter values.
In these situations, employing the original model, which is often referred to as \textit{\gls{fom}}, for each sample may easily exceed available resources or required feedback times. Here, the application of \gls{mor} can provide a remedy. Based on the original model, a \textit{\gls{rom}} with decreased computational complexity is constructed, while keeping its accuracy in the desired range. Common \gls{mor} techniques will be presented later in \Cref{subsec:methodologicalBackgroundROM}.
\bigskip\par
The work presented here addresses both types of complexity described above and entails a \gls{mor} approach for problems that are characterized by a transient solution field in a deforming domain with possibly changing topology. Regarding the question of deforming domain problems, we rely on the \textit{time-continuous space-time} setting. To significantly decrease the computational demands and/or response times of the underlying model in a next step, we apply a \textit{projection-based} \gls{mor} technique for which we make use of \gls{pod}. Key of the proposed approach is the particular combination of this \gls{mor} technique with the choice of a time-continuous space-time formulation. This connection allows us to construct a corresponding \gls{rom} for the aforementioned class of problems in a straightforward manner. In this way, the benefits of \gls{mor} are made easily accessible even when dealing with deforming domain problems that involve an unsteady solution and, if necessary, changes in the spatial topology. Alternative \gls{mor} approaches that exist and are applicable to transient or deformation-driven scenarios will also be discussed later, i.e., in \Cref{subsec:methodologicalBackgroundROM}.
\par
As examples for the type of problems in focus, we consider the simulation of two specific transient fluid flow scenarios like they may appear in engineering or biomedical applications. The domain deformation in these examples results from a moving valve plug or the narrowing of flexible artery walls. Furthermore, the parameterization is induced through a variation of the material properties of the fluid or of the boundary conditions.
\bigskip\par
The remainder of the article is structured as follows: in \Cref{sec:parametricFormulation}, we derive the parametric formulation for fluid flow problems in deforming domains using the space-time perspective. Next, \Cref{sec:fullOrderModel} contains the description of the respective \gls{fom}, which will be based on the \gls{fem}. Subsequently, the construction of the corresponding \gls{rom} using \gls{pod} and projection is outlined in \Cref{sec:reducedOrderModel}. Results for the two fluid flow test cases covering three- and four-dimensional space-time domains are presented in \Cref{sec:numericalResults} to illustrate the aptitude of the approach introduced in this work. Finally, our findings are summarized and discussed in \Cref{sec:conclusion}. 

\section{Parametric Problem Formulation for Fluid Flow in Deforming Domains}
\label{sec:parametricFormulation}
To lay the foundation for the subsequent \gls{mor} approach, we begin by deriving the \textit{parametric formulation} of transient flow problems defined in deforming domains. As mentioned before, parametric problems, e.g., may occur in the context of \gls{uq} or automatic control. Specifically, we consider problems that can be parametric in the \textit{material properties} involved or in their \textit{boundary conditions}. The variations in the material illustrate potential uncertainties in the process under investigation, which may be analyzed using \gls{uq}. Apart from that, adjustable boundary conditions prescribing, e.g., the inflow velocity, showcase the type of use in a context such as automatic control. 
As mentioned earlier, we focus on unsteady processes that take place in time-varying domains, which potentially undergo topology changes, too.

\bigskip\par
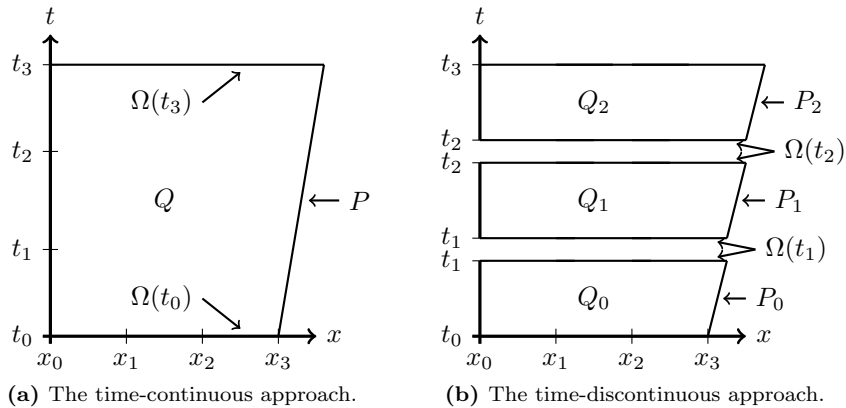
\begin{figure}
    \centering
    \subcaptionbox{The time-continuous approach.\label{fig:cst}}{
    \newcommand*\cols{2}
\newcommand*\rows{2}
\begin{tikzpicture}[
    scale=1,
    axis/.style={very thick, ->},
    important line/.style={thick},
    every node/.style={color=black}
    ]
    \draw[axis] (-0.1,0)  -- (3.5,0) node(xline)[right]{$x$};
    \foreach \x in {0,1,2,3}
    \draw (\x cm,0.1) -- (\x cm,-0.1) node[anchor=north] {$x_{\x}$};
    \draw[axis] (0,-0.1) -- (0,4) node(yline)[above]{$t$};
    \draw (0.1,0cm) -- (-0.1, 0cm) node[anchor=east] {$t_{0}$};
    \draw (0.1,1.15cm) -- (-0.1, 1.15cm) node[anchor=east] {$t_{1}$};
    \draw (0.1,2.45cm) -- (-0.1, 2.45cm) node[anchor=east] {$t_{2}$};
    \draw (0.1,3.6cm) -- (-0.1, 3.6cm) node[anchor=east] {$t_{3}$};
    \draw[important line] (3,0) coordinate (A)  -- (3.6,3.6) coordinate (B);
    \draw[important line] (0,3.6) coordinate (A)  -- (3.6,3.6) coordinate (B);
    \node (Q) at (1.5,1.8) {$Q$};
    \draw[thick,<-] (3.4,1.8cm) -- (3.8, 1.8) node[anchor=west] {$P$};
    \draw[thick,<-] (2.5,0.1cm) -- (2.0, 0.5) node[anchor=east] {$\Omega(t_0)$};
    \draw[thick,<-] (2.5,3.5cm) -- (2.0, 3.1) node[anchor=east] {$\Omega(t_3)$};
\end{tikzpicture}
    }
    \subcaptionbox{The time-discontinuous approach.\label{fig:dst}}{
    \newcommand*\cols{2}
\newcommand*\rows{3}
\begin{tikzpicture}[
    scale=1,
    axis/.style={very thick, ->},
    important line/.style={thick},
    newEdge/.style={thick, ->, >=circle},
    every node/.style={color=black}
    ]
    \draw[axis] (-0.1,0)  -- (3.5,0) node(xline)[right]{$x$};
    \foreach \x in {0,1,2,3}
    \draw (\x cm,0.1) -- (\x cm,-0.1) node[anchor=north] {$x_{\x}$};
    \draw[very thick] (0,-0.1) -- (0,1);
       \draw[very thick] (0,1.3) -- (0,2.3);
    \draw[axis] (0,2.6) -- (0,4) node(yline)[above]{$t$};
    \draw (0.1,0cm) -- (-0.1, 0cm) node[anchor=east] {$t_{0}$};
    \draw (0.1,1cm) -- (-0.1, 1cm) node[anchor=east] {$t_{1}$};
    \draw (0.1,1.3cm) -- (-0.1, 1.3cm) node[anchor=east] {$t_{1}$};
    \draw (0.1,2.3cm) -- (-0.1, 2.3cm) node[anchor=east] {$t_{2}$};
        \draw (0.1,2.6cm) -- (-0.1, 2.6cm) node[anchor=east] {$t_{2}$};
            \draw (0.1,3.6cm) -- (-0.1, 3.6cm) node[anchor=east] {$t_{3}$};
    \foreach \i / \row in {0/0, 1/1.3, 2/2.6} {
     		\foreach \col in {0, 1, ...,\cols} {
   				 \draw[important line] (\col,\row) coordinate (A)  -- (\col+1+0.25*\i,\row) coordinate (B);
   				 \draw[important line] (\col,\row+1) coordinate (A)  -- (\col+1+0.25*\i+0.25,\row+1) coordinate (B);

   		 }		  
   	    \draw[important line] (\cols+1+0.25*\i,\row) coordinate (A)  -- (\cols+1+0.25*\i+0.25,\row+1) coordinate (B);   		  \node (Q) at (1.5,0.5+\row) {$Q_{\i}$};		  
        \draw[thick,<-] (3.2+0.25*\i,0.5+\row) -- (3.5+0.25*\i, 0.5+\row) node[anchor=west] {$P_{\i}$};
    }
    \draw[thick,<-] (3.0+0.25*0.5,1.25cm) -- (3.5+0.25*0.5, 1.15) node[anchor=west] {$\Omega(t_1)$};
    \draw[thick,<-] (3.0+0.25*0.5,1.05cm) -- (3.5+0.25*0.5, 1.15);
    \draw[thick,<-] (3.0+0.25*1.5,2.35cm) -- (3.5+0.25*1.5, 2.45) node[anchor=west] {$\Omega(t_2)$};
    \draw[thick,<-] (3.0+0.25*1.5,2.55cm) -- (3.5+0.25*1.5, 2.45);
\end{tikzpicture}
    }

    \caption{Definition of the computational domain(s) in the space-time approach for deforming domain problems.}
    \label{fig:parametricProblemComputationalDomain}
\end{figure}

For the description of the time-dependent solution field in a deforming domain, it is convenient to treat space and time coordinates, which are denoted by $\x$ and $t$, in the same way. To that end, the \textit{time-continuous space-time domain} $\domainSpaceTime$ is introduced, which results from the time-dependent spatial domain $\domainSpace(t)$ and the time interval $[0,T]$:
\begin{align}
\domainSpaceTime = \domainSpace(t) \times [0,T],    
\end{align}
such that $(\x,t) \in \domainSpaceTime$. The respective space-time boundary is denoted as $\boundarySpaceTime$ where portions on which conditions of Dirichlet or Neumann type are prescribed are referred to as $\boundarySpaceTimeDirichlet$ and $\boundarySpaceTimeNeumann$, respectively.  Furthermore, initial conditions can be imposed over $\domainSpace(t=0)=\domainSpace(t_0)$. A sketch can be found in \Cref{fig:cst}. In contrast, in a \textit{time-discontinuous space-time} approach, one considers so-called time slabs $\domainSpaceTime_m$ with boundaries $\boundarySpaceTime_m$(see \Cref{fig:dst}). This results in a kind of time-stepping scheme with a modified computational domain in each time step. In that sense, the time-discontinuous approach is comparable to a semi-discrete description. A systematic comparison of time-continuous and time-discontinuous space-time formulations for advection-diffusion problems can be found in ~\cite{Danwitz2022}. However, as we would like to stress again, the time-continuous description is essential to directly apply established \gls{mor} techniques for deforming domain problems later.
\bigskip
\par
To be able to eventually demonstrate the presented approach, we consider a specific example of fluid flow equations here. In particular, we make use of the \textit{Stokes equations}, which can be used to model a variety of fluid flow of low Reynolds number, for example.
Based on the conservation of mass and momentum, the Stokes equations provide statements for the, in this case parametric, \textit{velocity} and \textit{pressure} fields of the fluid, which are referred to as $\vel(\x,t;\paramVec)$ and $\p(\x,t;\paramVec)$, respectively. Here, $\paramVec$ denotes the \textit{parameter vector} collecting parameter values to vary either material properties or boundary conditions.
Note that we will drop the arguments of $\vel$ and $\p$ for the sake of notation in the following. Nevertheless, dependencies on $\paramVec$ will be stressed if necessary.
\par
The resulting \gls{bvp} in the space-time setting reads as follows:
\begin{align}
        \dv{\vel(\paramVec)} &= 0 &&\text{ in } \domainSpaceTime,
        \\
        \density \lp \dd{\vel(\paramVec)}{t} - \bdf \rp - \dv{\cstressParam{\vel}{\p}} &= \vek{0} &&\text{ in } \domainSpaceTime,
        \\
        \vel(\paramVec) &= \velDirichlet(\paramVec) &&\text{ on } \boundarySpaceTimeDirichlet,
        \\
        \cstressParam{\vel}{\p}\cdot \vek{n} &= \tractionNeumann(\paramVec) &&\text{ on } \boundarySpaceTimeNeumann,
        \\
        \vel(\x,0;\paramVec) &= \vel_0(\paramVec) &&\text{ in } \domainSpace(t_0).
\end{align}
Here, $\density$ is the fluid \textit{density}, $\bdf$ is a \textit{body force}, while $\velDirichlet(\paramVec)$, $\vel_0(\paramVec)$, and $\tractionNeumann(\paramVec)$ are the parameter-dependent prescribed velocity values and a \textit{traction vector}, respectively. 
Moreover, we consider parametric material properties through the fluid \textit{viscosity} $\visc(\paramVec)$ that is included in the \textit{Cauchy stress tensor} given by
\begin{align}
    \cstressParam{\vel}{\p} = - \p\lp\paramVec\rp \mat{I} + 2 \visc\lp\paramVec\rp \ros{\vel\lp\paramVec\rp},
\end{align}
where
$
    \ros{\vel} = \frac{1}{2} ( \vek{\nabla} \vel + \vek{\nabla} \vel^T )
$
is the \textit{rate-of-strain tensor}.
In the following, we will consider fluid materials whose viscous properties cannot be appropriately described by a Newtonian model. Examples for such non-Newtonian fluids are plastics melt or blood. 
To account for shear-thinning effects of the fluid, we use the Carreau-Yasuda model \cite{Carreau1972, Yasuda1981, Bird1987} for the viscosity, which depends on the shear rate 
\begin{align}
\shearrate(\vel; \paramVec) = \sqrt{2 \ros{\vel(\paramVec)} : \ros{\vel(\paramVec)}}.   
\end{align}
To stress the dependency on the velocity field, we from now on write $\visc(\vel; \paramVec)$ given by
\begin{align}
        \visc(\vel; \paramVec)
        &= 
        \visc_{\infty} 
        + \left( \visc_0 - \visc_{\infty} \right) \left[ 1 + \left( \lambda \shearrate \right)^a \right]^{\frac{n-1}{a}},
    \label{eq:viscosityModel}
\end{align}
with zero-shear-rate viscosity $\visc_0$, infinite-shear-rate viscosity $\visc_{\infty}$, characteristic time $\lambda$, power-law index $n$, and a dimensionless parameter $a$ describing the transition region between the zero-shear-rate viscosity and the power-law regions.
In addition to the indirect parameter dependency of the viscosity through the parametric velocity field, we also take into account potential fluctuations in the model parameters. This is motivated by the fact that these model parameters can only be determined through a regression-based approach using measurement data associated with uncertainties. 
\bigskip\par
Since the \gls{fem} will be used later to construct the \gls{fom} in \Cref{sec:fullOrderModel}, which in turn will serve as the basis for the \gls{rom} presented in \Cref{sec:reducedOrderModel}, it is worthwhile to state the \textit{weak formulation} of the problem here. We will assume that appropriate trial and weighting spaces for the velocity and the pressure field in the space-time domain $\domainSpaceTime$ are given. For the velocity, the trial and weighting spaces are denoted as $\trialSpaceVelocity(\domainSpaceTime)$ and $\weightingSpaceVelocity(\domainSpaceTime)$, respectively. For the pressure, the trial space is also the weighting space and they are referred to as $\trialSpacePressure(\domainSpaceTime)$. The resulting weak formulation of the parametric problem defined in the space-time domain $\domainSpaceTime$ reads:
\bigskip\par\noindent
\textit{Find }$\lp \trialVelocityParam, \trialPressureParam \rp 
\in 
\trialSpaceVelocity\lp\domainSpaceTime\rp 
\times 
\trialSpacePressure\lp\domainSpaceTime\rp $\textit{, such that }$\forall \lp \weightVelocity, \weightPressure \rp 
\in 
\weightingSpaceVelocity\lp\domainSpaceTime\rp 
\times 
\weightingSpacePressure\lp\domainSpaceTime\rp $\textit{:}
\begin{align}
    &
    \intDomainSpaceTime
    {
        \weightVelocity \cdot \density 
        \lp
            \dd{\trialVelocityParam}{t} - \bdf
        \rp
    }
    -
    \intDomainSpaceTime
    {
        \trialPressureParam \lp \dv{\weightVelocity}\rp
    }
    \nonumber
    \\
    +
    &
    \intDomainSpaceTime
    {
        \ros{\weightVelocity} : 2 \visc\lp\vel,\paramVec\rp \ros{\trialVelocity}
    }
    +
    \intDomainSpaceTime
    {
        \weightPressure \lp \dv{\trialVelocityParam}\rp
    }
    =
    \intBoundarySpaceTimeNeumann
    {
        \weightVelocity \cdot \tractionNeumann\lp\paramVec\rp
    }
    .
    \label{eq:navierStokesWeakExtended}
\end{align}
To simplify notation, the following short forms will be used in the remainder of this work:
\begin{align*}
    \begin{aligned}[t]
        \ddt{\trialVelocity}\lp\paramVec\rp 
        &= 
        \dd{\trialVelocityParam}{t},
        \\
        \bilinearFormPressureParam{\trialVelocity}{\weightPressure}{\paramVec}_{\domainSpaceTime}  
        &=     
        \intDomainSpaceTime
        {
            \weightPressure \lp \dv{\trialVelocityParam}\rp
        },
        \\
        \lp \weightVelocity, \bdf \rp_{\domainSpaceTime} 
        &=
        \intDomainSpaceTime
        {
           \density \weightVelocity \cdot \bdf
        },
    \end{aligned}
    \qquad
    \begin{aligned}[t]
        \lp \weightVelocity, \ddt{\trialVelocity}; \paramVec \rp_{\domainSpaceTime} 
        &= 
        \intDomainSpaceTime
        {
           \density \weightVelocity \cdot \ddt{\trialVelocity}\lp\paramVec\rp
        },
        \\
        \bilinearFormViscousStressParam{\weightVelocity}{\trialVelocity}{\trialVelocityHom}{\paramVec}_{\domainSpaceTime}  
        &=     
        \intDomainSpaceTime
        {
            \ros{\weightVelocity} : 2 \visc\lp \trialVelocityHom,\paramVec\rp \ros{\trialVelocity}
        },
        \\
        \lp \weightVelocity, \tractionNeumann; \paramVec \rp_{\boundarySpaceTimeNeumann}
        &=     
        \intBoundarySpaceTimeNeumann
        {
            \weightVelocity \cdot \tractionNeumann\lp\paramVec\rp
        }
        .
    \end{aligned}
\end{align*}
In this notation, \Cref{eq:navierStokesWeakExtended} can be stated as
\begin{align}
    \lp \weightVelocity, \ddt{\trialVelocity}; \paramVec \rp_{\domainSpaceTime} 
    - \bilinearFormPressureParam
    {\weightVelocity}
    {\trialPressure}{\paramVec}_{\domainSpaceTime} 
    + \bilinearFormViscousStressParam
    {\weightVelocity}
    {\trialVelocity}
    {\vel}
    {\paramVec}_{\domainSpaceTime}
    +\bilinearFormPressureParam
    {\trialVelocity}
    {\weightPressure}{\paramVec}_{\domainSpaceTime}
    =
    \lp \weightVelocity, \bdf \rp_{\domainSpaceTime}  
    +\lp \weightVelocity, \tractionNeumann; \paramVec \rp_{\boundarySpaceTimeNeumann}
    .
\end{align}
\section{Full-Order Model for Fluid Flow in Deforming Domains}
\label{sec:fullOrderModel}
In this section, we will derive the \gls{fom} which results from applying the \gls{fem} to the parametric problem introduced beforehand. In particular, this model will be capable of simulating transient problems in deforming domains. To that end, \Cref{subsec:methodologicalBackgroundFOM}  will first contain a brief outline of existing techniques that address deforming domain problems in this context. Subsequently, we present the space-time \gls{fem} used in this work in \Cref{subsec:discreteFormulationFOM}.

\subsection{Methodological Background for Deforming Domains in the Full-Order Model}
\label{subsec:methodologicalBackgroundFOM}
In the following, we will roughly outline existing approaches for deforming domain problems, i.e., problems including moving boundaries or internal interfaces, without claiming completeness. First, one can make a division into interface capturing and interface tracking methods \cite{Elgeti2016}. Well-known examples for interface capturing methods, which employ an implicit description of domain boundaries or interfaces on a background mesh, are the level-set method \cite{Osher1988,Chang1996} or the volume-of-fluid method \cite{Hirt1981}. In contrast, interface tracking methods are based on boundary-conforming meshes. Thus, an update procedure that adapts the computational mesh according to the deforming domain is needed.
\par
There exist a vast number of strategies, ranging from global remeshing to elaborate methods specifically designed for certain applications or types of deformation. As examples for connectivity-preserving methods that update nodal coordinates, one can mention \gls{pde}-based methods like the \gls{emum} \cite{Tezduyar1992}, spring-based methods \cite{Batina1990} or techniques using radial basis functions \cite{DeBoer2007}. Further methods extend these strategies by a subsequent mesh optimization, e.g., through edge swapping or vertex smoothing operations \cite{Alauzet2013,Wang2015a}. Moreover, specialized mesh update methods use algebraic operations to control the mesh evolution for a-priori known boundary deformations~\cite{Helmig2019, Hinz2020}. Furthermore, local remeshing strategies limit the cost of maintaining a high quality boundary-conforming mesh under large deformations~\cite{Behr1999,Behr2001,Behr2003}. Broadening the scope of these methods, parts of the mesh can be activated and deactivated~\cite{Key2018} and topology changes of the computational domain can be handled elegantly~\cite{Gonzalez2023}.
\par
As an alternative, one can also follow weak domain coupling strategies to account for the moving boundary or interface by using composite grids and introducing additional conditions on the solution field. Examples are the Chimera method \cite{Steger1983, Benek1986, Steger1987} or sliding interface approaches \cite{Bazilevs2008, Takizawa2015, Helmig2020}. It is also worthwhile to mention the immersed boundary method \cite{Peskin1972}, where simulations are always performed on a cartesian background grid.
\par
Deforming domain problems are inherently transient and their solution requires a form of time discretization. Common approaches include time-stepping methods which separate spatial and temporal discretization and space-time methods which apply a combined discretization to the space-time domain. Time stepping schemes, e.g., the generalized $\alpha$-method~\cite{Chung1993}, require an arbitrary Lagrangian–Eulerian (ALE) formulation for moving-domain simulations~\cite{Hughes1981, Donea2004}, while the space-time formulation directly accounts for the (spatial) mesh deformation~\cite{Tezduyar1992a, Tezduyar1992b}.
\par
If the movement is known during mesh generation --- typically this means before the simulation start --- the time-continuous space-time approach allows to incorporate complex deformations of the spatial domain in a boundary-conforming space-time mesh. Even topology changes can be included as shown in finite volume and finite element simulations~\cite{Rendall2012, Danwitz2019}. For two-dimensional problems, standard mesh generation tools can be used to construct the three-dimensional space-time mesh. For three-dimensional problems, more sophisticated mesh generation and adaptation techniques are required. The common approach to generate an unstructured four-dimensional mesh is based on the extrusion of a tetrahedral mesh followed by a subdivision of the prismatic elements into pentatopes (4-simplices). The subdivision can either be achieved with an element-wise Delaunay triangulation~\cite{Behr2008} or with a predefined decomposition which requires a consistently numbered tetrahedral mesh~\cite{Karabelas2019}. Further techniques enable refinement and anisotropic adaptation of four-dimensional meshes~\cite{Neumueller2011, Wang2015a, Caplan2020}.  Four-dimensional meshes with complex deformations and topology changes of the three-dimensional spatial domain can be obtained with an elastic mesh update following extrusion based pentatope mesh generation~\cite{Danwitz2021}. Please, note that the additional effort for the space-time mesh generation completely replaces the special treatment of a deforming domain problem during the simulation. Nevertheless, our boundary-conforming space-time mesh approach is limited to 4D geometries that can be obtained by extrusion of a 3D geometry and a subsequent elastic deformation. Generating meshes for general 4D geometries of engineering scale is -- to the best of our knowledge -- an open research problem.
\subsection{Discrete Formulation for the Full-Order Model}
\label{subsec:discreteFormulationFOM}
Next, we will derive the \gls{fom}, which will be based on the \gls{fem}. Thus, we introduce corresponding finite dimensional subspaces for the trial and weighting spaces introduced in \Cref{sec:parametricFormulation}. Let $\trialSpaceVelocityDiscrete(\domainSpaceTime)$ and $\weightingSpaceVelocityDiscrete(\domainSpaceTime)$ be the finite dimensional subspaces of $\trialSpaceVelocity(\domainSpaceTime)$ and $\weightingSpaceVelocity(\domainSpaceTime)$, respectively. Furthermore, $\trialSpacePressureDiscrete(\domainSpaceTime)$ is the finite dimensional subspace of $\trialSpacePressure(\domainSpaceTime)$. We  stick to the \textit{\gls{csst}} approach, i.e., the computational mesh will be composed of simplex elements filling the entire space-time domain $\domainSpaceTime$. Furthermore, we use first-order polynomials as shape functions $\basisFunctionVelocityFOM{i}$ and $\basisFunctionPressureFOM{i}$ for velocity and pressure, respectively.
\par
To handle (parametric) Dirichlet boundary conditions, we introduce a \textit{velocity lifting function} $\liftingVelocityDiscrete(\paramVec)$ such that the \textit{discrete velocity trial function} is given as
\begin{align}
    \trialVelocityDiscrete(\paramVec) = \trialVelocityHomDiscrete(\paramVec) + \liftingVelocityDiscrete(\paramVec),
\end{align}
with the homogeneous portion $\trialVelocityHomDiscrete(\paramVec)$ and $\liftingVelocityDiscrete(\paramVec)\lvert_{\boundarySpaceTimeDirichlet}=\velDirichlet(\paramVec)$.
Note that additional lifting functions may be used to separate parameter-dependent and independent portions of the Dirichlet boundary conditions.
For the sake of notation, however, we only consider the case of a single lifting function here.
Consequently, it holds that $\trialVelocityDiscrete \in \trialSpaceVelocityDiscrete(\domainSpaceTime)$ and $ \trialVelocityHomDiscrete \in \weightingSpaceVelocityDiscrete(\domainSpaceTime)$. Moreover, the \textit{discrete pressure trial function} is denoted as $\trialPressureDiscreteParam \in \trialSpacePressureDiscrete(\domainSpaceTime)$
\par
We apply the \gls{gls} stabilization technique \cite{Hughes1987,Hughes1989,Mittal1992,Shakib1991}, which adds a least-squares form of the residual within each element to the original variational formulation of the problem. In the \gls{csst} formulation, this will apply to space-time elements denoted by $\domainSpaceTime^e$.
\bigskip\par
Following the description, e.g., from \cite{Behr1994,Franca1992}, the resulting space-time Galerkin formulation including the additional stabilization terms reads:
\bigskip\par\noindent
\textit{Find }$\lp \trialVelocityHomDiscreteParam, \trialPressureDiscreteParam \rp 
\in 
\weightingSpaceVelocityDiscrete\lp\domainSpaceTime\rp
\times 
\trialSpacePressureDiscrete\lp\domainSpaceTime\rp$\textit{, such that }$\forall \lp \weightVelocityDiscrete, \weightPressureDiscrete \rp 
\in
\weightingSpaceVelocityDiscrete\lp\domainSpaceTime\rp 
\times 
\weightingSpacePressureDiscrete\lp\domainSpaceTime\rp$\textit{:}
\begin{align}
    \lp \weightVelocityDiscrete, \ddt{\trialVelocityHomDiscrete}; \paramVec \rp_{\domainSpaceTime} 
    - \bilinearFormPressureParam
    {\weightVelocityDiscrete}
    {\trialPressureDiscrete}{\paramVec}_{\domainSpaceTime} 
    + \bilinearFormViscousStressParam
    {\weightVelocityDiscrete}
    {\trialVelocityHomDiscrete}{\trialVelocityDiscrete}{\paramVec}_{\domainSpaceTime}
    +\bilinearFormPressureParam
    {\trialVelocityHomDiscrete}
    {\weightPressureDiscrete}{\paramVec}_{\domainSpaceTime}&
    \nonumber
    \\
    +\bilinearFormStabMomStokesTemporalParam
    {\weightPressureDiscrete}
    {\trialVelocityDiscrete}
    {\trialPressureDiscrete}
    {\trialVelocityDiscrete}
    {\paramVec}_{\domainSpaceTime}&
    \nonumber
    \\
    =
    \lp \weightVelocityDiscrete, \bdf \rp_{\domainSpaceTime}  
    +\lp \weightVelocityDiscrete, \tractionNeumann; \paramVec \rp_{\boundarySpaceTimeNeumann}
    -\lp \weightVelocityDiscrete, \ddt{\liftingVelocityDiscrete}; \paramVec \rp_{\domainSpaceTime} 
    -\bilinearFormViscousStressParam
    {\weightVelocityDiscrete}
    {\liftingVelocityDiscrete}{\trialVelocityDiscrete}{\paramVec}_{\domainSpaceTime}
    -\bilinearFormPressureParam
    {\liftingVelocityDiscrete}
    {\weightPressureDiscrete}{\paramVec}_{\domainSpaceTime}&
    ,
    \label{eq:navierStokesGalerkinStabilizedST}
\end{align}
with
\begin{align}
\bilinearFormStabMomStokesTemporalParam
    {\weightPressureDiscrete}
    {\trialVelocityDiscrete}
    {\trialPressureDiscrete}
    {\trialVelocityDiscrete}
    {\paramVec}
    _{\domainSpaceTime} 
    =
    s^{\trialVelocityHom}\lp\weightPressureDiscrete,\trialVelocityHomDiscrete;\trialVelocityDiscrete,\paramVec\rp_{\domainSpaceTime}
    +
    s^{\vek{l}}\lp\weightPressureDiscrete,\liftingVelocityDiscrete;\trialVelocityDiscrete,\paramVec\rp_{\domainSpaceTime}
    +
    s^{\trialPressure}\lp\weightPressureDiscrete,
    \trialPressureDiscrete;\trialVelocityDiscrete,\paramVec\rp_{\domainSpaceTime},
\end{align}
and
\begin{align}
    s^{\trialVelocityHom}
    \lp
    \weightPressureDiscrete,
    \trialVelocityHomDiscrete;
    \trialVelocityDiscrete,
    \paramVec
    \rp_{\domainSpaceTime}
    &=
    \sum_{e} \intElementSpaceTime{\tMom \frac{1}{\rho}
    \lp
    -\gr{\weightPressureDiscrete} 
    \rp
    \cdot 
    \lp
    \rho \ddt{\trialVelocityHomDiscrete}\lp\paramVec\rp
    \rp
    },
    \\
    s^{\vek{l}}
    \lp
    \weightPressureDiscrete,
    \liftingVelocityDiscrete;
    \trialVelocityDiscrete,
    \paramVec
    \rp_{\domainSpaceTime}
    &=
    \sum_{e} \intElementSpaceTime{\tMom \frac{1}{\rho}
    \lp
    -\gr{\weightPressureDiscrete} 
    \rp
    \cdot 
    \lp
    \rho \ddt{\liftingVelocityDiscrete}\lp\paramVec\rp
    \rp
    },
    \\
    s^{\trialPressure}
    \lp
    \weightPressureDiscrete,
    \trialPressureDiscrete;
    \trialVelocityDiscrete,
    \paramVec
    \rp_{\domainSpaceTime}
    &=
    \sum_{e} \intElementSpaceTime{\tMom \frac{1}{\rho}
    \gr{\weightPressureDiscrete} 
    \cdot 
    \gr{\trialPressureDiscreteParam} 
    }.
\end{align}
\par
The stabilization parameter $\tMom$ is chosen as presented in \cite{Pauli2016a}. Although the formulation in detail is not of great importance for this work, note that it depends both on the parametric velocity $\trialVelocityDiscrete(\paramVec)$ and the parametric viscosity $\visc(\trialVelocityDiscrete, \paramVec)$ in a non-linear way. Furthermore, the second-order derivatives of the velocity weighting and trial functions, which appear in the original formulation of the momentum stabilization, are omitted due to the first-order linear basis functions in use.
\bigskip\par
As a foundation for the description of the \gls{rom} in the following section, we present next the \textit{algebraic formulation} of the problem. The vectors of coefficients are denoted as $\velocityDOFVector \in \mathbb{R}^{\nBasisVelocityFOM}$ and $\pressureDOFVector\in \mathbb{R}^{\nBasisPressureFOM}$ for the homogeneous velocity $\trialVelocityHomDiscrete$ and pressure field $\trialPressureDiscrete$, respectively. Here, $\nBasisVelocityFOM$ and $\nBasisPressureFOM$ stand for the number of \glspl{dof} in the \gls{fom}. The solution can then be computed by solving the following non-linear system for $\velocityDOFVectorParam$ and $\pressureDOFVectorParam$:
\begin{align}
    \lb
    \begin{array}{cc}
        \matrixTemporal + \matrixViscousStressParam         &  -\matrixPressureTrans 
        \\
        \matrixPressure + \matrixStabMomStokesTemporalParam & \matrixStabMomStokesParam
    \end{array}
    \rb
    \lb
    \begin{array}{c}
        \velocityDOFVectorParam\\
        \pressureDOFVectorParam
    \end{array}
    \rb
    =
    \lb
    \begin{array}{c}
         \vectorRHSTemporal + \vectorRHSVelocity + \vectorRHSVelocityNonlinearParam\\
         \vectorRHSPressure + \vectorRHSStabMomStokesTemporalParam
    \end{array}
    \rb&
    ,
    \label{eq:algebraicSystemFOM}
\end{align}
where the \gls{lhs} matrices for $i,j = 1, \dots, \nBasisVelocityFOM $ and for $k,l = 1, \dots, \nBasisPressureFOM $ are given as
\begin{align*}
    \begin{aligned}[t]
        \matrixTemporal &= \lb \matrixTemporalSymbol_{i,j} \rb, 
        \text{ with } 
        \matrixTemporalSymbol_{i,j} 
        = 
        \lp \basisFunctionVelocityFOM{i}, \dd{\basisFunctionVelocityFOM{j}}{t} \rp_{\domainSpaceTime},
        \\
        \matrixPressure &= \lb\matrixPressureSymbol_{k,j}\rb,
        \text{ with }
        \matrixPressureSymbol_{k,j} 
        = 
        \bilinearFormPressure
        {\basisFunctionVelocityFOM{j}}
        {\basisFunctionPressureFOM{k}}_{\domainSpaceTime},
        \\
        \matrixStabMomStokes &= \lb \matrixStabMomStokesSymbol_{k,l}\rb,
        \text{ with }
        \matrixStabMomStokesSymbol_{k,l}
        = 
        s^{\trialPressure}
        \lp
        \basisFunctionPressureFOM{k},
        \basisFunctionPressureFOM{l};
        \trialVelocityDiscrete,
        \paramVec
        \rp_{\domainSpaceTime},
    \end{aligned}
    \quad
    \begin{aligned}[t]
        \matrixViscousStress &= \lb A_{i,j}\rb,
        \text{ with } 
        \matrixViscousStressSymbol_{i,j} 
        = 
        \bilinearFormViscousStressParam
        {\basisFunctionVelocityFOM{i}}
        {\basisFunctionVelocityFOM{j}}
        {\trialVelocityDiscrete}
        {\paramVec},    
        \\
        \matrixStabMomStokesTemporal &= \lb \matrixStabMomStokesTemporalSymbol_{k,j}\rb,
        \text{ with }
        \matrixStabMomStokesTemporalSymbol_{k,j}
        = 
        s^{\trialVelocityHom}
        \lp
        \basisFunctionPressureFOM{k},
        \basisFunctionVelocityFOM{j};
        \trialVelocityDiscrete,
        \paramVec
        \rp_{\domainSpaceTime},
    \end{aligned}
\end{align*}
and the \gls{rhs} vectors read
\begin{align*}
    \begin{aligned}[t]
        \vectorRHSTemporal &= \{ \vectorRHSTemporalSymbol_i\}, 
        \text{ with }
        \vectorRHSTemporalSymbol_i 
        =
        - \lp \basisFunctionVelocityFOM{i}, \ddt{\liftingVelocityDiscrete};\paramVec \rp_{\domainSpaceTime}, 
        \\
        \vectorRHSVelocityNonlinear &= \{\vectorRHSVelocityNonlinearSymbol_i\},
        \text{ with }
        \vectorRHSVelocityNonlinearSymbol_i = 
        - \bilinearFormViscousStressParam
        {\basisFunctionVelocityFOM{i}}
        {\liftingVelocityDiscrete}
        {\trialVelocityDiscrete}
        {\paramVec},
        \\
        \vectorRHSStabMomStokesTemporal &= \{ \vectorRHSStabMomStokesTemporalSymbol_{k}\}
        \text{, with }
        \vectorRHSStabMomStokesTemporalSymbol_{k}
        = 
        -s^{\vek{l}}
        \lp
        \basisFunctionPressureFOM{k},
        \liftingVelocityDiscrete;
        \trialVelocityDiscrete,
        \paramVec
        \rp_{\domainSpaceTime}.
    \end{aligned}
    \quad
    \begin{aligned}[t]
        \vectorRHSVelocity &= \{\vectorRHSVelocitySymbol_i\},
        \text{ with }
        \vectorRHSVelocitySymbol_i
        =
        \bilinearFormNeumannSpaceTimeParam
        {\basisFunctionVelocityFOM{i}}
        {\nval}
        {\paramVec},
        \\
        \vectorRHSPressure &= \{G_k\},
        \text{ with }
        G_k
        =
        -\bilinearFormPressureParam
        {\liftingVelocityDiscrete}
        {\basisFunctionPressureFOM{k}}
        {\paramVec}_{\domainSpaceTime},
    \end{aligned} 
\end{align*}
For convenience, the dimensions of the matrices and vectors are summarized in \Cref{tab:dimensionsFOM}.
\begin{table}
    \centering
    \begin{tabular}{lclc}
        \toprule
        \acrshort{lhs} Matrices & Dimensions & \acrshort{rhs} Vectors & Dimensions\\
        \midrule
        $\matrixTemporal$, $\matrixViscousStress$ & $\mathbb{R}^{\nBasisVelocityFOM \times \nBasisVelocityFOM}$ & $\vectorRHSTemporal$, $\vectorRHSVelocity$, $\vectorRHSVelocityNonlinear$ & $\mathbb{R}^{\nBasisVelocityFOM}$  \\
        $\matrixPressure$, $\matrixStabMomStokesTemporal$ & $\mathbb{R}^{\nBasisPressureFOM \times \nBasisVelocityFOM}$ & $\vectorRHSPressure$, $\vectorRHSStabMomStokesTemporal$ & $\mathbb{R}^{\nBasisPressureFOM}$  \\
        $\matrixStabMomStokes$ & $\mathbb{R}^{\nBasisPressureFOM \times \nBasisPressureFOM}$ \\
        \bottomrule
    \end{tabular}
    \caption{Dimensions of \acrshort{lhs} matrices and \acrshort{rhs} vectors for the \acrshort{fom}.}
    \label{tab:dimensionsFOM}
\end{table}

\section{Reduced-Order Model for Fluid Flow in Deforming Domains}
\label{sec:reducedOrderModel}
Now that the \gls{fom} has been defined, we can turn to the \gls{mor}. We will start with a description of the underlying ideas of \gls{mor} in \Cref{subsec:methodologicalBackgroundROM} before the \gls{rom} is eventually constructed via \gls{pod} with subsequent projection in \Cref{subsec:discreteFormulationROM}.

\subsection{Methodological Background for the Reduced-Order Model}
\label{subsec:methodologicalBackgroundROM}
In the context of \gls{mor} for numerical schemes, one can distinguish between \textit{interpolation}- and \textit{projection}-based approaches \cite{Manzoni2012}. Methods of the former type try to directly capture input-output relations on the basis of data that comes from numerical simulations or measurements. The latter project the underlying governing equations onto a reduced space that has been constructed before. The class of projection-based methods further entails \textit{certified \gls{rb}} methods \cite{Hesthaven2015} or \textit{\gls{pod}-projection} methods \cite{Chinesta2011,Hesthaven2015,Manzoni2012}. It is also worthwhile to mention the \textit{\gls{pgd}}~\cite{Chinesta2011} here.
\par
In view of the fact that we would like to focus on the class of unsteady problems in deforming domains, one can note that the application of \gls{mor} techniques requires some thorough treatment in this case.
\par
To account for the unsteady character of the problem in the reduction process, there exist two alternative ways of handling time \cite{Glas2017}. The first one follows the classical time-stepping approach, leading to so-called greedy or \gls{pod}-greedy methods, which have been applied to linear as well as to nonlinear problems \cite{Drohmann2012, Fick2018, Grepl2005a, Grepl2012, Haasdonk2008, Haasdonk2013, Sleeman2022}. The second one is rather based on space-time formulations for which reduction techniques for steady problems have been extended appropriately \cite{Urban2012, Urban2013, Tamma2018, Yano2014a, Yano2014b}. This perspective has been used to tackle, e.g., time-dependent optimal control problems \cite{Strazzullo2020, Strazzullo2022}. In addition, a space-time \gls{rom} for sub-intervals of the total simulation time has been presented in \cite{Fritzen2018a}. In distinction to the former methods, we will apply the \gls{pod}-projection approach to the time-continuous space-time formulation stated beforehand.
\par
The construction of a \gls{rom} for deforming domain problems is usually quite involved and requires some careful treatment of the deformations applied to the computational mesh \cite{Anttonen2001,Forti2014}. This is due to the fact that the \gls{fem} function spaces are inherently linked with the geometry of the underlying grid. A strategy based on a mapping functional to relate the time-dependent solution in one- and two-dimensional deforming domains to fixed reference domains has been presented in \cite{Izadi2013}. Furthermore, temporally-local eigenfunctions have been used for \gls{mor} in \cite{Narasingam2017}. 
To this respect, the \gls{csst} approach (cf. \Cref{subsec:discreteFormulationFOM}) offers an appealing but straightforward alternative, since all deformations --- as long as they are prescribed or known a-priori --- are already integrated in the computational mesh and, thereby, considered in the original function spaces. Also, it is worthwhile to note that this approach even works in the presence of spatial topology changes and for two- and three-dimensional spatial domains without any further adaptions.

\subsection{Discrete Formulation for the Reduced-Order Model}
\label{subsec:discreteFormulationROM}
Before we can construct an efficient \gls{rom} whose complexity is independent of the \gls{fom}, we will have to address the non-linearity appearing in our problem through the viscosity model as well as through the formulation of the stabilization parameter.  
For this purpose, we make use of the \gls{eim} \cite{Barrault2004, Grepl2007} here and introduce the following approximations for the viscosity and the stabilization parameter:
\begin{align*}
    \visc\left(\trialVelocityDiscrete, \paramVec\right) 
    \approx
    \sum_{q=1}^{Q_{\visc}} c_q^{\visc}\lp\trialVelocityDiscrete, \paramVec\rp h_q^{\visc}\lp\x\rp,
    &&
    \tMom
    \approx
    \sum_{q=1}^{Q_{\tau}} c_q^{\tau}\lp\trialVelocityDiscrete, \paramVec\rp  h_q^{\tau}\lp\x\rp
    .
\end{align*}
Note that each of the $Q_{\star}$ terms consists of parameter-dependent coefficients $c_q^{\star}\lp\trialVelocity, \paramVec\rp$ and parameter-independent basis functions $h_q^{\star}\lp\x\rp$, where $\star\in\{\visc,\tau\}$.
Consequently, the viscous term is replaced by
\begin{align}
    \bilinearFormViscousStressParam{\weightVelocityDiscrete}{\trialVelocityHomDiscrete}{\trialVelocityDiscrete}{\paramVec}_{\domainSpaceTime}
    \approx
    \sum_{q=1}^{Q_{\visc}} c_q^{\visc}\lp\trialVelocityDiscrete, \paramVec\rp 
    a_q\lp\weightVelocityDiscrete, \trialVelocityHomDiscrete\rp_{\domainSpaceTime},
\end{align}
with
\begin{align}
    a_q\lp\weightVelocityDiscrete, \trialVelocityHomDiscrete\rp_{\domainSpaceTime}
    =
    \intDomainSpaceTime
    {
        \ros{\weightVelocityDiscrete} : 2 h_q^{\visc}\lp\x\rp \ros{\trialVelocityHomDiscrete}
    }.
\end{align}
Moreover, the stabilization term $s^{\trialVelocityHom}$ is approximated via
\begin{align}
    s^{\trialVelocityHom}\lp\weightPressureDiscrete,\trialVelocityHomDiscrete;\trialVelocityDiscrete, \paramVec\rp_{\domainSpaceTime}
    \approx
    \sum_{q=1}^{Q_{\tau}} c_q^{\tau}\left(\trialVelocityDiscrete, \paramVec\right)
    s^{\trialVelocityHom}_q\lp\weightPressureDiscrete,\trialVelocityHomDiscrete\rp_{\domainSpaceTime},
\end{align}
where
\begin{align}
    s^{\trialVelocityHom}_q\lp\weightPressureDiscrete,\trialVelocityHomDiscrete\rp_{\domainSpaceTime}
    =
    \sum_{e} \intElementSpaceTime{\basisFunctionEIM{q}{}^{\tau}\lp\x\rp \frac{1}{\rho}
    \lp
    -\gr{\weightPressureDiscrete}
    \rp
    \cdot
    \lp
    \rho \ddt{\trialVelocityHomDiscrete}
    \rp
    }.
\end{align}
Analogously, the remaining stabilization terms 
$s^{\vek{l}}
    (\weightPressureDiscrete,
    \liftingVelocityDiscrete;
    \trialVelocityDiscrete,
    \paramVec)_{\domainSpaceTime}$ 
    and
$s^{\trialPressure}(
    \weightPressureDiscrete,
    \trialPressureDiscrete;
    \trialVelocityDiscrete
    \paramVec)_{\domainSpaceTime}$ 
can be approximated using parameter-independent terms
$s_q^{\vek{l}}
    (\weightPressureDiscrete,
    \liftingVelocityDiscrete)_{\domainSpaceTime}$
    and
$s_q^{\trialPressure}(
    \weightPressureDiscrete,
    \trialPressureDiscrete)_{\domainSpaceTime}$,
respectively. As a consequence, all the matrices and vectors in the algebraic system of the \gls{fom} relying on these forms, i.e., $\matrixViscousStress$, $\matrixStabMomStokesTemporal$, $\matrixStabMomStokes$, $\vectorRHSVelocityNonlinear$, and $\vectorRHSStabMomStokesTemporal$ are replaced by approximations using $\matrixViscousStress^q$, $\matrixStabMomStokesTemporal^q$, $\matrixStabMomStokes^q$, $\vectorRHSVelocityNonlinear^q$, and $\vectorRHSStabMomStokesTemporal^q$ corresponding to the respective terms above. For example, this means $\matrixViscousStress^q = \lb \matrixViscousStressSymbol^q_{i,j}\rb,
        \text{ with } 
        \matrixViscousStressSymbol^q_{i,j} 
        = 
        a_q\lp\basisFunctionVelocityFOM{i},\basisFunctionVelocityFOM{j}\rp_{\domainSpaceTime}$.
\bigskip\par
To perform the projection step later, a basis spanning the reduced spaces is needed. To that end, we apply the \gls{pod} using the method of snapshots~\cite{Sirovich1987}, i.e., solutions of the \gls{fom}. In particular, this is done individually for the homogeneous velocity $\trialVelocityHomDiscrete$ and the pressure $\trialPressureDiscrete$ leading to the \textit{reduced finite-dimensional function spaces} $\weightingSpaceVelocityReduced \subset \weightingSpaceVelocityDiscrete$ and $\weightingSpacePressureReduced \subset \trialSpacePressureDiscrete$.
As a result of the \gls{pod}, we obtain $\nBasisVelocityHomROM$ and $\nBasisPressureROM$ basis functions for the reduced representation $(\trialVelocityHomReduced,\trialPressureReduced)$ of the homogeneous velocity and pressure field, respectively. To account for the Dirichlet boundary conditions, the basis for $\weightingSpaceVelocityReduced$ is augmented with the lifting function(s) $\liftingVelocityDiscrete$ yielding the reduced space $\trialSpaceVelocityReduced \subset \trialSpaceVelocityDiscrete$ for the reduced velocity $\trialVelocityReduced$. Therefore, we use $\nBasisVelocityROM \ge \nBasisVelocityHomROM$ to denote the number of basis functions for $\trialVelocityReduced$. We sort the basis functions in descending order of significance --- indicated by the magnitude of the corresponding eigenvalues --- while the lifting functions are always leading to ensure that the Dirichlet boundary conditions are met, even if we only use a subset of these basis functions.
\par
For all basis functions, we collect the coefficients with respect to the \gls{fom} function spaces in the so-called \textit{basis function matrices} $\basisFunctionMatrixVelocity\in\mathbb{R}^{\nBasisVelocityFOM \times \nBasisVelocityROM}$ and $\basisFunctionMatrixPressure\in\mathbb{R}^{\nBasisPressureFOM \times \nBasisPressureROM}$. These matrices are multiplied with those from the algebraic system of the \gls{fom} given in \Cref{eq:algebraicSystemFOM}, which yields the projection of the corresponding operators onto the reduced space. Thus, one can formulate the algebraic system for the vectors of unknowns $\velocityDOFVectorReduced\in\mathbb{R}^{\nBasisVelocityROM}$ and $\pressureDOFVectorReduced\in\mathbb{R}^{\nBasisPressureROM}$, where $\nBasisVelocityROM$ and $\nBasisPressureROM$ are the number of unknowns of the \gls{rom}. Key assumption for an effective reduction is that it holds that $\nBasisROM = \nBasisVelocityROM + \nBasisPressureROM \ll \nBasisFOM = \nBasisVelocityFOM + \nBasisPressureFOM$. The system finally reads:
\begin{align}
    \lb
    \begin{array}{cc}
        \matrixTemporalReduced + \matrixViscousStressReducedParam
         & -\matrixPressureReducedTrans
        \\
        \matrixPressureReduced + \matrixStabMomStokesTemporalReducedParam& \matrixStabMomStokesReducedParam
    \end{array}
    \rb
    \lb
    \begin{array}{c}
        \velocityDOFVectorReducedParam\\
        \pressureDOFVectorReducedParam
    \end{array}
    \rb&
    \nonumber
    \\
    =
    \lb
    \begin{array}{c}
        \vectorRHSTemporalReduced +
        \vectorRHSVelocityReduced + \vectorRHSVelocityNonlinearReducedParam \\
        \vectorRHSPressureReduced
        + \vectorRHSStabMomStokesTemporalReducedParam
    \end{array}
    \rb&,
    \label{eq:ROMStokesShearThinningDeformingDomainsAlgebraicSystem}
\end{align}
where the \gls{lhs} matrices read
\begin{align*}
    \begin{aligned}[t]
        \matrixTemporalReduced
        &=
        \basisFunctionMatrixVelocityTrans
        \matrixTemporal
        \basisFunctionMatrixVelocity,
        \\
        \matrixPressureReduced
        &=
        \basisFunctionMatrixPressureTrans
        \matrixPressure
        \basisFunctionMatrixVelocity,
    \end{aligned}
    \quad
    \begin{aligned}[t]
        \matrixViscousStressReducedParam
        &=
        \sum_{q=1}^{Q_{\visc}}
        \coeffEIM{q}^{\visc}\lp\trialVelocityDiscrete,\paramVec\rp
        \mat{A}^q_{\superscriptROM},
        \text{ and }
        \mat{A}^q_{\superscriptROM}
        =
        \basisFunctionMatrixVelocityTrans
        \mat{A}^q
        \basisFunctionMatrixVelocity,
        \\
        \matrixStabMomStokesTemporalReducedParam
        &=
        \sum_{q=1}^{Q_{\tau}}
        \coeffEIM{q}^{\tau}\lp\trialVelocityDiscrete,\paramVec\rp
        \matrixStabMomStokesTemporalReduced^q,
        \text{ and }
        \matrixStabMomStokesTemporalReduced^q
        =
        \basisFunctionMatrixPressureTrans
        \matrixStabMomStokesTemporal^q
        \basisFunctionMatrixVelocity,
        \\
        \matrixStabMomStokesReducedParam
        &=
        \sum_{q=1}^{Q_{\tau}}
        \coeffEIM{q}^{\tau}\lp\trialVelocityDiscrete,\paramVec\rp
        \mat{S}^q_{\superscriptROM},
        \text{ and }
        \mat{S}^q_{\superscriptROM}
        =
        \basisFunctionMatrixPressureTrans
        \mat{S}^q
        \basisFunctionMatrixPressure,
    \end{aligned}
\end{align*}
and the \gls{rhs} vectors are given as
\begin{align*}
    \begin{aligned}[t]
        \vectorRHSTemporalReduced
        &=
        \basisFunctionMatrixVelocityTrans
        \vectorRHSTemporal,
        \\
        \vectorRHSVelocityReduced
        &=
        \basisFunctionMatrixVelocityTrans
        \vectorRHSVelocity,
        \\
        \vectorRHSPressureReduced
        &=
        \basisFunctionMatrixPressureTrans
        \vectorRHSPressure,
    \end{aligned}
    \quad
    \begin{aligned}[t]
        \vectorRHSVelocityNonlinearReducedParam
        &=
        \sum_{q=1}^{Q_{\visc}}
        \coeffEIM{q}^{\visc}\lp\trialVelocityDiscrete,\paramVec\rp
        \vectorRHSVelocityNonlinearReduced^q,
        \text{ and }
        \vectorRHSVelocityNonlinearReduced^q
        =
        \basisFunctionMatrixVelocityTrans
        \vectorRHSVelocityNonlinear^q,
        \\
        \vectorRHSStabMomStokesTemporalReducedParam
        &=
        \sum_{q=1}^{Q_{\tau}}
        \coeffEIM{q}^{\tau}\lp\trialVelocityDiscrete,\paramVec\rp
        \vectorRHSStabMomStokesTemporalReduced^q,
        \text{ and }
        \vectorRHSStabMomStokesTemporalReduced^q
        =
        \basisFunctionMatrixPressureTrans
        \vectorRHSStabMomStokesTemporal^q.
    \end{aligned}
\end{align*}
\Cref{tab:dimensionsROM} gives an overview over the dimensions of the matrices and vectors for the \gls{rom}. Due to the reduced dimensions, the solution of the algebraic system of the \gls{rom} for any new parameter sample $\paramVec$ is possible with significantly less computational resources. Thanks to the \gls{eim}, the same holds for the assembly process of this system. Note, however, that due to the specific implementation of the full-order solver, we keep the dependency on $\trialVelocityDiscrete$ in the \gls{eim} coefficients.
\begin{table}
    \centering
    \begin{tabular}{lclc}
        \toprule
        \acrshort{lhs} Matrices & Dimensions & \acrshort{rhs} Vectors & Dimensions\\
        \midrule
        $\matrixTemporalReduced$, $\matrixViscousStressReduced$ & $\mathbb{R}^{\nBasisVelocityROM \times \nBasisVelocityROM}$ & $\vectorRHSTemporalReduced$, $\vectorRHSVelocityReduced$, $\vectorRHSVelocityNonlinearReduced$ & $\mathbb{R}^{\nBasisVelocityROM}$  \\
        $\matrixPressureReduced$, $\matrixStabMomStokesTemporalReduced$ & $\mathbb{R}^{\nBasisPressureROM \times \nBasisVelocityROM}$ & $\vectorRHSPressureReduced$, $\vectorRHSStabMomStokesTemporalReduced$ & $\mathbb{R}^{\nBasisPressureROM}$  \\
        $\matrixStabMomStokesReduced$ & $\mathbb{R}^{\nBasisPressureROM \times \nBasisPressureROM}$ \\
        \bottomrule
    \end{tabular}
    \caption{Dimensions of \acrshort{lhs} matrices and \acrshort{rhs} vectors for the \acrshort{rom}.}
    \label{tab:dimensionsROM}
\end{table}
\section{Numerical Results}
\label{sec:numericalResults}
In the following, we will illustrate how the proposed \gls{rom} approach, which makes use of a time-continuous space-time setting, can help to significantly decrease the required computational resources for time-dependent parametric problems defined in deforming domains. To that end, we present error and performance analysis results for two test cases. The test cases are representative of applications that may arise from engineering or biomedical problems, for example. The first test case involves a two-dimensional valve-like geometry, resulting in a 3D space-time geometry. In space, the valve plug moves over time and even closes off parts of the geometry, which means that the spatial topology changes. In the second test case, we consider a three-dimensional artery-like geometry. Consequently, the resulting space-time domain is 4D. In the center region, the geometry is compressed over time, which yields a deforming domain problem.
\subsection{Valve-Like Geometry with Topology Changes}
\begin{figure}
\centering
\subcaptionbox{$t=0.0$}{
\begin{tikzpicture}[
    axis/.style={ -latex},
    every node/.style={color=black}
    ]
       \node at (0,0) [inner sep=0pt, anchor=south west] {\includegraphics[width=0.4\textwidth,trim={0cm 0cm 15cm 0cm},clip]{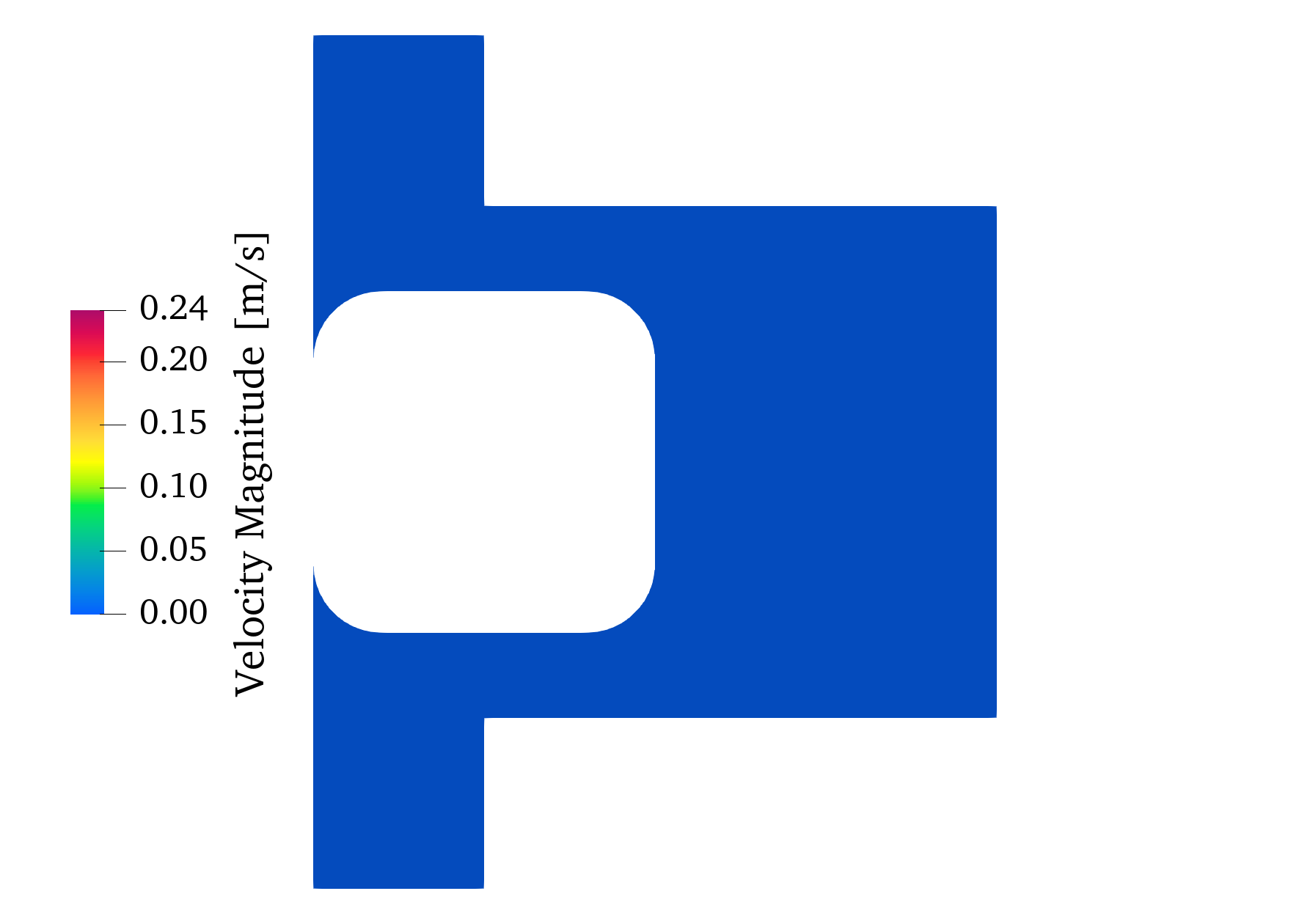}};
       \draw[axis] (0.1,0.2)  -- (0.6,0.2) node[above]{$x$};
       \draw[axis] (0.1,0.2)  -- (0.1,0.7) node[left]{$y$};
\end{tikzpicture}

}
\subcaptionbox{$t=0.45$}{
\begin{tikzpicture}[
    axis/.style={ -latex},
    every node/.style={color=black}
    ]
       \node at (0,0) [inner sep=0pt, anchor=south west] {\includegraphics[width=0.4\textwidth,trim={0cm 0cm 15cm 0cm},clip]
{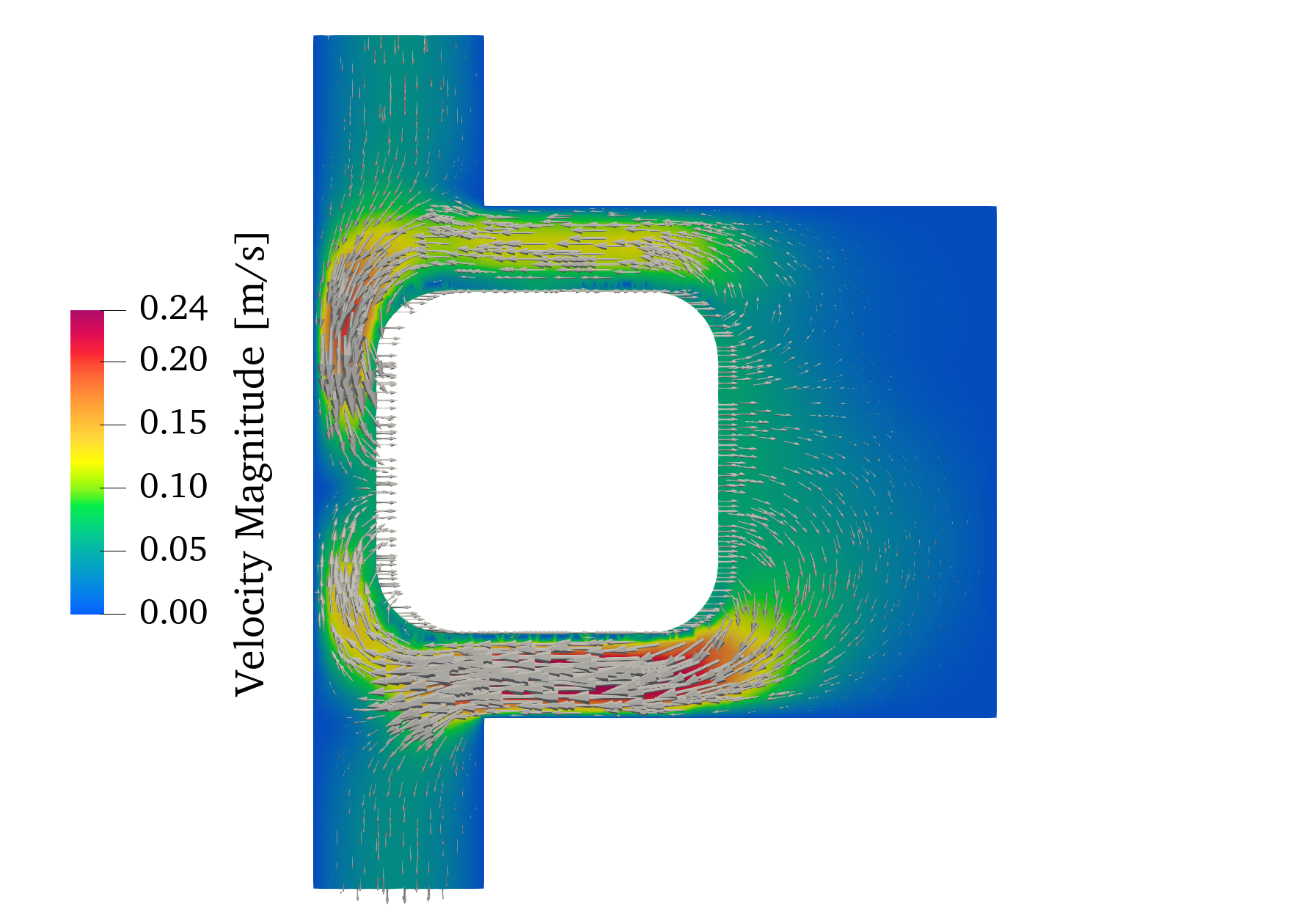}};
       \draw[axis] (0.1,0.2)  -- (0.6,0.2) node[above]{$x$};
       \draw[axis] (0.1,0.2)  -- (0.1,0.7) node[left]{$y$};
\end{tikzpicture}

}
\\
\subcaptionbox{$t=0.9$}{
\begin{tikzpicture}[
    axis/.style={ -latex},
    every node/.style={color=black}
    ]
       \node at (0,0) [inner sep=0pt, anchor=south west] {\includegraphics[width=0.4\textwidth,trim={0cm 0cm 15cm 0cm},clip]
{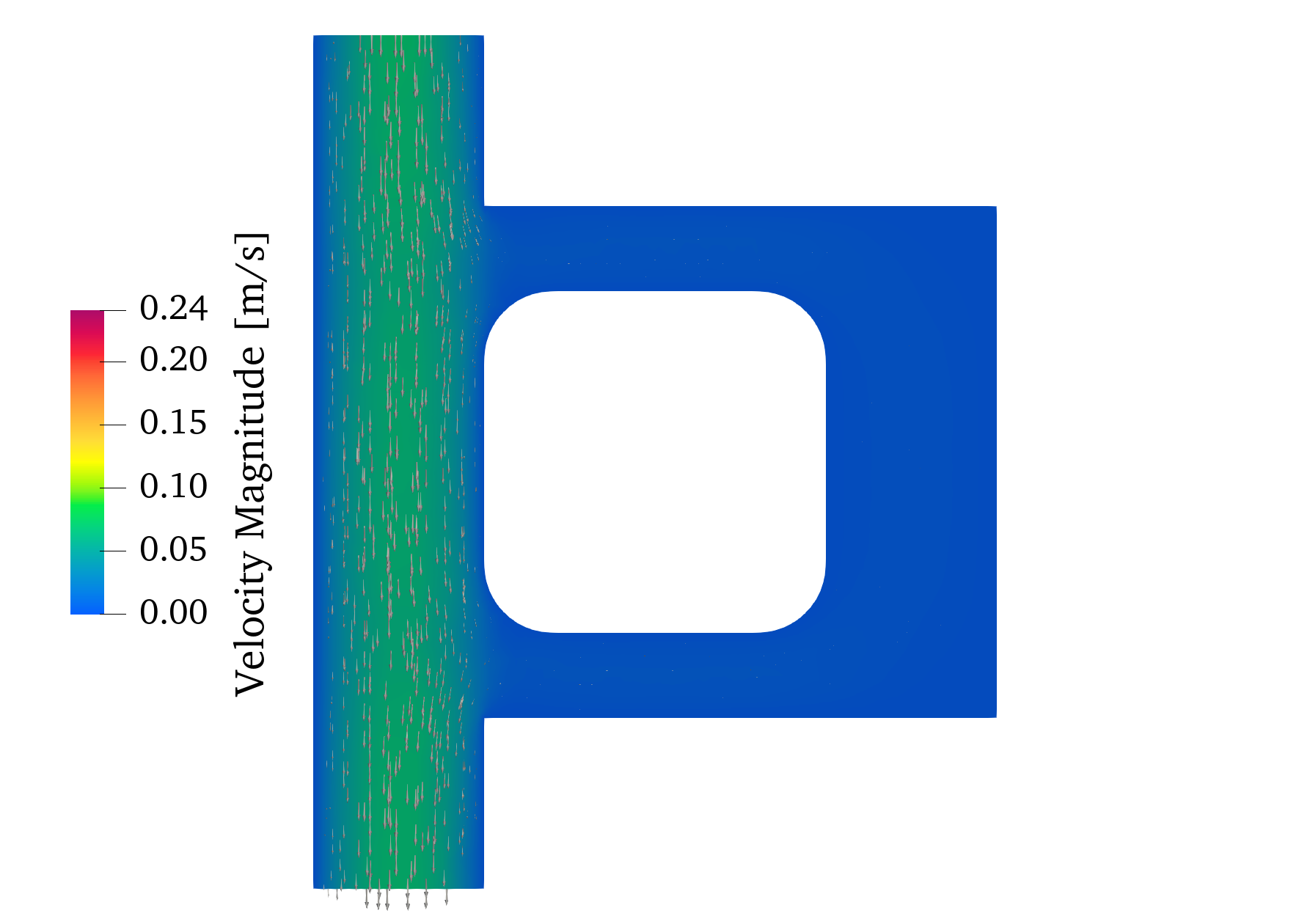}};
       \draw[axis] (0.1,0.2)  -- (0.6,0.2) node[above]{$x$};
       \draw[axis] (0.1,0.2)  -- (0.1,0.7) node[left]{$y$};
\end{tikzpicture}

}
\\
\subcaptionbox{$t=1.35$}{
\begin{tikzpicture}[
    axis/.style={ -latex},
    every node/.style={color=black}
    ]
       \node at (0,0) [inner sep=0pt, anchor=south west] {\includegraphics[width=0.4\textwidth,trim={0cm 0cm 15cm 0cm},clip]{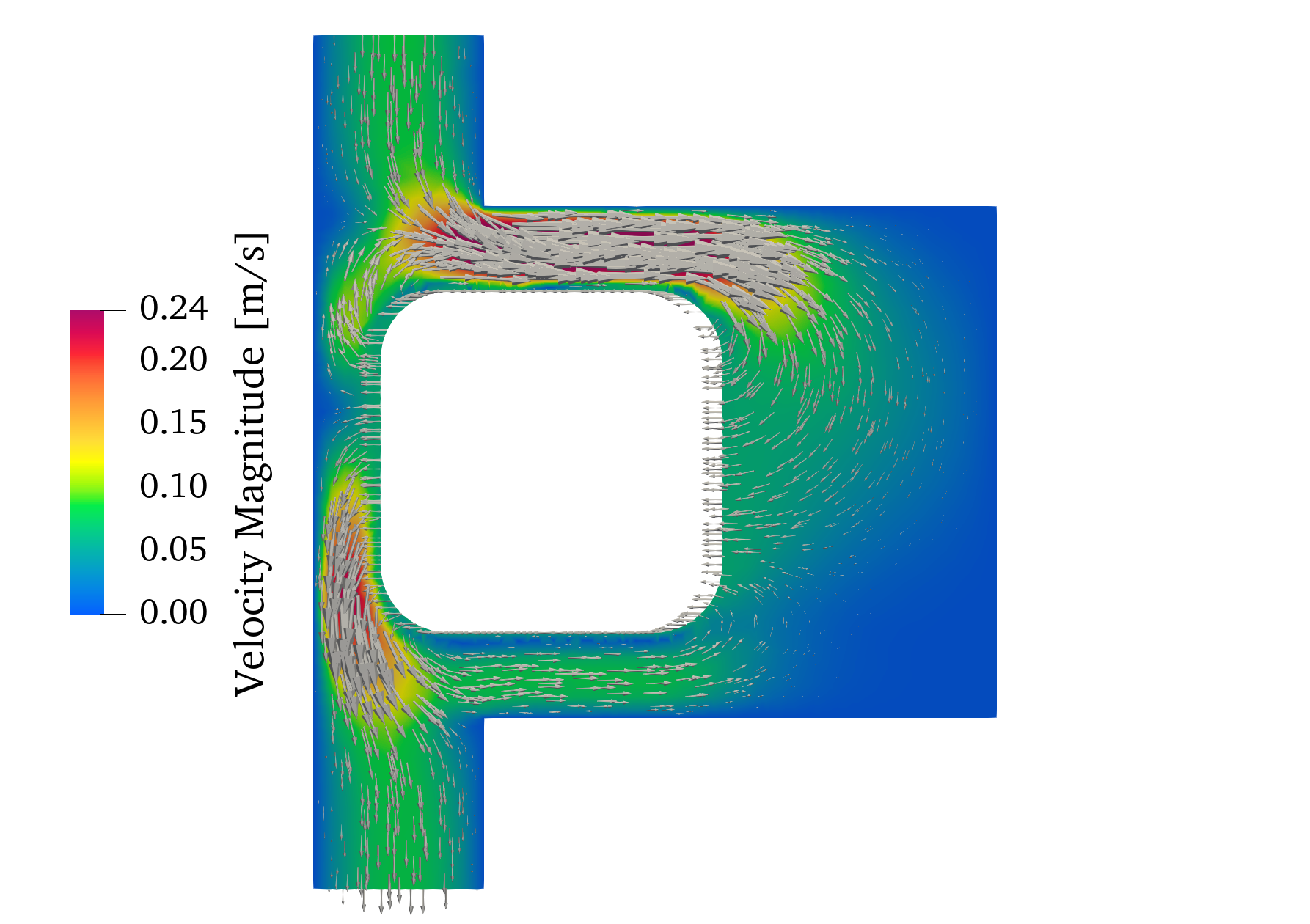}};
       \draw[axis] (0.1,0.2)  -- (0.6,0.2) node[above]{$x$};
       \draw[axis] (0.1,0.2)  -- (0.1,0.7) node[left]{$y$};
\end{tikzpicture}

}
\subcaptionbox{$t=1.8$}{
\begin{tikzpicture}[
    axis/.style={ -latex},
    every node/.style={color=black}
    ]
       \node at (0,0) [inner sep=0pt, anchor=south west] {\includegraphics[width=0.4\textwidth,trim={0cm 0cm 15cm 0cm},clip]
{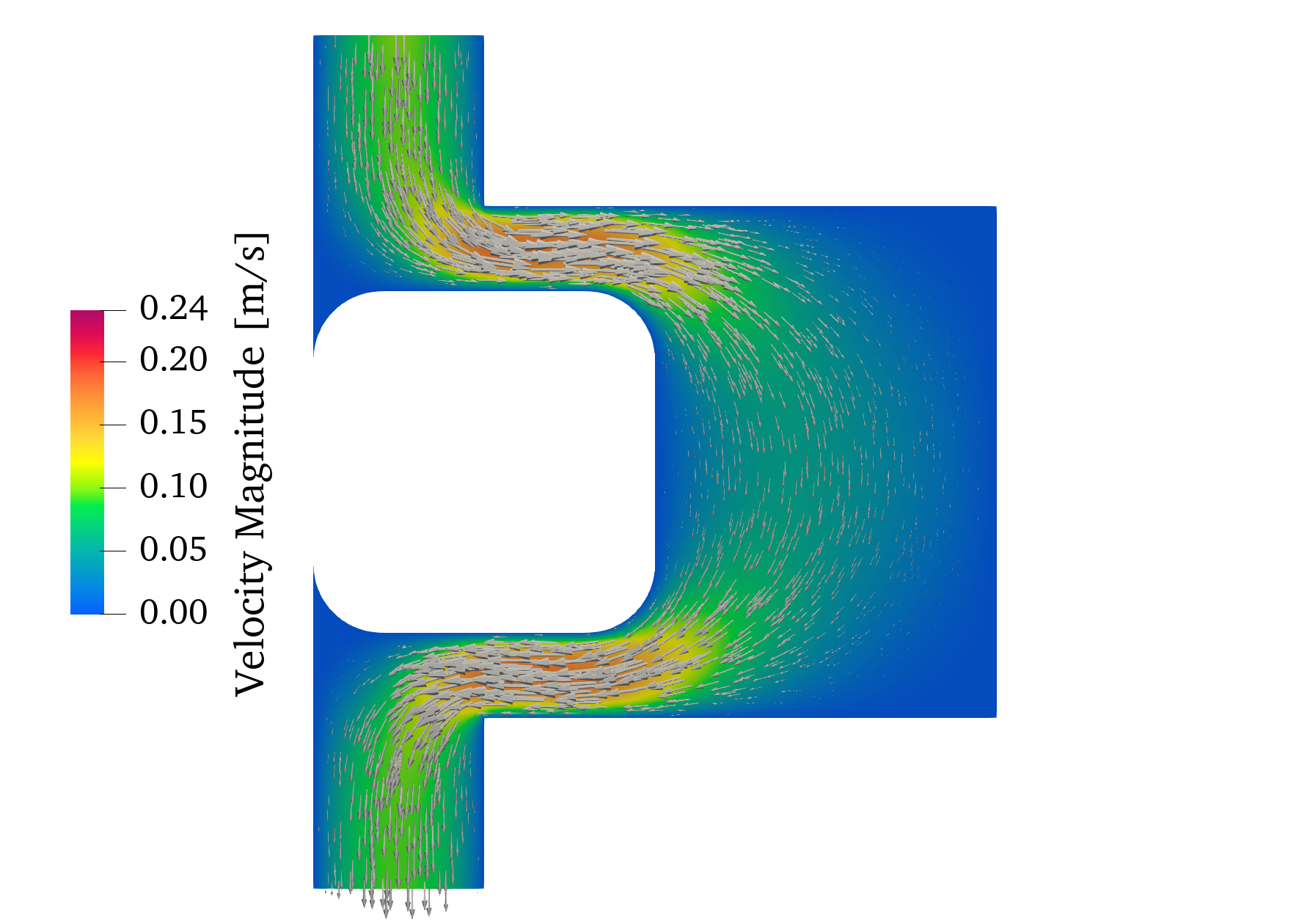}};
       \draw[axis] (0.1,0.2)  -- (0.6,0.2) node[above]{$x$};
       \draw[axis] (0.1,0.2)  -- (0.1,0.7) node[left]{$y$};
\end{tikzpicture}

}
\caption{Valve-like test case: flow field for different points in time.}
\label{fig:valveFlowField}
\end{figure}
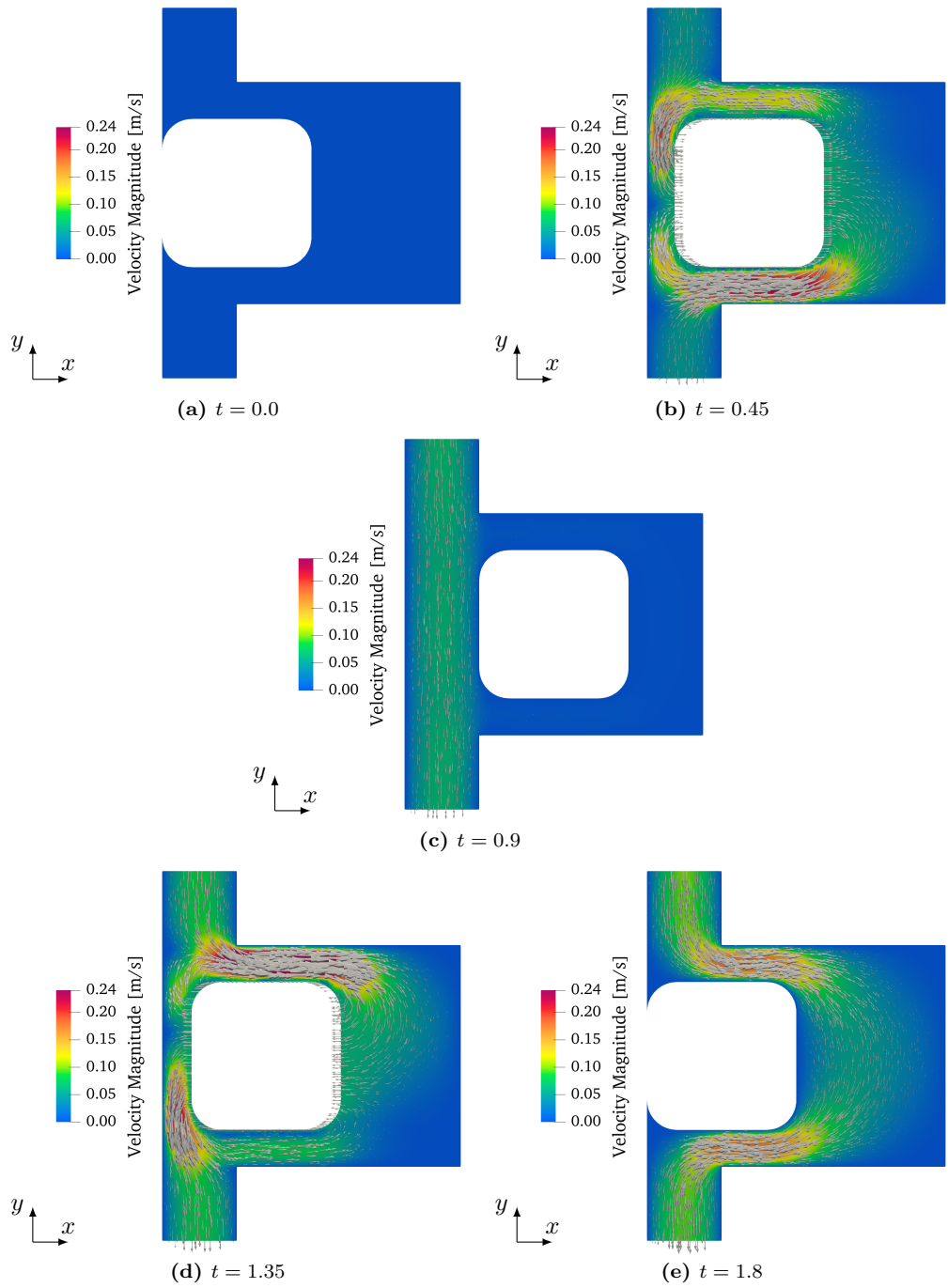

\newcommand{\myD}{0.9\textwidth}
In this section, we consider a two-dimensional valve-like geometry composed of a valve plug encased in a flow channel. Its initial configuration for $t=\SI{0}{\second}$ is depicted in \Cref{subfig:valveSpatialDomain}. At the top, fluid can enter the channel through an inlet with width $L_{\text{inlet}}=\SI{0.025}{\meter}$. The outlet is located at the bottom. All remaining boundary portions are impermeable. We consider a time interval such that $t \in [0,1.8]\,\SI{}{\second}$. With the passage of time, the plug starts to move inward for $t \in [0.3,0.7)\,\SI{}{\second}$, opening a second branch for the flow on the left-hand side. After reaching the center of the casing, the plug stays at rest for $t \in [0.7,1.1)\,\SI{}{\second}$ and, still a bit later, the movement is reversed for $t \in [1.1,1.5)\,\SI{}{\second}$ to close the emerged branch for the remainder of the simulation. The resulting three-dimensional space-time domain is shown in \Cref{subfig:valveSpaceTimeDomain}, while the spatial geometry for different time instances can be seen in \Cref{fig:valveFlowField}. The velocity of the plug in $x$-direction is $u_{\text{plug}} = \pm \SI{0.0625}{\meter\per\second}$ for the inward and outward movement, respectively.
\newcommand{\myw}{0.25}

\newcommand{\myY}{4.7}
\newcommand{\xZero}{1.1}
\newcommand{\deltaX}{0.71}
\newcommand{\deltaXGap}{1.35}

\newcommand{\myX}{4.5}
\newcommand{\yZero}{0.54}
\newcommand{\deltaY}{0.71}

\begin{figure}
\captionsetup[sub]{position=bottom}
\centering
\subcaptionbox{Spatial domain for $t=0.0$.\label{subfig:valveSpatialDomain}}{
\begin{tikzpicture}[
    axis/.style={ -latex},
    every node/.style={color=black}
    ]
       \node at (0,0) [inner sep=0pt, anchor=south west] {\includegraphics[width=\myw\textwidth,trim={0cm 0cm 12cm 0cm},clip]{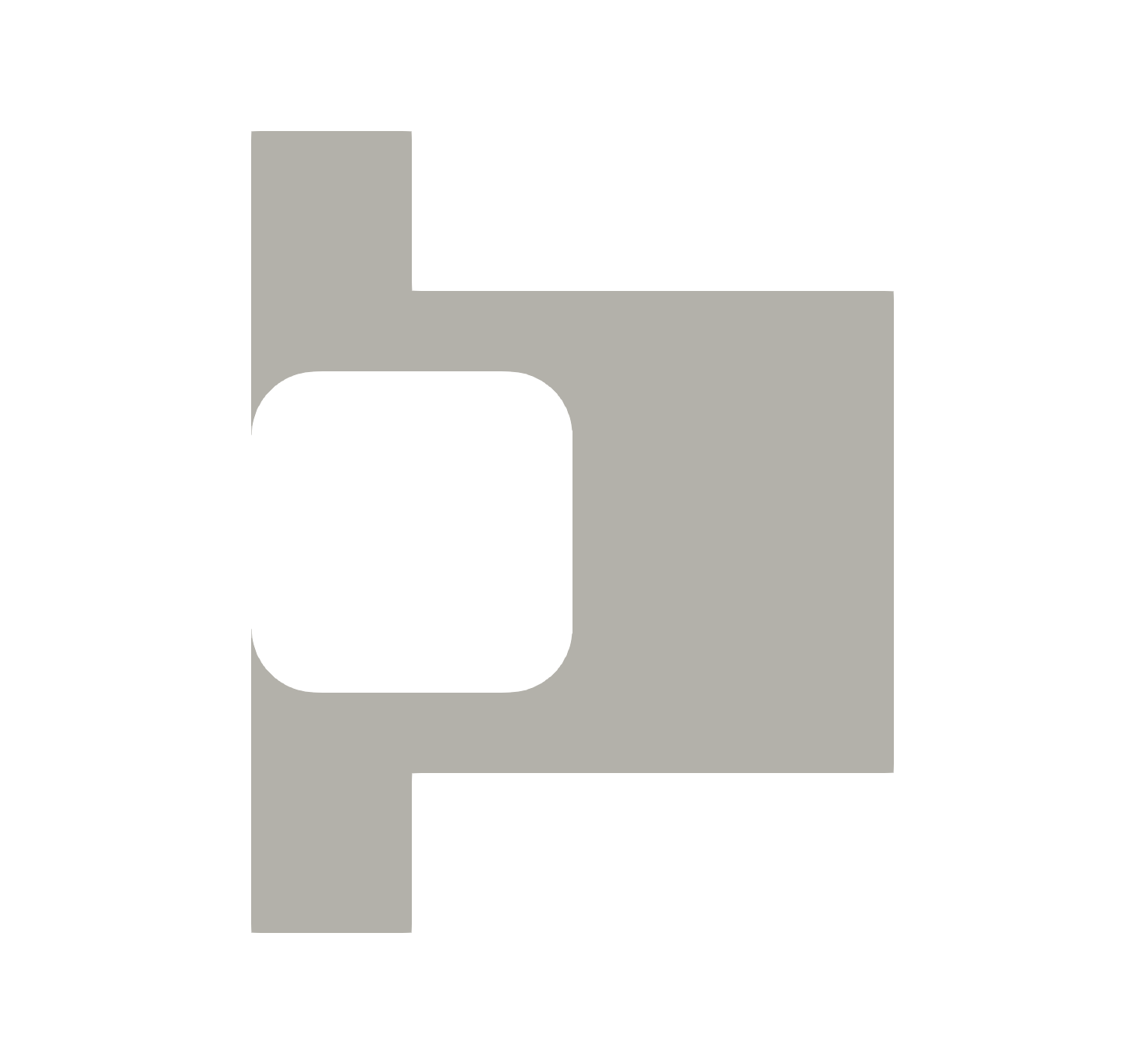}};
       \draw[axis] (0.1,0.2)  -- (0.6,0.2) node[above]{$x$};
       \draw[axis] (0.1,0.2)  -- (0.1,0.7) node[left]{$y$};      

       \dimline[line style={line width=0.75pt},label style={below=2.25em,rotate=-90,anchor=center},extension start length=-0.6,extension end length=-0.6]{(\myX,\yZero)}{(\myX,\yZero+\deltaY)}{\footnotesize \SI{0.025}{\meter}};
       \dimline[line style={line width=0.75pt},label style={below=2.5em,rotate=-90,anchor=center},extension start length=-0.3,extension end length=-0.3]{(\myX,\yZero+\deltaY)}{(\myX,\yZero+1.5*\deltaY)}{\footnotesize \SI{0.0125}{\meter}};
       \dimline[line style={line width=0.75pt},label style={below=2.1em,rotate=-90,anchor=center},extension start length=-0.3,extension end length=-0.3]      {(\myX,\yZero+1.5*\deltaY)}{(\myX,\yZero+3.5*\deltaY)}{\footnotesize \SI{0.05}{\meter}};
       \dimline[line style={line width=0.75pt},label style={below=2.5em,rotate=-90,anchor=center},extension start length=-0.3,extension end length=-0.3]{(\myX,\yZero+3.5*\deltaY)}{(\myX,\yZero+4*\deltaY)}{\footnotesize \SI{0.0125}{\meter}};
     \dimline[line style={line width=0.75pt},label style={below=2.25em,rotate=-90,anchor=center},extension start length=-0.6,extension end length=-0.6]{(\myX,\yZero+4*\deltaY)}{(\myX,\yZero+5*\deltaY)}{\footnotesize \SI{0.025}{\meter}};

       \dimline[line style={line width=0.75pt}, 
       label style={rotate=45, above = 6 ex, anchor=north},
       extension start length=0.6,
       extension end length=0.6]{(\xZero,\myY)}{(\xZero+\deltaX,\myY)}{\footnotesize \SI{0.025}{\meter}};
       
       \dimline[line style={line width=0.75pt}, label style={rotate=45, above = 6 ex, anchor=north},extension start length=-0,extension end length=0]{(\xZero+\deltaX,\myY)}{(\xZero+2*\deltaX,\myY)}{\footnotesize\SI{0.025}{\meter}};
       
       \dimline[line style={line width=0.75pt},label style={rotate=45, above = 6 ex, anchor=north},extension start length=0.3,extension end length=0.3]      {(\xZero+2*\deltaX,\myY)}{(\xZero+2*\deltaX+\deltaXGap,\myY)}{\footnotesize \SI{0.05}{\meter}};

       \dimline[line style={line width=0.75pt},label style={above=2ex}, extension start length=-2,extension end length=-2]{(1.14, 1.74)}{(1.3, 1.9)}{\footnotesize $r=\SI{0.01}{\meter}$};
       
       \draw[axis] (\xZero+\deltaX,\yZero+2.5*\deltaY)  -- (\xZero+2.5*\deltaX,\yZero+2.5*\deltaY) node[above]{$u_{\text{plug}}$};
     
\end{tikzpicture}
}
\subcaptionbox{Space-time domain.\label{subfig:valveSpaceTimeDomain}}{
\begin{tikzpicture}[
    axis/.style={ -latex},
    every node/.style={color=black}
    ]
       \node at (0,0) [inner sep=0pt, anchor=south west] {\includegraphics[width=0.45\textwidth,trim={0cm 0cm 0cm 0cm},clip]{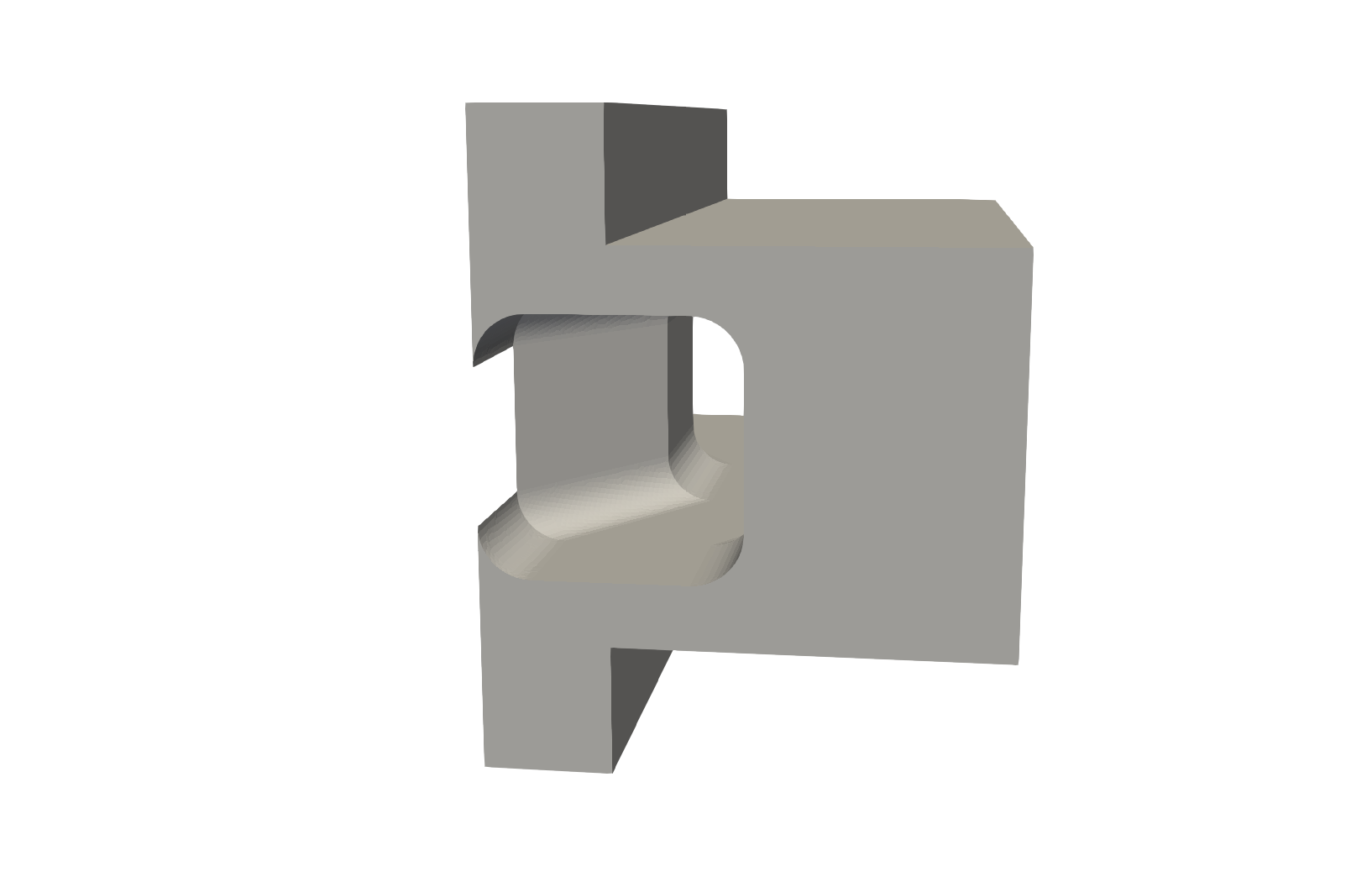}};
      \draw[axis] (1.0,0.1)  -- (1.8,0.05) node(xline)[above]{$x$};
      \draw[axis] (1.0,0.1)  -- (1.0,0.8) node(xline)[right]{$y$};
      \draw[axis] (1.0,0.1)  -- (0.7,-0.2) node(xline)[left]{$t$};

\end{tikzpicture}

}

\caption{Computational domain of valve-like test case.}
\label{fig:valveGeometry}
\end{figure}
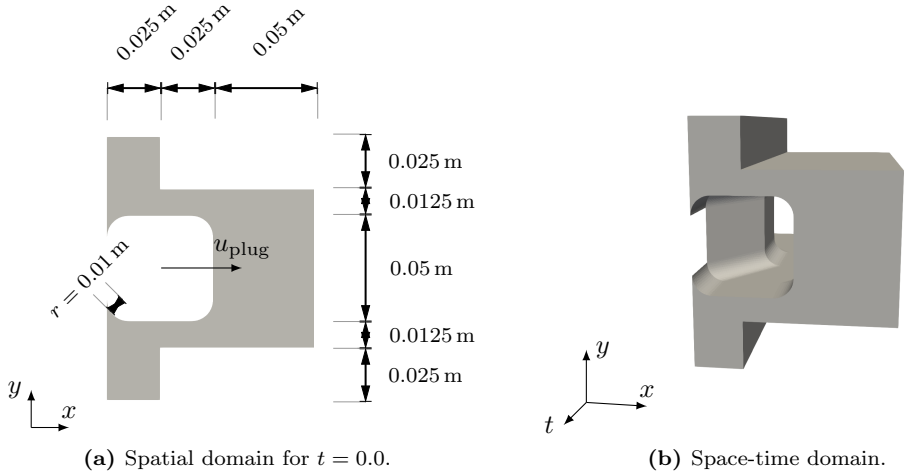
\par
The fluid is supposed to be plastics melt as present, e.g., in various polymer processing techniques for thermoplastics. We use material properties for an exemplary \gls{pc} offered by Covestro Makrolon\textsuperscript{\textregistered}. 
We use the density $\density = \SI{1200}{\kilo\gram \per\cubic\meter}$ and, following \cite{Schroeder2020}, the parameters for the viscosity model from~\Cref{eq:viscosityModel} as
\begin{align}
        \visc_0 = \SI{270}{\pascal\second},\,
        \visc_{\infty} = \SI{0}{\pascal\second},\,
        \lambda = \SI{1.2e-3}{\second},\,
        a = \SI{1}{},\,
        n = \SI{0.775}{}.
\end{align}
In the following, the velocity vector $\vel = \lb \velX, \velY \rb^T$ will collect the velocity components $\velX$ and $\velY$ for the $x$- and $y$-direction, respectively. We prescribe a time-dependent inflow profile given as
\begin{align*}
\velX_{\text{in}}(x,t) &= \SI{0}{\meter\per\second},\\
\velY_{\text{in}}(x,t)& = v_{\text{in}}^0(x(x-\SI{0.025}{\meter}))\sqrt{\frac{t}{\SI{1.8}{\second}}}.   
\end{align*}
Choosing $v_{\text{in}}^0 = \SI{640}{\per\meter\second}$, this leads to a Reynolds number of approximately $Re =  \SI{1e-2}{}$ 
such that the Stokes equations are appropriate to model this creeping flow. For the walls of the casing, no-slip boundary conditions are assumed. At the outlet, parallel outflow is enforced by setting the velocity in $x$-direction to zero, i.e., $\velX_{\text{out}} = \SI{0}{\meter\per\second}$. The resulting flow field is visualized in \Cref{fig:valveFlowField} for several points in time.
\subsubsection{\gls{rom} for the Flow of Plastics Melt in a Valve-Like Geometry}
\label{subsubsec:valve-likeROM}
In the following, we will consider uncertainties in the model parameters from the viscosity model stated in \Cref{eq:viscosityModel}. Therefore, we collect two of those parameters in the parameter vector, i.e., $\paramVec = \lb \lambda,n\rb$. Assuming that variations of $\pm5\%$ may occur, it follows that $\paramVec \in \lb0.95\bar{\lambda}, 1.05\bar{\lambda}\rb \times \lb0.95\bar{n}, 1.05\bar{n}\rb$, with $\bar{\lambda} = \SI{1.2e-3}{\second}$ and $\bar{n} = \SI{0.775}{}$.
\bigskip
\par
\begin{figure}
    \captionsetup[sub]{position=bottom}
    \centering
    \Large
    \subcaptionbox{Distribution of the eigenvalues from the \gls{pod}.\label{fig:valve-likeEigenvalues}}{
        \resizebox{0.44\textwidth}{!}{
        \input{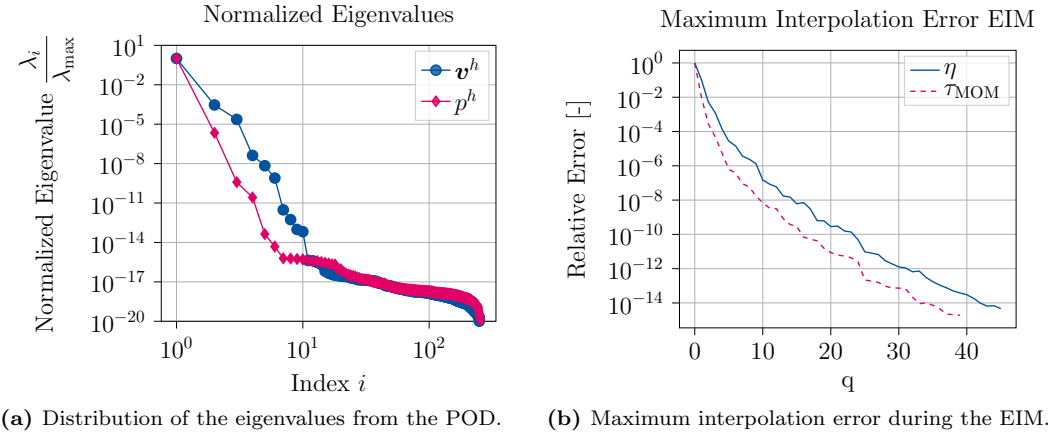}
        }
    }
    \subcaptionbox{Maximum interpolation error during the \gls{eim}.\label{fig:valve-likeEIMErrors}}{
        \resizebox{0.44\textwidth}{!}{
\begin{tikzpicture}


\definecolor{color0}{rgb}{0,0.329411764705882,0.623529411764706}
\definecolor{color1}{rgb}{0.890196078431372,0,0.4}

\begin{axis}[
legend cell align={left},
legend style={fill opacity=0.8, draw opacity=1, text opacity=1, draw=white!80!black},
log basis y={10},
tick align=outside,
tick pos=left,
title={Maximum Interpolation Error EIM},
x grid style={white!69.0196078431373!black},
xlabel={q},
xmajorgrids,
xmin=-2.25, xmax=47.25,
xtick style={color=black},
xtick={-10,0,10,20,30,40,50},
xticklabels={
  \(\displaystyle -10\),
  \(\displaystyle 0\),
  \(\displaystyle 10\),
  \(\displaystyle 20\),
  \(\displaystyle 30\),
  \(\displaystyle 40\),
  \(\displaystyle 50\)
},
y grid style={white!69.0196078431373!black},
ylabel={Relative Error [-]},
ymajorgrids,
ymin=3.51471410776371e-16, ymax=5.44389825525944,
ymode=log,
ytick style={color=black},
ytick={1e-18,1e-16,1e-14,1e-12,1e-10,1e-08,1e-06,0.0001,0.01,1,100,10000},
yticklabels={
  \(\displaystyle 10^{-18}\),
  \(\displaystyle 10^{-16}\),
  \(\displaystyle 10^{-14}\),
  \(\displaystyle 10^{-12}\),
  \(\displaystyle 10^{-10}\),
  \(\displaystyle 10^{-8}\),
  \(\displaystyle 10^{-6}\),
  \(\displaystyle 10^{-4}\),
  \(\displaystyle 10^{-2}\),
  \(\displaystyle 10^{0}\),
  \(\displaystyle 10^{2}\),
  \(\displaystyle 10^{4}\)
}
]
\addplot [thick, color0]
table {%
0 1
1 0.0987673453717658
2 0.00515599137631726
3 0.00124900271315819
4 0.000142463479281955
5 2.74366574975223e-05
6 1.44975700969233e-05
7 3.73748485858245e-06
8 2.40757494835904e-06
9 1.34040458242958e-06
10 1.47604939852668e-07
11 8.49109837914781e-08
12 5.91845957000327e-08
13 1.74224683009756e-08
14 1.51054092822489e-08
15 6.07435279360086e-09
16 7.16579464736882e-09
17 3.01910092646183e-09
18 6.27830508569281e-10
19 6.15982655829241e-10
20 2.83140648375761e-10
21 3.03192480175683e-10
22 1.57506587196689e-10
23 1.36631157416829e-10
24 4.91399776836732e-11
25 9.63959118288821e-12
26 7.98281605427662e-12
27 6.6676277667823e-12
28 2.91438313811172e-12
29 1.90614931246576e-12
30 1.26634505351022e-12
31 1.10023595172808e-12
32 6.60646845871169e-13
33 7.2149035717033e-13
34 3.11375616648648e-13
35 1.74109286658845e-13
36 1.10318338826161e-13
37 7.70544122335398e-14
38 5.03169522508634e-14
39 3.8106143754838e-14
40 3.02112244686146e-14
41 1.85267439250041e-14
42 9.68443432443396e-15
43 6.52646660994462e-15
44 6.84226338139356e-15
45 4.63168598125102e-15
};
\addlegendentry{$\eta$}
\addplot [thick, dashed, color1]
table {%
0 1
1 0.00932238039336628
2 0.000266943760459508
3 4.90949857038112e-05
4 5.599960695831e-06
5 6.12824410476193e-07
6 3.86211523959076e-07
7 9.04010167944746e-08
8 5.82744632540993e-08
9 1.70168470184833e-08
10 7.04436726456646e-09
11 3.36213438130085e-09
12 3.11207201674506e-09
13 8.38932040992142e-10
14 3.78263928821779e-10
15 2.82163428907021e-10
16 6.95533441664873e-11
17 5.24268227949716e-11
18 4.19121093616782e-11
19 1.54370265480439e-11
20 8.64749152149743e-12
21 6.24447738563931e-12
22 5.44681941832098e-12
23 4.3396799607452e-12
24 2.47660431675892e-12
25 2.18757513175438e-13
26 1.86329799904232e-13
27 1.57026268596924e-13
28 9.55790405605825e-14
29 7.48084623843346e-14
30 7.36848530815293e-14
31 6.64897673464922e-14
32 2.24123930998593e-14
33 8.85932301724094e-15
34 1.03441814307043e-14
35 7.2349477058683e-15
36 4.64890234819224e-15
37 2.2322703665489e-15
38 2.07282248322398e-15
39 1.91337459989906e-15
};
\addlegendentry{$\tMomPlain$}
\end{axis}

\end{tikzpicture}
        }
    }
    \caption{Valve-like test case: results from the construction of the \gls{rom}.}
    \label{fig:valve-likeOffline}
\end{figure}

Next, we will construct the \gls{rom}. To create the snapshots, we generate $\nTrain = 256$ training samples that are equidistantly spaced over the parameter domain. These snapshots are used both for the \gls{pod} and the \gls{eim}. Applying the \gls{pod} results in the distribution of eigenvalues $\lambda_i$ shown in \Cref{fig:valve-likeEigenvalues}. Based on a threshold for the so-called retained energy, we choose $\nBasisVelocityHomROM = 1,\dots,10$ and $\nBasisPressureROM = 1,\dots,6$ in the following. Since all Dirichlet boundary conditions are parameter-independent in this case, we use only one lifting function resulting in $\nBasisVelocityROM = 1,\dots,11$. For comparison, the \gls{fom} includes a total number of \glspl{dof} of $\nBasisFOM = 334,762$. As has been mentioned before, the non-linearity in the viscosity $\visc$ and in the stabilization parameter $\tMomPlain$ is tackled by the \gls{eim}. The maximum interpolation error for both quantities in the course of the greedy search is depicted in \Cref{fig:valve-likeEIMErrors}. Setting a tolerance of $\SI{5e-15}{}$ and $\SI{2e-15}{}$, we use $Q_{\visc} = 45$ and $Q_{\tau} = 39$ terms for the approximation.
\bigskip
\par
\begin{figure}
    \captionsetup[sub]{position=bottom}
    \centering
    \Large
    \subcaptionbox{Velocity.}{
        \resizebox{0.44\textwidth}{!}{
\begin{tikzpicture}

\definecolor{color0}{rgb}{1,0.901960784313726,0.0980392156862745}
\definecolor{color1}{rgb}{1,0.701960784313725,0.298039215686275}
\definecolor{color2}{rgb}{1,0.498039215686275,0.501960784313725}
\definecolor{color3}{rgb}{1,0.298039215686275,0.701960784313725}
\definecolor{color4}{rgb}{1,0.0980392156862745,0.901960784313726}

\begin{axis}[
tick align=outside,
tick pos=left,
title={$\relErrorVelocity$},
x grid style={white!69.0196078431373!black},
xlabel={\(\displaystyle N_p\)},
xmin=1, xmax=6,
xtick style={color=black},
xtick={1,2,3,4,5,6},
xticklabels={
  \(\displaystyle 1\),
  \(\displaystyle 2\),
  \(\displaystyle 3\),
  \(\displaystyle 4\),
  \(\displaystyle 5\),
  \(\displaystyle 6\)
},
y grid style={white!69.0196078431373!black},
ylabel={\(\displaystyle N_{\mathbf{u}}\)},
ymin=2, ymax=11,
ytick style={color=black},
ytick={2,3,4,5,6,7,8,9,10,11},
yticklabels={
  \(\displaystyle 2\),
  \(\displaystyle 3\),
  \(\displaystyle 4\),
  \(\displaystyle 5\),
  \(\displaystyle 6\),
  \(\displaystyle 7\),
  \(\displaystyle 8\),
  \(\displaystyle 9\),
  \(\displaystyle 10\),
  \(\displaystyle 11\)
}
]
\addplot [draw=none, fill=color0]
table{%
x  y
5 4.96599787446714
6 4.94558403088085
6 5
6 6
6 7
6 8
6 9
6 10
6 11
5 11
4.45150563015397 11
4.44439519224548 10
4.39739935034553 9
4.36824322490229 8
4.22374978864749 7
4.20914210167462 6
4.53690125971689 5
5 4.96599787446714
};
\addplot [draw=none, fill=color1]
table{%
x  y
4 3.88534108136274
5 3.883984583036
6 3.88463500185776
6 4
6 4.94558403088085
5 4.96599787446714
4.53690125971689 5
4.20914210167462 6
4.22374978864749 7
4.36824322490229 8
4.39739935034553 9
4.44439519224548 10
4.45150563015397 11
4 11
3.79209881156784 11
3.7920675208487 10
3.76903747323132 9
3.76391793282701 8
3.76005659547684 7
3.70941213226149 6
3.69604662399699 5
3.72869415494111 4
4 3.88534108136274
};
\addplot [draw=none, fill=color2]
table{%
x  y
2 2.93124679617498
3 2.74974998839047
4 2.73430167488701
5 2.73596752053407
6 2.73665469600107
6 3
6 3.88463500185776
5 3.883984583036
4 3.88534108136274
3.72869415494111 4
3.69604662399699 5
3.70941213226149 6
3.76005659547684 7
3.76391793282701 8
3.76903747323132 9
3.7920675208487 10
3.79209881156784 11
3 11
2.22128768144098 11
2.20846393838137 10
2.18801105582443 9
2.15495112430777 8
2.14262579296783 7
2.12827313741416 6
2 5.37581418374384
1.99920268694824 5
1.99894567142086 4
1.99867460329951 3
2 2.93124679617498
};
\addplot [draw=none, fill=color3]
table{%
x  y
2 2
3 2
4 2
5 2
6 2
6 2.73665469600107
5 2.73596752053407
4 2.73430167488701
3 2.74974998839047
2 2.93124679617498
1.99867460329951 3
1.99894567142086 4
1.99920268694824 5
2 5.37581418374384
2.12827313741416 6
2.14262579296783 7
2.15495112430777 8
2.18801105582443 9
2.20846393838137 10
2.22128768144098 11
2 11
1.89268464986793 11
1.89260074197068 10
1.88972363516085 9
1.88904143611823 8
1.88914080571655 7
1.88263020803887 6
1.88085188940186 5
1.8801798059559 4
1.86130344507112 3
1.8749030879246 2
2 2
};
\addplot [draw=none, fill=color4]
table{%
x  y
1.86130344507112 3
1.8801798059559 4
1.88085188940186 5
1.88263020803887 6
1.88914080571655 7
1.88904143611823 8
1.88972363516085 9
1.89260074197068 10
1.89268464986793 11
1 11
1 10
1 9
1 8
1 7
1 6
1 5
1 4
1 3
1 2
1.8749030879246 2
1.86130344507112 3
};
\path [draw=black, semithick]
(axis cs:6,4.94558403088085)
--(axis cs:5,4.96599787446714)
--(axis cs:4.53690125971689,5)
--(axis cs:4.20914210167462,6)
--(axis cs:4.21459959660165,6.37360431786082);

\path [draw=black, semithick]
(axis cs:4.30925730464647,7.59177439623068)
--(axis cs:4.36824322490229,8)
--(axis cs:4.39739935034553,9)
--(axis cs:4.44439519224548,10)
--(axis cs:4.45150563015397,11);

\path [draw=black, semithick]
(axis cs:6,3.88463500185776)
--(axis cs:5,3.883984583036)
--(axis cs:4,3.88534108136274)
--(axis cs:3.72869415494111,4)
--(axis cs:3.69604662399699,5)
--(axis cs:3.70103918778138,5.37354088491024);

\path [draw=black, semithick]
(axis cs:3.74092044627085,6.6221472597187)
--(axis cs:3.76005659547684,7)
--(axis cs:3.76391793282701,8)
--(axis cs:3.76903747323132,9)
--(axis cs:3.7920675208487,10)
--(axis cs:3.79209881156784,11);

\path [draw=black, semithick]
(axis cs:6,2.73665469600107)
--(axis cs:5,2.73596752053407)
--(axis cs:4,2.73430167488701)
--(axis cs:3,2.74974998839047)
--(axis cs:2,2.93124679617498)
--(axis cs:1.99867460329951,3)
--(axis cs:1.99894567142086,4)
--(axis cs:1.9990415933667,4.3732145945417);

\path [draw=black, semithick]
(axis cs:2.04619596162303,5.60060687685705)
--(axis cs:2.12827313741416,6)
--(axis cs:2.14262579296783,7)
--(axis cs:2.15495112430777,8)
--(axis cs:2.18801105582443,9)
--(axis cs:2.20846393838137,10)
--(axis cs:2.22128768144098,11);

\path [draw=black, semithick]
(axis cs:1.8749030879246,2)
--(axis cs:1.86130344507112,3)
--(axis cs:1.8801798059559,4)
--(axis cs:1.88085188940186,5)
--(axis cs:1.88263020803887,6)
--(axis cs:1.88506056188794,6.37329197247242);

\path [draw=black, semithick]
(axis cs:1.88907852229237,7.62678550819736)
--(axis cs:1.88904143611823,8)
--(axis cs:1.88972363516085,9)
--(axis cs:1.89260074197068,10)
--(axis cs:1.89268464986793,11);

\draw (axis cs:4.22374978864749,7) node[
  scale=0.8,
  text=black,
  rotate=79.4
]{$10^{-6}$};
\draw (axis cs:3.70941213226149,6) node[
  scale=0.8,
  text=black,
  rotate=85.6
]{$10^{-5}$};
\draw (axis cs:1.99920268694824,5) node[
  scale=0.8,
  text=black,
  rotate=86.6
]{$10^{-4}$};
\draw (axis cs:1.88914080571655,7) node[
  scale=0.8,
  text=black,
  rotate=89.6
]{$10^{-3}$};
\end{axis}

\end{tikzpicture}}
    }
    \subcaptionbox{Pressure.}{
        \resizebox{0.44\textwidth}{!}{
\begin{tikzpicture}

\definecolor{color0}{rgb}{1,0.917647058823529,0.0823529411764706}
\definecolor{color1}{rgb}{1,0.749019607843137,0.250980392156863}
\definecolor{color2}{rgb}{1,0.584313725490196,0.415686274509804}
\definecolor{color3}{rgb}{1,0.415686274509804,0.584313725490196}
\definecolor{color4}{rgb}{1,0.250980392156863,0.749019607843137}
\definecolor{color5}{rgb}{1,0.0823529411764706,0.917647058823529}

\begin{axis}[
tick align=outside,
tick pos=left,
title={$\relErrorPressure$},
x grid style={white!69.0196078431373!black},
xlabel={\(\displaystyle N_p\)},
xmin=1, xmax=6,
xtick style={color=black},
xtick={1,2,3,4,5,6},
xticklabels={
  \(\displaystyle 1\),
  \(\displaystyle 2\),
  \(\displaystyle 3\),
  \(\displaystyle 4\),
  \(\displaystyle 5\),
  \(\displaystyle 6\)
},
y grid style={white!69.0196078431373!black},
ylabel={\(\displaystyle N_{\mathbf{u}}\)},
ymin=2, ymax=11,
ytick style={color=black},
ytick={2,3,4,5,6,7,8,9,10,11},
yticklabels={
  \(\displaystyle 2\),
  \(\displaystyle 3\),
  \(\displaystyle 4\),
  \(\displaystyle 5\),
  \(\displaystyle 6\),
  \(\displaystyle 7\),
  \(\displaystyle 8\),
  \(\displaystyle 9\),
  \(\displaystyle 10\),
  \(\displaystyle 11\)
}
]
\addplot [draw=none, fill=color0]
table{%
x  y
6 3.99175330071759
6 4
6 5
6 6
6 7
6 8
6 9
6 10
6 11
5.46729690810002 11
5.4687824974526 10
5.47061512667462 9
5.46648881531239 8
5.462665299272 7
5.52745867768282 6
5.47368849457939 5
5.78987892924434 4
6 3.99175330071759
};
\addplot [draw=none, fill=color1]
table{%
x  y
4 3.92571089090235
5 3.77144437107509
6 3.72000400255733
6 3.99175330071759
5.78987892924434 4
5.47368849457939 5
5.52745867768282 6
5.462665299272 7
5.46648881531239 8
5.47061512667462 9
5.4687824974526 10
5.46729690810002 11
5 11
4 11
3.98968441585694 11
3.9897173028144 10
3.98931863181225 9
3.98920930286623 8
3.98916868509722 7
3.98912255106139 6
3.98904468895928 5
3.98947162075893 4
4 3.92571089090235
};
\addplot [draw=none, fill=color2]
table{%
x  y
4 2.30445672025347
5 2.33250776366727
6 2.432164097066
6 3
6 3.72000400255733
5 3.77144437107509
4 3.92571089090235
3.98947162075893 4
3.98904468895928 5
3.98912255106139 6
3.98916868509722 7
3.98920930286623 8
3.98931863181225 9
3.9897173028144 10
3.98968441585694 11
3.34814960283007 11
3.34857374765808 10
3.34125290185621 9
3.34333442502001 8
3.34329408616309 7
3.33351737131764 6
3.33473861614855 5
3.33513260167389 4
3.30994634836988 3
4 2.30445672025347
};
\addplot [draw=none, fill=color3]
table{%
x  y
2 2
3 2
4 2
5 2
6 2
6 2.432164097066
5 2.33250776366727
4 2.30445672025347
3.30994634836988 3
3.33513260167389 4
3.33473861614855 5
3.33351737131764 6
3.34329408616309 7
3.34333442502001 8
3.34125290185621 9
3.34857374765808 10
3.34814960283007 11
3 11
2 11
1.97891596034193 11
1.97895047703883 10
1.97895983627377 9
1.97900237407811 8
1.97899390040099 7
1.9789633073209 6
1.97896592073396 5
1.97895199710449 4
1.97893120263934 3
1.98082715218374 2
2 2
};
\addplot [draw=none, fill=color4]
table{%
x  y
1.97893120263934 3
1.97895199710449 4
1.97896592073396 5
1.9789633073209 6
1.97899390040099 7
1.97900237407811 8
1.97895983627377 9
1.97895047703883 10
1.97891596034193 11
1.63837733657086 11
1.63809841851355 10
1.63788986849578 9
1.63801429212607 8
1.63802543214768 7
1.63723278139985 6
1.63723507236599 5
1.63724318905758 4
1.63705647212741 3
1.63833996554076 2
1.98082715218374 2
1.97893120263934 3
};
\addplot [draw=none, fill=color5]
table{%
x  y
1.63705647212741 3
1.63724318905758 4
1.63723507236599 5
1.63723278139985 6
1.63802543214768 7
1.63801429212607 8
1.63788986849578 9
1.63809841851355 10
1.63837733657086 11
1 11
1 10
1 9
1 8
1 7
1 6
1 5
1 4
1 3
1 2
1.63833996554076 2
1.63705647212741 3
};
\path [draw=black, semithick]
(axis cs:6,3.99175330071759)
--(axis cs:5.78987892924434,4)
--(axis cs:5.47368849457939,5)
--(axis cs:5.52745867768282,6)
--(axis cs:5.50278838494497,6.38075330138575);

\path [draw=black, semithick]
(axis cs:5.46506172157529,7.62675879425428)
--(axis cs:5.46648881531239,8)
--(axis cs:5.47061512667462,9)
--(axis cs:5.4687824974526,10)
--(axis cs:5.46729690810002,11);

\path [draw=black, semithick]
(axis cs:6,3.72000400255733)
--(axis cs:5,3.77144437107509)
--(axis cs:4,3.92571089090235)
--(axis cs:3.98947162075893,4)
--(axis cs:3.98904468895928,5)
--(axis cs:3.98907374822361,5.37321448483213);

\path [draw=black, semithick]
(axis cs:3.98915146720713,6.62678552236211)
--(axis cs:3.98916868509722,7)
--(axis cs:3.98920930286623,8)
--(axis cs:3.98931863181225,9)
--(axis cs:3.9897173028144,10)
--(axis cs:3.98968441585694,11);

\path [draw=black, semithick]
(axis cs:6,2.432164097066)
--(axis cs:5,2.33250776366727)
--(axis cs:4,2.30445672025347)
--(axis cs:3.30994634836988,3)
--(axis cs:3.33513260167389,4)
--(axis cs:3.33473861614855,5)
--(axis cs:3.33428282657094,5.37321720106552);

\path [draw=black, semithick]
(axis cs:3.33964356650891,6.62661080824277)
--(axis cs:3.34329408616309,7)
--(axis cs:3.34333442502001,8)
--(axis cs:3.34125290185621,9)
--(axis cs:3.34857374765808,10)
--(axis cs:3.34814960283007,11);

\path [draw=black, semithick]
(axis cs:1.98082715218374,2)
--(axis cs:1.97893120263934,3)
--(axis cs:1.97895199710449,4)
--(axis cs:1.97896592073396,5)
--(axis cs:1.9789633073209,6)
--(axis cs:1.97897472510124,6.37321447545737);

\path [draw=black, semithick]
(axis cs:1.97899921157916,7.62678552612284)
--(axis cs:1.97900237407811,8)
--(axis cs:1.97895983627377,9)
--(axis cs:1.97895047703883,10)
--(axis cs:1.97891596034193,11);

\path [draw=black, semithick]
(axis cs:1.63833996554076,2)
--(axis cs:1.63705647212741,3)
--(axis cs:1.63724318905758,4)
--(axis cs:1.63723507236599,5)
--(axis cs:1.63723278139985,6)
--(axis cs:1.63752861104227,6.37321562268217);

\path [draw=black, semithick]
(axis cs:1.63801844974337,7.62678552602721)
--(axis cs:1.63801429212607,8)
--(axis cs:1.63788986849578,9)
--(axis cs:1.63809841851355,10)
--(axis cs:1.63837733657086,11);

\draw (axis cs:5.462665299272,7) node[
  scale=0.8,
  text=black,
  rotate=274.2
]{$10^{-7}$};
\draw (axis cs:3.98912255106139,6) node[
  scale=0.8,
  text=black,
  rotate=90.0
]{$10^{-6}$};
\draw (axis cs:3.33351737131764,6) node[
  scale=0.8,
  text=black,
  rotate=89.4
]{$10^{-5}$};
\draw (axis cs:1.97899390040099,7) node[
  scale=0.8,
  text=black,
  rotate=90.0
]{$10^{-4}$};
\draw (axis cs:1.63802543214768,7) node[
  scale=0.8,
  text=black,
  rotate=89.9
]{$10^{-3}$};
\end{axis}

\end{tikzpicture}}
    }
    \caption{Valve-like test case: maximum relative error of the \gls{rom} over all testing samples.}
    \label{fig:valveRelativeErrorMax}
\end{figure}
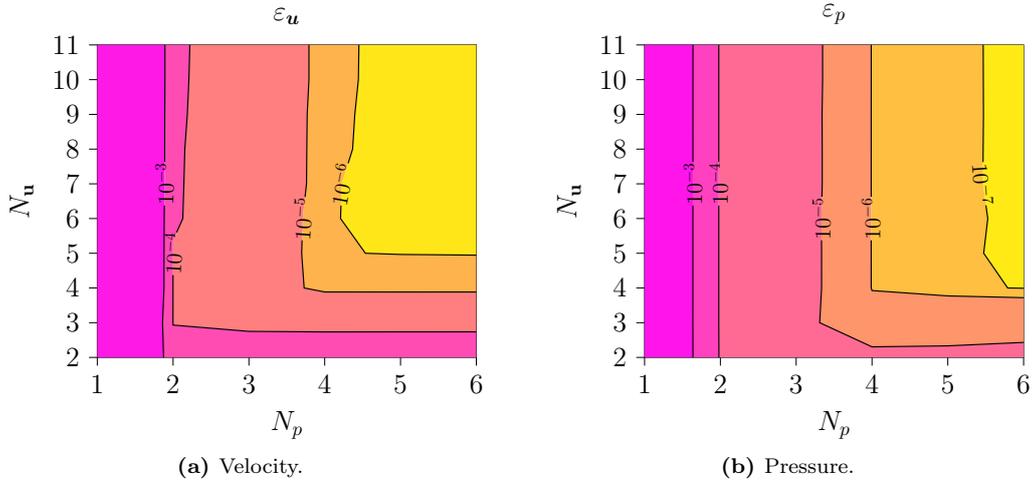

To quantify the quality of the \gls{rom}, we perform an error and performance analysis. To that end, we create $\nTest = 50$ testing samples drawn from a uniform distribution. For each sample, we compare the solutions from the \gls{fom} and the \gls{rom} as well as the respective runtimes.
To evaluate the accuracy of the \gls{rom}, we use the following error definitions:
\begin{align}
    \relErrorVelocity
    =
    \frac
    {|\trialVelocityDiscrete-\trialVelocityReduced|_{\sobolevSpaceVector}}
    {|\trialVelocityDiscrete|_{\sobolevSpaceVector}},
    \,
    \relErrorPressure
    =
    \frac
    {||\trialPressureDiscrete-\trialPressureReduced||_{L^2}}
    {||\trialPressureDiscrete||_{L^2}},
\end{align}
where $|\cdot|_{\sobolevSpaceVector}$ and $||\cdot||_{L^2}$ are discrete measures for the $\sobolevSpaceVector$ semi-norm and the $L^2$ norm, respectively. The maximal error over all testing samples is shown in \Cref{fig:valveRelativeErrorMax}, using different numbers of basis functions for the velocity and the pressure field. Note that we ignore the lifting function in these plots, since it does not represent a solution, but is only intended to ensure the Dirichlet boundary conditions.
For the velocity, the error ranges from values smaller than \SI{1e-2}{} to values smaller than \SI{1e-6}{} when increasing the number of basis functions. Similarly, the maximum error in the pressure is limited by \SI{1e-2}{} and drops below \SI{1e-7}{} for $\nBasisVelocityROM$ and $\nBasisPressureROM$ large enough. These results suggest that the error introduced by the \gls{rom} is in a reasonable range for typical engineering applications and, furthermore, that it can be controlled by choosing the number of basis functions such that a desired accuracy is achieved.
\par
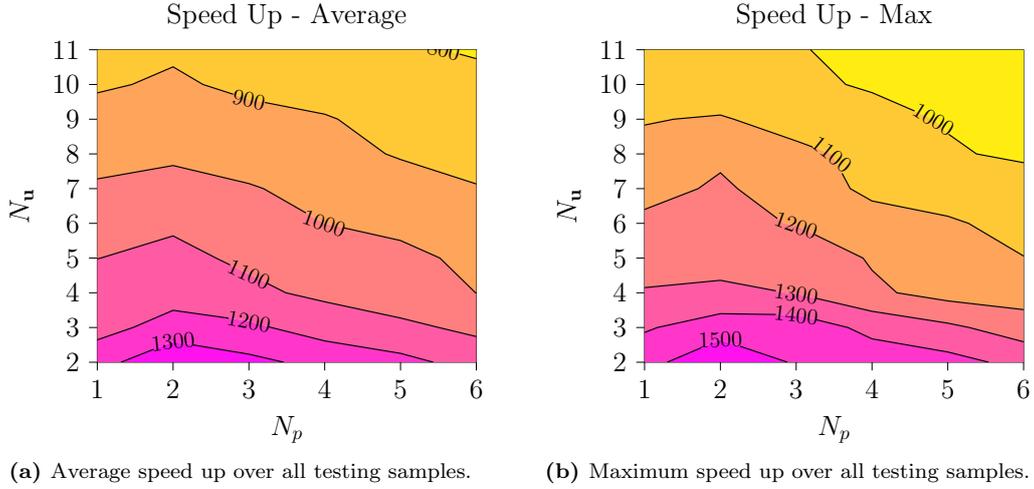
\begin{figure}
    \captionsetup[sub]{position=bottom}
    \centering
    \Large
    \subcaptionbox{Average speed up over all testing samples.}{
        \resizebox{0.44\textwidth}{!}{
\begin{tikzpicture}

\definecolor{color0}{rgb}{1,0.929411764705882,0.0705882352941176}
\definecolor{color1}{rgb}{1,0.788235294117647,0.211764705882353}
\definecolor{color2}{rgb}{1,0.643137254901961,0.356862745098039}
\definecolor{color3}{rgb}{1,0.498039215686275,0.501960784313725}
\definecolor{color4}{rgb}{1,0.356862745098039,0.643137254901961}
\definecolor{color5}{rgb}{1,0.211764705882353,0.788235294117647}
\definecolor{color6}{rgb}{1,0.0705882352941176,0.929411764705882}

\begin{axis}[
tick align=outside,
tick pos=left,
title={Speed Up - Average},
x grid style ={white!69.0196078431373!black},
xlabel={\(\displaystyle N_p\)},
xmin=1, xmax=6,
xtick style={color=black},
xtick={1,2,3,4,5,6},
xticklabels={
  \(\displaystyle 1\),
  \(\displaystyle 2\),
  \(\displaystyle 3\),
  \(\displaystyle 4\),
  \(\displaystyle 5\),
  \(\displaystyle 6\)
},
y grid style={white!69.0196078431373!black},
ylabel={\(\displaystyle N_{\mathbf{u}}\)},
ymin=2, ymax=11,
ytick style={color=black},
ytick={2,3,4,5,6,7,8,9,10,11},
yticklabels={
  \(\displaystyle 2\),
  \(\displaystyle 3\),
  \(\displaystyle 4\),
  \(\displaystyle 5\),
  \(\displaystyle 6\),
  \(\displaystyle 7\),
  \(\displaystyle 8\),
  \(\displaystyle 9\),
  \(\displaystyle 10\),
  \(\displaystyle 11\)
}
]
\addplot [draw=none, fill=color0]
table{%
x  y
6 10.7426292953028
6 11
5.57031648131517 11
6 10.7426292953028
};
\addplot [draw=none, fill=color1]
table{%
x  y
5 7.83894797094878
6 7.13509322303163
6 8
6 9
6 10
6 10.7426292953028
5.57031648131517 11
5 11
4 11
3 11
2 11
1 11
1 10
1 9.75767609564334
1.45734684097972 10
2 10.5078469602098
2.39394636086238 10
3 9.57223390198292
4 9.14460119922037
4.17406044843618 9
4.80284005228388 8
5 7.83894797094878
};
\addplot [draw=none, fill=color2]
table{%
x  y
6 3.99479159633512
6 4
6 5
6 6
6 7
6 7.13509322303163
5 7.83894797094878
4.80284005228388 8
4.17406044843618 9
4 9.14460119922037
3 9.57223390198292
2.39394636086238 10
2 10.5078469602098
1.45734684097972 10
1 9.75767609564334
1 9
1 8
1 7.28016385651156
2 7.66532012449406
3 7.14668088289607
3.19339138167989 7
3.99244708274187 6
4 5.99487469286416
5 5.5057151600678
5.51991847795867 5
5.98951589728088 4
6 3.99479159633512
};
\addplot [draw=none, fill=color3]
table{%
x  y
6 2.73993166459977
6 3
6 3.99479159633512
5.98951589728088 4
5.51991847795867 5
5 5.5057151600678
4 5.99487469286416
3.99244708274187 6
3.19339138167989 7
3 7.14668088289607
2 7.66532012449406
1 7.28016385651156
1 7
1 6
1 5
1 4.97488346272111
1.04462887398337 5
2 5.63867372454922
2.56658441191742 5
3 4.53904904661717
3.48697363850504 4
4 3.74124459044498
5 3.27370866683256
5.50365841864754 3
6 2.73993166459977
};
\addplot [draw=none, fill=color4]
table{%
x  y
1.47135757322684 3
2 3.4953795716924
3 3.18636001561186
3.29511815097591 3
4 2.61445841715845
5 2.25590772313131
5.43288407044525 2
6 2
6 2.73993166459977
5.50365841864754 3
5 3.27370866683256
4 3.74124459044498
3.48697363850504 4
3 4.53904904661717
2.56658441191742 5
2 5.63867372454922
1.04462887398337 5
1 4.97488346272111
1 4
1 3
1 2.6452570407776
1.47135757322684 3
};
\addplot [draw=none, fill=color5]
table{%
x  y
2 2.59435916397521
3 2.23167347223914
3.45984522193356 2
4 2
5 2
5.43288407044525 2
5 2.25590772313131
4 2.61445841715845
3.29511815097591 3
3 3.18636001561186
2 3.4953795716924
1.47135757322684 3
1 2.6452570407776
1 2
1.31309890311084 2
2 2.59435916397521
};
\addplot [draw=none, fill=color6]
table{%
x  y
2 2
3 2
3.45984522193356 2
3 2.23167347223914
2 2.59435916397521
1.31309890311084 2
2 2
};
\path [draw=black, semithick]
(axis cs:6,10.7426292953028)
--(axis cs:5.77831481102616,10.8754137071456);

\path [draw=black, semithick]
(axis cs:6,7.13509322303163)
--(axis cs:5,7.83894797094878)
--(axis cs:4.80284005228388,8)
--(axis cs:4.17406044843618,9)
--(axis cs:4,9.14460119922037)
--(axis cs:3.21101633736785,9.48199641530725);

\path [draw=black, semithick]
(axis cs:2.79430370332184,9.71741890958685)
--(axis cs:2.39394636086238,10)
--(axis cs:2,10.5078469602098)
--(axis cs:1.45734684097972,10)
--(axis cs:1,9.75767609564334);

\path [draw=black, semithick]
(axis cs:6,3.99479159633512)
--(axis cs:5.98951589728088,4)
--(axis cs:5.51991847795867,5)
--(axis cs:5,5.5057151600678)
--(axis cs:4.25588965990414,5.86970382637804);

\path [draw=black, semithick]
(axis cs:3.7536619573203,6.29883414273149)
--(axis cs:3.19339138167989,7)
--(axis cs:3,7.14668088289607)
--(axis cs:2,7.66532012449406)
--(axis cs:1,7.28016385651156);

\path [draw=black, semithick]
(axis cs:6,2.73993166459977)
--(axis cs:5.50365841864754,3)
--(axis cs:5,3.27370866683256)
--(axis cs:4,3.74124459044498)
--(axis cs:3.48697363850504,4)
--(axis cs:3.24448165327576,4.26842330490708);

\path [draw=black, semithick]
(axis cs:2.75387070400582,4.80081520538842)
--(axis cs:2.56658441191742,5)
--(axis cs:2,5.63867372454922)
--(axis cs:1.04462887398337,5)
--(axis cs:1,4.97488346272111);

\path [draw=black, semithick]
(axis cs:5.43288407044525,2)
--(axis cs:5,2.25590772313131)
--(axis cs:4,2.61445841715845)
--(axis cs:3.29511815097591,3)
--(axis cs:3.26018534483042,3.02205922704894);

\path [draw=black, semithick]
(axis cs:2.73324507418725,3.26879250436881)
--(axis cs:2,3.4953795716924)
--(axis cs:1.47135757322684,3)
--(axis cs:1,2.6452570407776);

\path [draw=black, semithick]
(axis cs:3.45984522193356,2)
--(axis cs:3,2.23167347223914)
--(axis cs:2.26594788992643,2.49790366955149);

\path [draw=black, semithick]
(axis cs:1.74682377931879,2.37529180836322)
--(axis cs:1.31309890311084,2);

\draw (axis cs:5.57031648131517,11) node[
  scale=0.8,
  text=black,
  rotate=346.1
]{800};
\draw (axis cs:3,9.57223390198292) node[
  scale=0.8,
  text=black,
  rotate=346.8
]{900};
\draw (axis cs:3.99244708274187,6) node[
  scale=0.8,
  text=black,
  rotate=340.5
]{1000};
\draw (axis cs:3,4.53904904661717) node[
  scale=0.8,
  text=black,
  rotate=335.8
]{1100};
\draw (axis cs:3,3.18636001561186) node[
  scale=0.8,
  text=black,
  rotate=349.0
]{1200};
\draw (axis cs:2,2.59435916397521) node[
  scale=0.8,
  text=black,
  rotate=5.6
]{1300};
\end{axis}

\end{tikzpicture}}
    }
    \subcaptionbox{Maximum speed up over all testing samples.}{
        \resizebox{0.44\textwidth}{!}{
\begin{tikzpicture}

\definecolor{color0}{rgb}{1,0.929411764705882,0.0705882352941176}
\definecolor{color1}{rgb}{1,0.788235294117647,0.211764705882353}
\definecolor{color2}{rgb}{1,0.643137254901961,0.356862745098039}
\definecolor{color3}{rgb}{1,0.498039215686275,0.501960784313725}
\definecolor{color4}{rgb}{1,0.356862745098039,0.643137254901961}
\definecolor{color5}{rgb}{1,0.211764705882353,0.788235294117647}
\definecolor{color6}{rgb}{1,0.0705882352941176,0.929411764705882}

\begin{axis}[
tick align=outside,
tick pos=left,
title={Speed Up - Max},
x grid style={white!69.0196078431373!black},
xlabel={\(\displaystyle N_p\)},
xmin=1, xmax=6,
xtick style={color=black},
xtick={1,2,3,4,5,6},
xticklabels={
  \(\displaystyle 1\),
  \(\displaystyle 2\),
  \(\displaystyle 3\),
  \(\displaystyle 4\),
  \(\displaystyle 5\),
  \(\displaystyle 6\)
},
y grid style={white!69.0196078431373!black},
ylabel={\(\displaystyle N_{\mathbf{u}}\)},
ymin=2, ymax=11,
ytick style={color=black},
ytick={2,3,4,5,6,7,8,9,10,11},
yticklabels={
  \(\displaystyle 2\),
  \(\displaystyle 3\),
  \(\displaystyle 4\),
  \(\displaystyle 5\),
  \(\displaystyle 6\),
  \(\displaystyle 7\),
  \(\displaystyle 8\),
  \(\displaystyle 9\),
  \(\displaystyle 10\),
  \(\displaystyle 11\)
}
]
\addplot [draw=none, fill=color0]
table{%
x  y
6 7.7479711514141
6 8
6 9
6 10
6 11
5 11
4 11
3.18907468301226 11
3.6542422633017 10
4 9.76546811290731
4.81305624881071 9
5 8.58317161530536
5.37825869829586 8
6 7.7479711514141
};
\addplot [draw=none, fill=color1]
table{%
x  y
6 5.05893301083366
6 6
6 7
6 7.7479711514141
5.37825869829586 8
5 8.58317161530536
4.81305624881071 9
4 9.76546811290731
3.6542422633017 10
3.18907468301226 11
3 11
2 11
1 11
1 10
1 9
1 8.83063460943087
1.38479402207256 9
2 9.11992187081993
2.17634453203126 9
3 8.37325226941018
3.46894123267777 8
3.71282666270448 7
4 6.64468357965588
5 6.20637290372433
5.27372208134803 6
6 5.05893301083366
};
\addplot [draw=none, fill=color2]
table{%
x  y
5 3.76620239625192
6 3.5152227366982
6 4
6 5
6 5.05893301083366
5.27372208134803 6
5 6.20637290372433
4 6.64468357965588
3.71282666270448 7
3.46894123267777 8
3 8.37325226941017
2.17634453203126 9
2 9.11992187081993
1.38479402207256 9
1 8.83063460943087
1 8
1 7
1 6.4006276576681
1.70789166385819 7
2 7.45653812498959
2.22605929964209 7
2.83667132641095 6
3 5.93475023455696
3.88150504652744 5
4 4.64913882141652
4.32762128598689 4
5 3.76620239625192
};
\addplot [draw=none, fill=color3]
table{%
x  y
6 2.58864568119381
6 3
6 3.5152227366982
5 3.76620239625192
4.32762128598689 4
4 4.64913882141652
3.88150504652744 5
3 5.93475023455696
2.83667132641095 6
2.22605929964209 7
2 7.45653812498959
1.70789166385819 7
1 6.4006276576681
1 6
1 5
1 4.1493491416269
2 4.35928013384988
2.8690525340667 4
3 3.97509391102203
4 3.46660600393282
5 3.12549197127818
5.2821679122398 3
6 2.58864568119381
};
\addplot [draw=none, fill=color4]
table{%
x  y
1.16844559164342 3
2 3.39818634001089
3 3.36615353640909
3.68170447217772 3
4 2.67387601038452
5 2.29505899140663
5.53977807598436 2
6 2
6 2.58864568119381
5.2821679122398 3
5 3.12549197127818
4 3.46660600393282
3 3.97509391102203
2.8690525340667 4
2 4.35928013384988
1 4.1493491416269
1 4
1 3
1 2.85772492129179
1.16844559164342 3
};
\addplot [draw=none, fill=color5]
table{%
x  y
2 2.65103963463062
2.88624927701573 2
3 2
4 2
5 2
5.53977807598436 2
5 2.29505899140663
4 2.67387601038452
3.68170447217772 3
3 3.36615353640909
2 3.39818634001089
1.16844559164341 3
1 2.85772492129179
1 2
1.29856987770096 2
2 2.65103963463062
};
\addplot [draw=none, fill=color6]
table{%
x  y
2 2
2.88624927701573 2
2 2.65103963463062
1.29856987770096 2
2 2
};
\path [draw=black, semithick]
(axis cs:6,7.7479711514141)
--(axis cs:5.37825869829586,8)
--(axis cs:5.01223893980926,8.56430251241269);

\path [draw=black, semithick]
(axis cs:4.56248612438337,9.23590426932565)
--(axis cs:4,9.76546811290731)
--(axis cs:3.6542422633017,10)
--(axis cs:3.18907468301226,11);

\path [draw=black, semithick]
(axis cs:6,5.05893301083366)
--(axis cs:5.27372208134803,6)
--(axis cs:5,6.20637290372433)
--(axis cs:4,6.64468357965588)
--(axis cs:3.71282666270448,7)
--(axis cs:3.60544684353425,7.44028796291141);

\path [draw=black, semithick]
(axis cs:3.21352108345597,8.20330084818034)
--(axis cs:3,8.37325226941017)
--(axis cs:2.17634453203126,9)
--(axis cs:2,9.11992187081993)
--(axis cs:1.38479402207256,9)
--(axis cs:1,8.83063460943087);

\path [draw=black, semithick]
(axis cs:6,3.5152227366982)
--(axis cs:5,3.76620239625192)
--(axis cs:4.32762128598689,4)
--(axis cs:4,4.64913882141652)
--(axis cs:3.88150504652744,5)
--(axis cs:3.24624676649357,5.67362952555419);

\path [draw=black, semithick]
(axis cs:2.75110180493382,6.14013730114346)
--(axis cs:2.22605929964209,7)
--(axis cs:2,7.45653812498959)
--(axis cs:1.70789166385819,7)
--(axis cs:1,6.4006276576681);

\path [draw=black, semithick]
(axis cs:6,2.58864568119381)
--(axis cs:5.2821679122398,3)
--(axis cs:5,3.12549197127818)
--(axis cs:4,3.46660600393282)
--(axis cs:3.26316149966656,3.84127947083013);

\path [draw=black, semithick]
(axis cs:2.7334480862735,4.05606103456807)
--(axis cs:2,4.35928013384988)
--(axis cs:1,4.1493491416269);

\path [draw=black, semithick]
(axis cs:5.53977807598436,2)
--(axis cs:5,2.29505899140663)
--(axis cs:4,2.67387601038452)
--(axis cs:3.68170447217772,3)
--(axis cs:3.26251753533931,3.2251514924782);

\path [draw=black, semithick]
(axis cs:2.73109479512054,3.3747673240245)
--(axis cs:2,3.39818634001089)
--(axis cs:1.16844559164341,3)
--(axis cs:1,2.85772492129179);

\path [draw=black, semithick]
(axis cs:2.88624927701573,2)
--(axis cs:2.2572943121112,2.46203096709922);

\path [draw=black, semithick]
(axis cs:1.74896474965724,2.41803866637041)
--(axis cs:1.29856987770096,2);

\draw (axis cs:4.81305624881071,9) node[
  scale=0.8,
  text=black,
  rotate=328.0
]{1000};
\draw (axis cs:3.46894123267777,8) node[
  scale=0.8,
  text=black,
  rotate=321.1
]{1100};
\draw (axis cs:3,5.93475023455696) node[
  scale=0.8,
  text=black,
  rotate=340.5
]{1200};
\draw (axis cs:3,3.97509391102203) node[
  scale=0.8,
  text=black,
  rotate=350.8
]{1300};
\draw (axis cs:3,3.36615353640909) node[
  scale=0.8,
  text=black,
  rotate=353.4
]{1400};
\draw (axis cs:2,2.65103963463062) node[
  scale=0.8,
  text=black,
  rotate=2.1
]{1500};
\end{axis}

\end{tikzpicture}}
    }
    \caption{Valve-like test case: performance results.}
    \label{fig:valvePerformance}
\end{figure}

Finally, we investigate the performance of the \gls{rom} by comparing the CPU time needed for a single evaluation of the \gls{rom} to that for the \gls{fom}. Note that the former has been run on a single core, whereas the latter used 64 cores. Also in this analysis, we skip the lifting function for the same reason as before. \Cref{fig:valvePerformance} presents the speed up, i.e., the ratio between the CPU time of \gls{fom} and \gls{rom} evaluations, for different number of basis functions $\nBasisVelocityROM$ and $\nBasisPressureROM$.  Both, the average and the maximum speed up, are in the order of 1000 and indicate that a significant reduction of the runtime for the \gls{rom} is realized. 
\par
Taken together, the error and performance analysis confirms the effectiveness of the approach for this two-dimensional deforming domain problem with topology changes. In particular, the accuracy of the \gls{rom} is acceptable as well as controllable while a significant reduction of the  computational demands is achieved, too. This qualifies the \gls{rom} as a surrogate model, e.g., in one of the aforementioned many query scenarios.   
\subsection{Artery-Like Geometry with Compression}
After presenting results for a spatially deforming two-dimensional geometry, we will now demonstrate the aptitude of the proposed approach also for the three-dimensional case, resulting in a four-dimensional space-time domain. The geometry is inspired by an artery that locally undergoes compression over time. The initial spatial geometry can be seen in \Cref{fig:arteryGeoInitialFront,fig:arteryGeoInitialSide,fig:arteryGeoInitialTop}. It has a length of $L = \SI{60e-3}{\meter}$ and a radius of $r_0 = \SI{5e-3}{\meter}$.
We are interested in the internal flow over a time period of $\SI{1}{\second}$ where fluid enters on the left-hand side, i.e., at $x_{\text{min}}$ .
The local narrowing of the artery happens according to the following expression for the upper and lower parts of the moving boundary:
\begin{align*}
    y(t) = \pm  \left[ 0.2 + 0.2 \cdot\left(\cos{(\pi\, t\,\si{\per\second})}+1\right)\right] \cdot r_0 .
\end{align*}
The final state of the deformed spatial geometry is depicted in \Cref{fig:arteryGeoFinalSide,fig:arteryGeoFinalTop}.
\par
To mimic blood, the density is set to $\density = \SI{1058}{\kilo\gram\per\cubic\meter}$ and the parameters for the viscosity model are chosen as presented in \cite{Cho1991}:
\begin{align*}
    \visc_0 = \SI{0.056}{Pa\cdot s}, \,
    \visc_{\infty} = \SI{0.00345}{Pa\cdot s}, \,
    \lambda = \SI{1.902}{s}, \,
    a = \SI{1.25}{}, \,
    n = \SI{0.22}{}.
\end{align*}
\par
For the inflow velocity, we prescribe the following time-dependent profile:
\begin{align*}
    \begin{aligned}
        \velX = \velX_{\text{in}}^0 \lp 1 - \frac{(y^2+z^2)}{r_0^2} \rp  \cdot 
        \begin{cases}  
            \sqrt{\frac{t}{\SI{0.2}{\second}}}    &\text{ for } t < \SI{0.2}{\second},
            \\
            1   &\text{ for } t \ge \SI{0.2}{\second},
        \end{cases}
    \end{aligned}
    \quad
    \begin{aligned}
        \velY = \SI{0}{\meter\per\second}, \quad \velZ = \SI{0}{\meter\per\second},       
    \end{aligned}
    \label{eq:artery-likeInflow}
\end{align*}
with the velocity vector $\trialVelocity = \lb \velX, \velY, \velZ \rb^T$ and $\velX_{\text{in}}^0=\SI{0.1}{\meter\per\second}$. At the outlet, a parallel outflow is enforced, i.e., $\velY=\velZ=0$. Along the walls, no-slip conditions are set. To account for the narrowing, we apply the following boundary conditions for the velocity in $y$-direction on the horizontal and rounded parts:
\begin{align}
   \velY = \frac{\partial y(t)}{\partial t} = \mp \pi \sin(\pi t)\times 10^{-3}\,\si{\meter\per\second}.
\end{align}
The sign of this term depends on whether you consider the upper or the lower part. In particular, the negative and the positive sign correspond to a downward movement for $y>0$ and an upward movement for $y<0$, respectively. 
For this geometry, a locally refined boundary-conforming simplex space-time mesh is constructed using a four-dimensional elastic mesh update~\cite{Danwitz2021}.
\Cref{fig:artery-likeVelocityInitial,fig:artery-likeVelocityFinal} show the resulting velocity field along the artery in its center plane for the initial and final state.
\newcommand{\myW}{0.48\textwidth} 
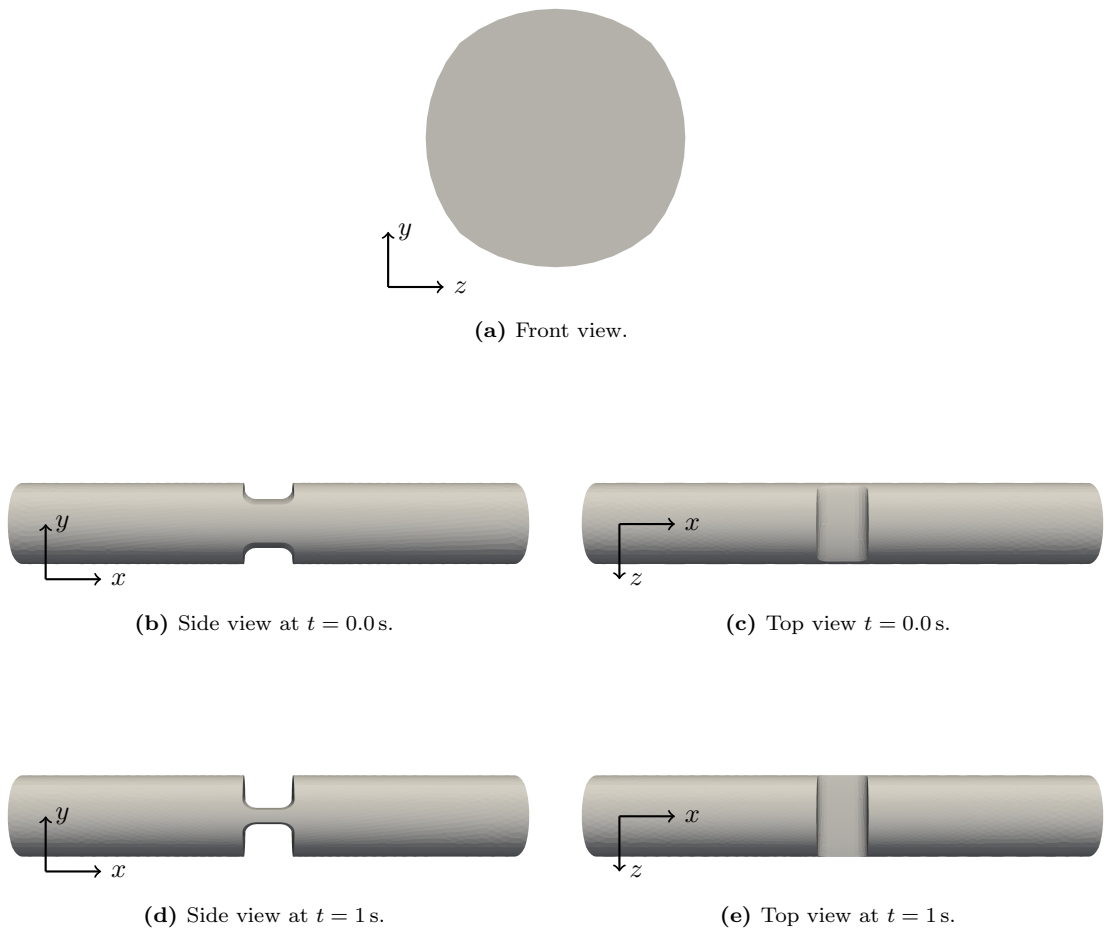
\begin{figure}
\captionsetup[sub]{position=bottom}
\centering
\subcaptionbox{Front view.\label{fig:arteryGeoInitialFront}}{
   \begin{tikzpicture}[
    axis/.style={thick, ->},
    every node/.style={color=black}
    ]
       \node at (0,0) [inner sep=0pt, anchor=south west] {\includegraphics[width=\myW,trim={0cm 0cm 0cm 0cm},clip]{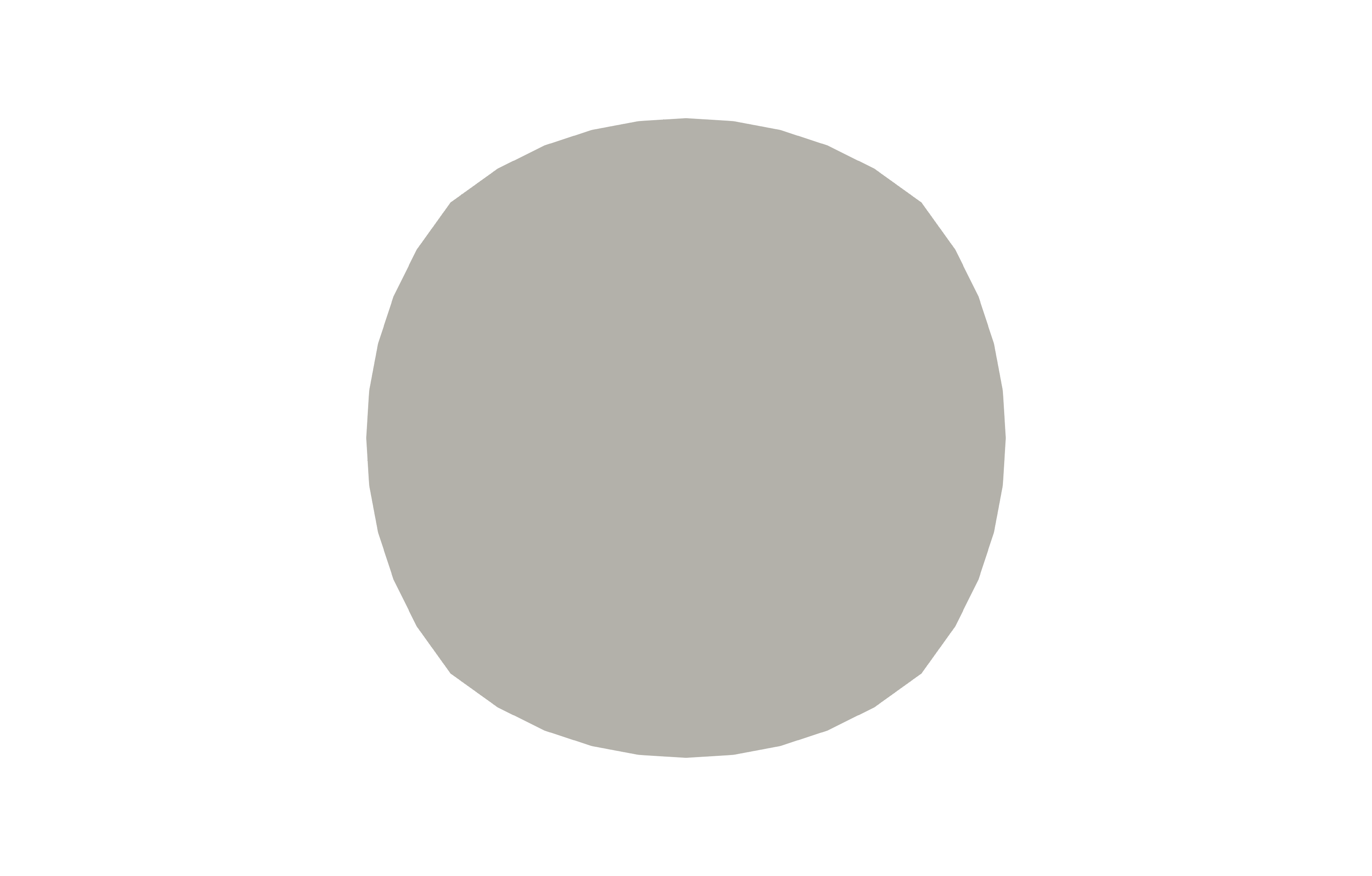}};
       \draw[axis] (0.2*\myW,0.05*\myW)  -- (0.3*\myW,0.05*\myW) node(xline)[right]{$z$};
       \draw[axis] (0.2*\myW,0.05*\myW)  -- (0.2*\myW,0.15*\myW) node(xline)[right]{$y$};
\end{tikzpicture}
}
\subcaptionbox{Side view at $t=0.0\,$s.\label{fig:arteryGeoInitialSide}}{
\begin{tikzpicture}[
    axis/.style={thick, ->},
    every node/.style={color=black}
    ]
       \node at (0,0) [inner sep=0pt, anchor=south west] {\includegraphics[width=\myW,trim={0cm 20cm 0cm 0cm},clip]{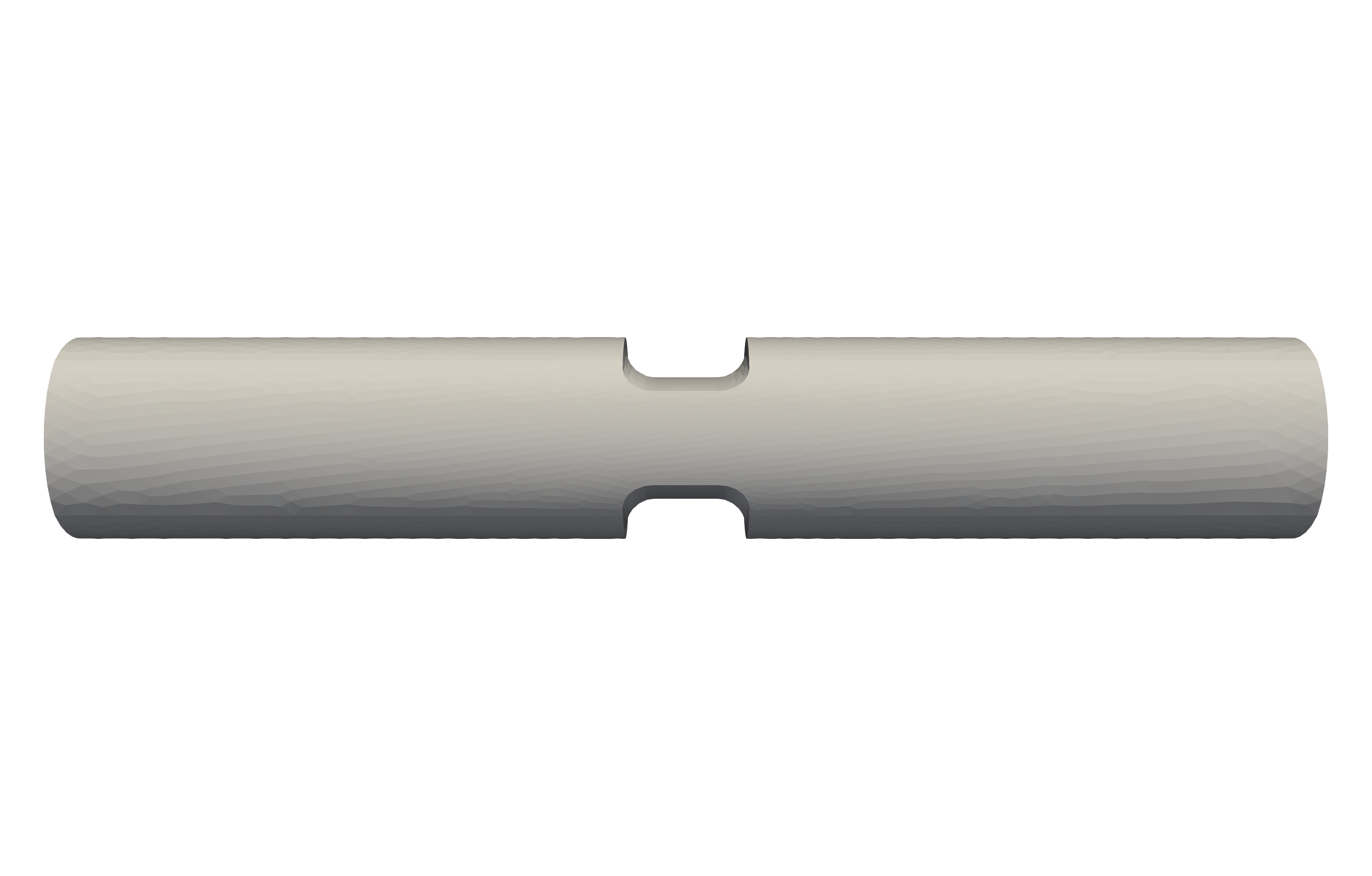}};
       \draw[axis] (0.1*\myW,0.05*\myW)  -- (0.2*\myW,0.05*\myW) node(xline)[right]{$x$};
       \draw[axis] (0.1*\myW,0.05*\myW)  -- (0.1*\myW,0.15*\myW) node(xline)[right]{$y$};
\end{tikzpicture}
}\subcaptionbox{Top view $t=0.0\,$s.\label{fig:arteryGeoInitialTop}}{
\begin{tikzpicture}[
    axis/.style={thick, ->},
    every node/.style={color=black}
    ]
       \node at (0,0) [inner sep=0pt, anchor=south west] {\includegraphics[width=\myW,trim={0cm 20cm 0cm 0cm},clip]{
       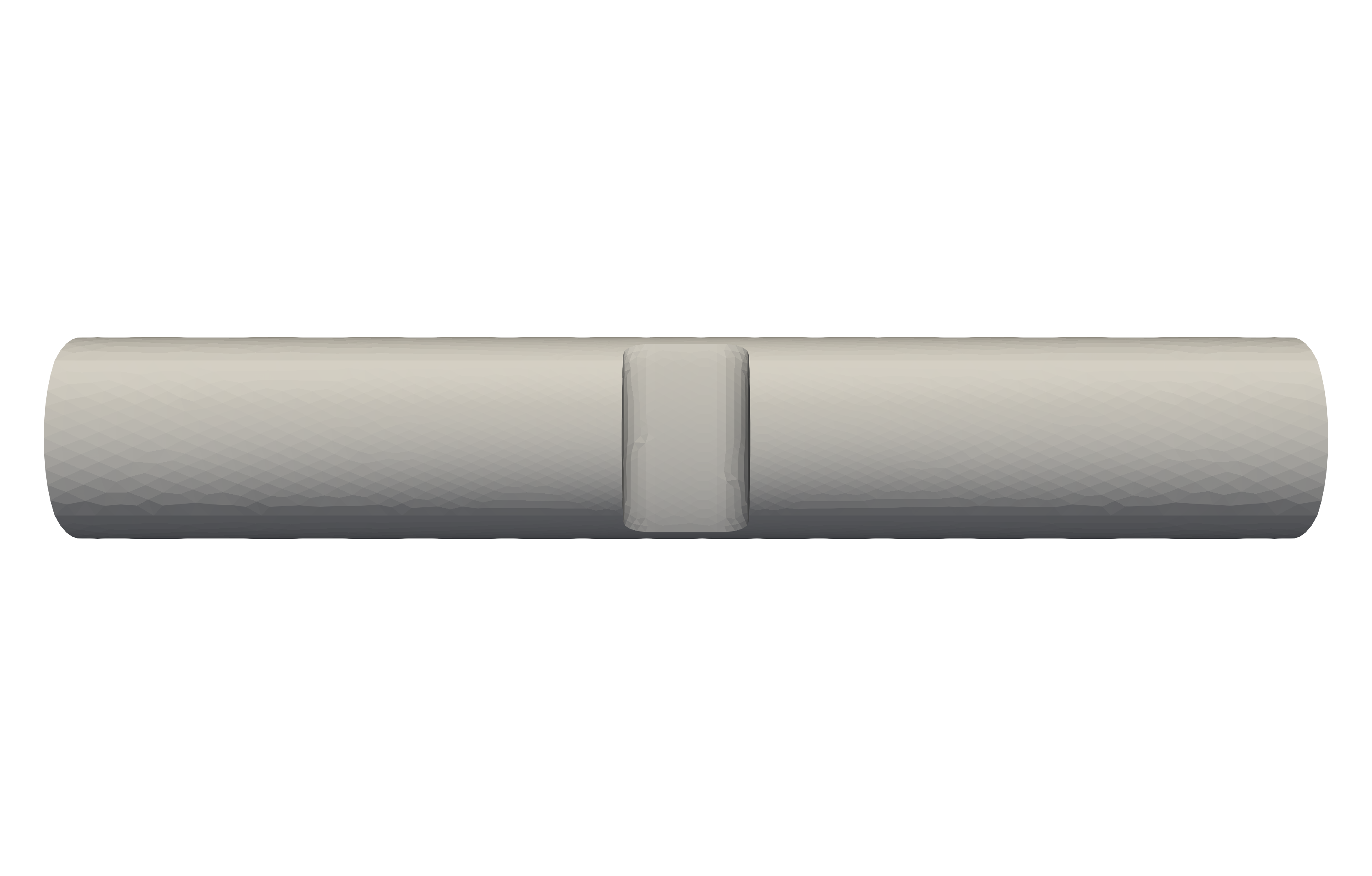}};
       \draw[axis] (0.1*\myW,0.15*\myW)  -- (0.2*\myW,0.15*\myW) node(xline)[right]{$x$};
       \draw[axis] (0.1*\myW,0.15*\myW)  -- (0.1*\myW,0.05*\myW) node(xline)[right]{$z$};
\end{tikzpicture}
}\\
\subcaptionbox{Side view at $t=\SI{1}{\second}$.\label{fig:arteryGeoFinalSide}}{
\begin{tikzpicture}[
    axis/.style={thick, ->},
    every node/.style={color=black}
    ]
       \node at (0,0) [inner sep=0pt, anchor=south west] {\includegraphics[width=\myW,trim={0cm 20cm 0cm 0cm},clip]{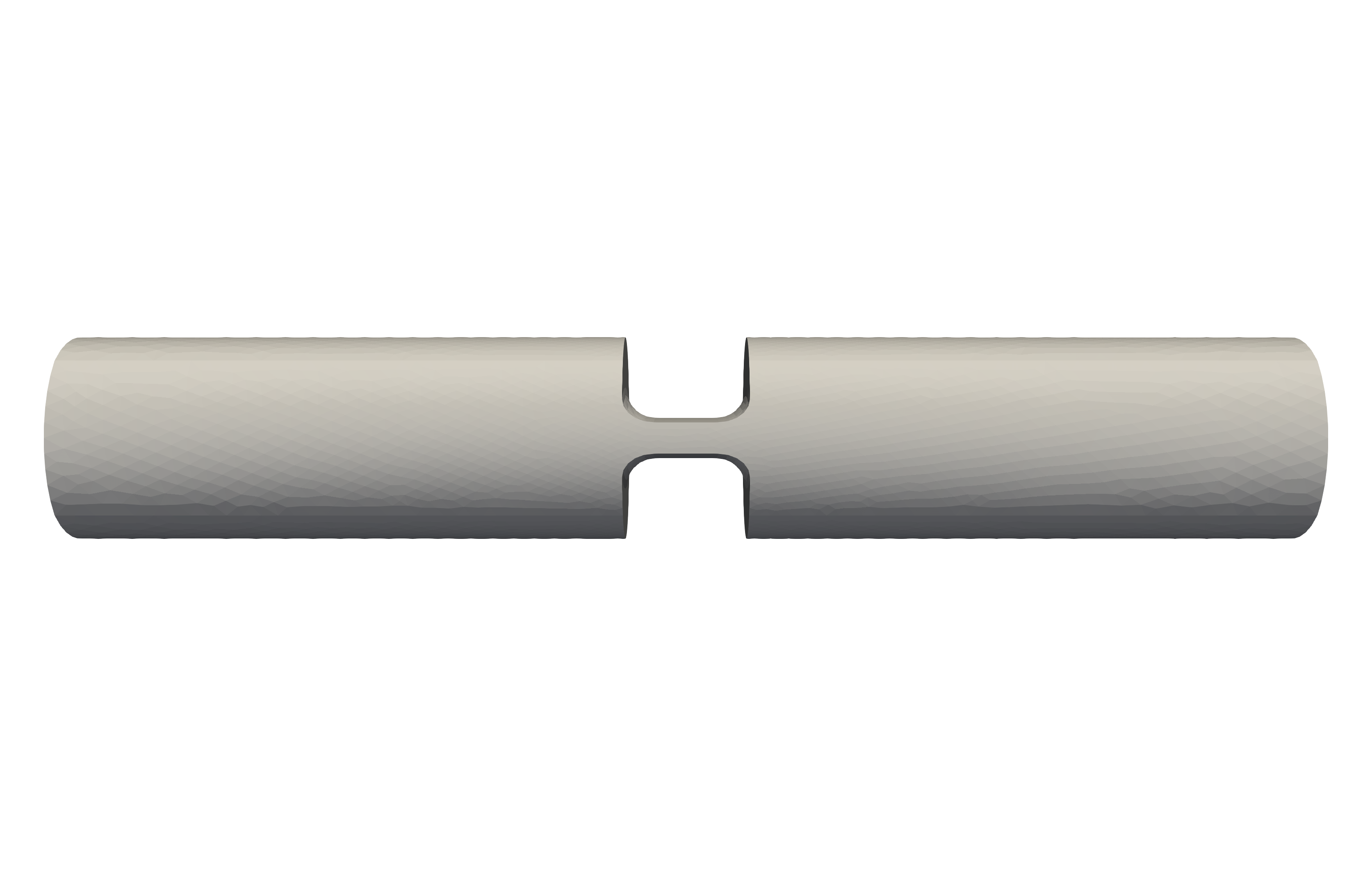}};
       \draw[axis] (0.1*\myW,0.05*\myW)  -- (0.2*\myW,0.05*\myW) node(xline)[right]{$x$};
       \draw[axis] (0.1*\myW,0.05*\myW)  -- (0.1*\myW,0.15*\myW) node(xline)[right]{$y$};
\end{tikzpicture}
}\subcaptionbox{Top view at $t=\SI{1}{\second}$.\label{fig:arteryGeoFinalTop}}{
 \begin{tikzpicture}[
    axis/.style={thick, ->},
    every node/.style={color=black}
    ]
       \node at (0,0) [inner sep=0pt, anchor=south west] {\includegraphics[width=\myW,trim={0cm 20cm 0cm 0cm},clip]{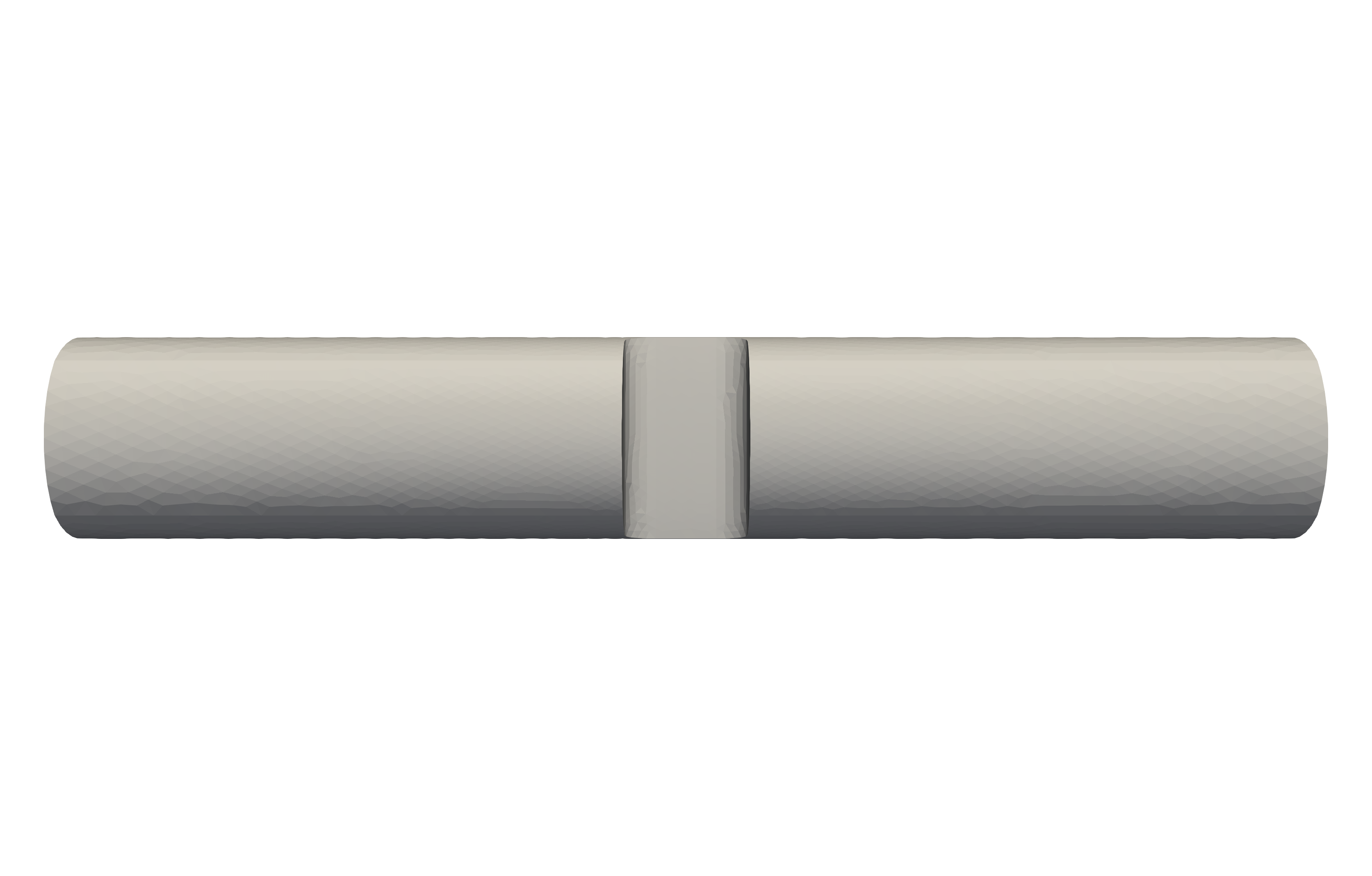}};
       \draw[axis] (0.1*\myW,0.15*\myW)  -- (0.2*\myW,0.15*\myW) node(xline)[right]{$x$};
       \draw[axis] (0.1*\myW,0.15*\myW)  -- (0.1*\myW,0.05*\myW) node(xline)[right]{$z$};
\end{tikzpicture}
}
\caption{Artery-like test case: geometry.}
\label{fig:arteryGeo}
\end{figure}

\renewcommand{\myW}{0.44\textwidth} 
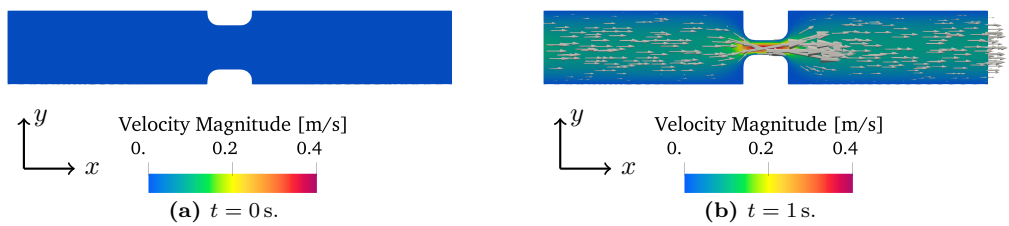
\begin{figure}
\captionsetup[sub]{position=bottom}
\centering

\subcaptionbox{$t=\SI{0}{\second}$.\label{fig:artery-likeVelocityInitial}}{
\begin{tikzpicture}[
    axis/.style={thick, ->},
    every node/.style={color=black}
    ]
       \node at (0,0) [inner sep=0pt, anchor=south west] {\includegraphics[width=\myW,trim={0cm 8cm 0cm 0cm},clip]{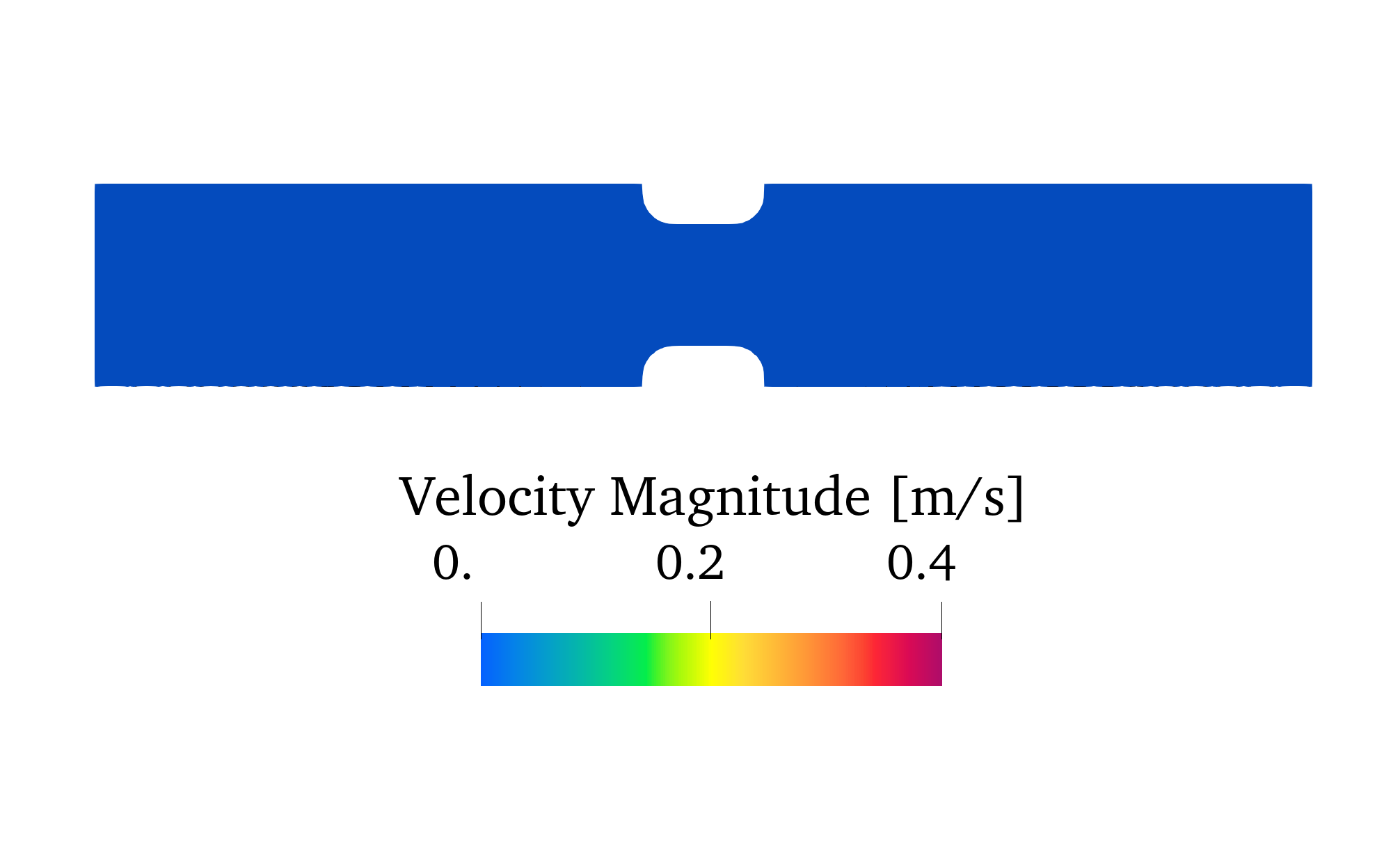}};
       \draw[axis] (0.1*\myW,0.05*\myW)  -- (0.2*\myW,0.05*\myW) node(xline)[right]{$x$};
       \draw[axis] (0.1*\myW,0.05*\myW)  -- (0.1*\myW,0.15*\myW) node(xline)[right]{$y$};
\end{tikzpicture}
}
\subcaptionbox{$t=\SI{1}{\second}$.\label{fig:artery-likeVelocityFinal}}{
\begin{tikzpicture}[
    axis/.style={thick, ->},
    every node/.style={color=black}
    ]
       \node at (0,0) [inner sep=0pt, anchor=south west] {\includegraphics[width=\myW,trim={0cm 8cm 0cm 0cm},clip]{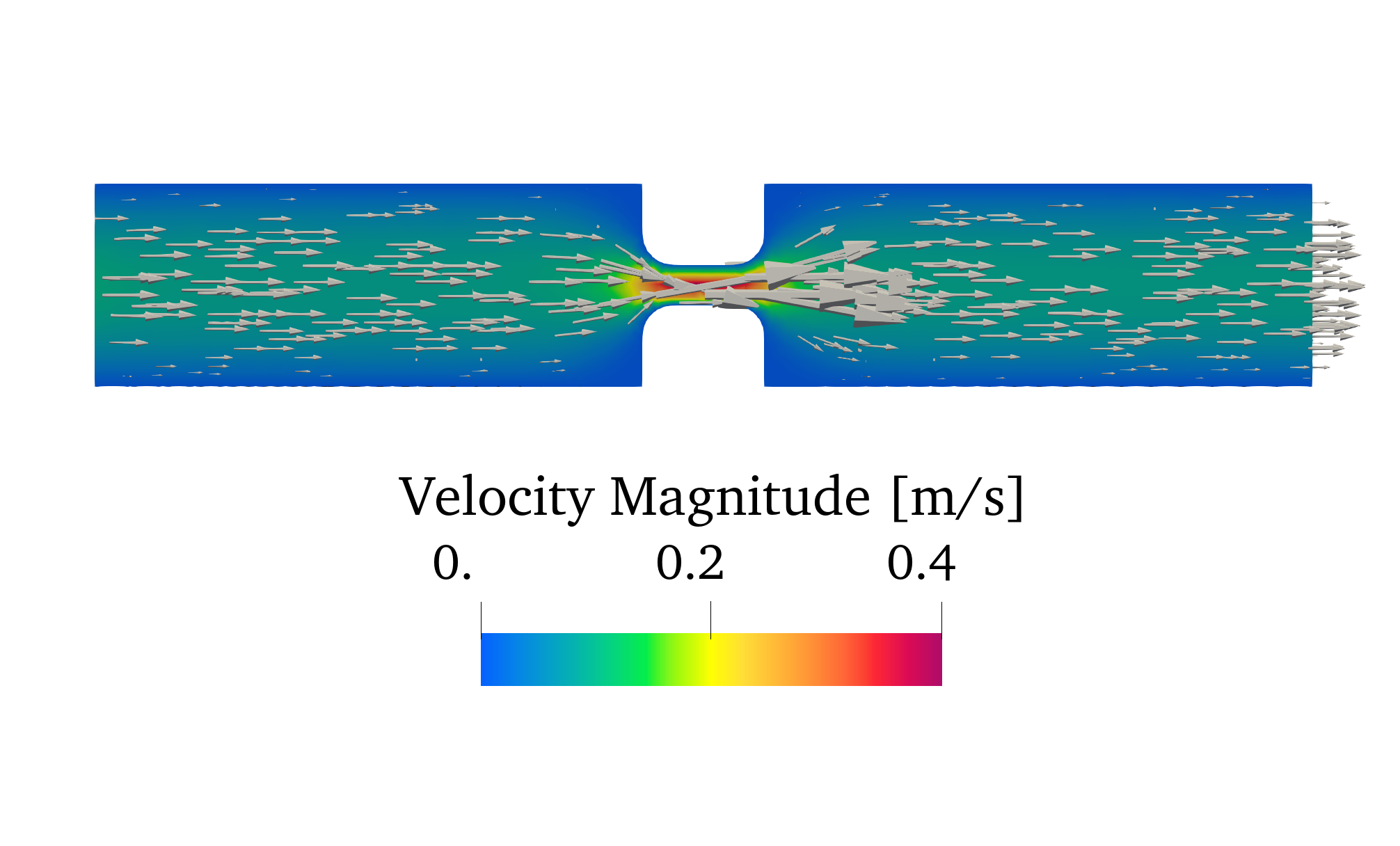}};
       \draw[axis] (0.1*\myW,0.05*\myW)  -- (0.2*\myW,0.05*\myW) node(xline)[right]{$x$};
       \draw[axis] (0.1*\myW,0.05*\myW)  -- (0.1*\myW,0.15*\myW) node(xline)[right]{$y$};
\end{tikzpicture}
}
\caption{Artery-like test case: velocity field for the initial and final state.}
\label{fig:artery-likeVelocity}
\end{figure}

\subsubsection{\gls{rom} for the Flow of Blood in an Artery-Like Geometry}
In the following, a \gls{rom} is constructed for a variation of the prescribed inflow velocity, i.e., $\paramVec = [u^0_{\text{in}}]$ with $\paramVec \in [0.95u^0_{\text{in}}, 1.05u^0_{\text{in}}]$. In this case, we have to distinguish between the inlet boundary portion with a parameter-dependent Dirichlet boundary condition and the moving artery walls with prescribed values that are non-zero but parameter-independent. Thus, we make use of two lifting functions here.
\par
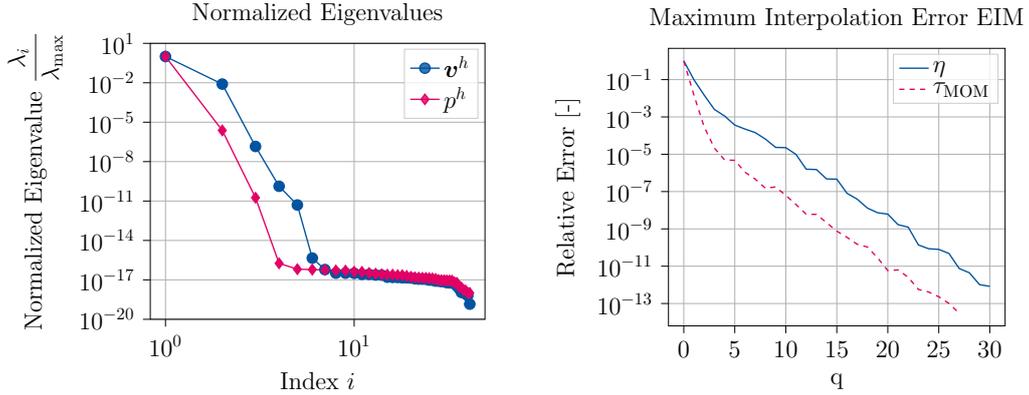
\begin{figure}
    \captionsetup[sub]{position=bottom}
    \centering
    \Large
    \subcaptionbox{Distribution of the eigenvalues from the \gls{pod}.\label{fig:artery-likeEigenvalues}}{
        \resizebox{0.44\textwidth}{!}{
\begin{tikzpicture}

\definecolor{color0}{rgb}{0,0.329411764705882,0.623529411764706}
\definecolor{color1}{rgb}{0.890196078431372,0,0.4}

\begin{axis}[
legend cell align={left},
legend style={fill opacity=0.8, draw opacity=1, text opacity=1, draw=white!80!black},
log basis x={10},
log basis y={10},
tick align=outside,
tick pos=left,
title={Normalized Eigenvalues},
x grid style={white!69.0196078431373!black},
xlabel={Index \(\displaystyle i\)},
xmajorgrids,
xmin=0.830540485032707, xmax=49.3654442364545,
xmode=log,
xtick style={color=black},
xtick={0.01,0.1,1,10,100,1000},
xticklabels={
  \(\displaystyle 10^{-2}\),
  \(\displaystyle 10^{-1}\),
  \(\displaystyle 10^{0}\),
  \(\displaystyle 10^{1}\),
  \(\displaystyle 10^{2}\),
  \(\displaystyle 10^{3}\)
},
y grid style={white!69.0196078431373!black},
ylabel={Normalized Eigenvalue \(\displaystyle \frac{\lambda_i}{\lambda_{\text{max}}}\)},
ymajorgrids,
ymin=1e-20, ymax=10,
ymode=log,
ytick style={color=black},
ytick={1e-23,1e-20,1e-17,1e-14,1e-11,1e-08,1e-05,0.01,10,10000},
yticklabels={
  \(\displaystyle 10^{-23}\),
  \(\displaystyle 10^{-20}\),
  \(\displaystyle 10^{-17}\),
  \(\displaystyle 10^{-14}\),
  \(\displaystyle 10^{-11}\),
  \(\displaystyle 10^{-8}\),
  \(\displaystyle 10^{-5}\),
  \(\displaystyle 10^{-2}\),
  \(\displaystyle 10^{1}\),
  \(\displaystyle 10^{4}\)
}
]
\addplot [thick, color0, mark=*, mark size=3, mark options={solid}]
table {%
1 1
2 0.00794975858415845
3 1.42582257382325e-07
4 1.32791974403868e-10
5 5.00366283890635e-12
6 4.54114155240754e-16
7 6.04573529497411e-17
8 3.26201133022688e-17
9 3.22574638852698e-17
10 3.15358446366046e-17
11 2.62454621731955e-17
12 2.60840954654917e-17
13 2.33242093571454e-17
14 2.28813585566599e-17
15 1.60218882337199e-17
16 1.54306862706516e-17
17 1.52807821353466e-17
18 1.44864745846082e-17
19 1.43851705211444e-17
20 1.38677194525787e-17
21 1.21734236816605e-17
22 1.14897279038485e-17
23 1.13956935268721e-17
24 1.12564920102658e-17
25 9.63808041350083e-18
26 9.58979066065725e-18
27 8.13046505587439e-18
28 8.00740202533402e-18
29 7.46188080562777e-18
30 7.26718042554605e-18
31 6.11314351785897e-18
32 5.93273799548898e-18
33 5.84970067761768e-18
34 5.29939613242502e-18
35 3.47549460184894e-18
36 2.1905261247255e-18
37 1.12228744856004e-18
38 1.01406198292573e-18
39 9.59373162820783e-19
40 6.26226288149959e-19
41 1.42523827201686e-19
};
\addlegendentry{$\trialVelocityHomDiscrete$}
\addplot [thick, color1, mark=diamond*, mark size=3, mark options={solid}]
table {%
1 1
2 2.37328846603967e-06
3 1.78817929292869e-11
4 1.81154657334976e-16
5 6.51678466039161e-17
6 6.02235881107004e-17
7 5.69247083103718e-17
8 5.20065737906182e-17
9 5.01348143603465e-17
10 4.20958043774671e-17
11 4.07622974109276e-17
12 3.32697934377648e-17
13 2.98511733192585e-17
14 2.64421964531683e-17
15 2.52583071694472e-17
16 2.48386648445979e-17
17 2.27569426656641e-17
18 2.27247846357554e-17
19 1.99821382141988e-17
20 1.70491132142761e-17
21 1.54624364513944e-17
22 1.49112670601225e-17
23 1.45406520529215e-17
24 1.40935012359226e-17
25 1.35426769923973e-17
26 1.2976308799412e-17
27 1.18179500955933e-17
28 9.80114361457536e-18
29 9.51423374282225e-18
30 9.36239599182889e-18
31 9.08584610553105e-18
32 8.96571848518543e-18
33 8.78366822577326e-18
34 6.40001828288954e-18
35 5.93546015704658e-18
36 2.97856912325833e-18
37 2.64216217074413e-18
38 1.67829393129063e-18
39 1.53977850304594e-18
40 1.02797999983454e-18
41 1.00473640799576e-18
};
\addlegendentry{$p^h$}
\end{axis}

\end{tikzpicture}
        }
    }
    \subcaptionbox{Maximum interpolation error during the \gls{eim}.\label{fig:artery-likeEIMErrors}}{
        \resizebox{0.44\textwidth}{!}{
\begin{tikzpicture}

\definecolor{color0}{rgb}{0,0.329411764705882,0.623529411764706}
\definecolor{color1}{rgb}{0.890196078431372,0,0.4}

\begin{axis}[
legend cell align={left},
legend style={fill opacity=0.8, draw opacity=1, text opacity=1, draw=white!80!black},
log basis y={10},
tick align=outside,
tick pos=left,
title={Maximum Interpolation Error EIM},
x grid style={white!69.0196078431373!black},
xlabel={q},
xmajorgrids,
xmin=-1.5, xmax=31.5,
xtick style={color=black},
xtick={-5,0,5,10,15,20,25,30,35},
xticklabels={
  \(\displaystyle -5\),
  \(\displaystyle 0\),
  \(\displaystyle 5\),
  \(\displaystyle 10\),
  \(\displaystyle 15\),
  \(\displaystyle 20\),
  \(\displaystyle 25\),
  \(\displaystyle 30\),
  \(\displaystyle 35\)
},
y grid style={white!69.0196078431373!black},
ylabel={Relative Error [-]},
ymajorgrids,
ymin=6.32307985340846e-15, ymax=4.74401718941179,
ymode=log,
ytick style={color=black},
ytick={1e-17,1e-15,1e-13,1e-11,1e-09,1e-07,1e-05,0.001,0.1,10,1000},
yticklabels={
  \(\displaystyle 10^{-17}\),
  \(\displaystyle 10^{-15}\),
  \(\displaystyle 10^{-13}\),
  \(\displaystyle 10^{-11}\),
  \(\displaystyle 10^{-9}\),
  \(\displaystyle 10^{-7}\),
  \(\displaystyle 10^{-5}\),
  \(\displaystyle 10^{-3}\),
  \(\displaystyle 10^{-1}\),
  \(\displaystyle 10^{1}\),
  \(\displaystyle 10^{3}\)
}
]
\addplot [thick, color0]
table {%
0 1
1 0.0973104687076621
2 0.0147633179066779
3 0.00252896805406514
4 0.00111793965375263
5 0.000366212967264887
6 0.000223970653068922
7 0.000145396874338459
8 6.52046804767942e-05
9 2.33321752580709e-05
10 2.26000046016262e-05
11 9.83915149206451e-06
12 1.58107750750885e-06
13 1.51597207604619e-06
14 4.75253511891816e-07
15 4.62299885501037e-07
16 8.06828170742323e-08
17 3.91006462079267e-08
18 1.30274119864675e-08
19 7.34041907118899e-09
20 6.15222072195395e-09
21 1.68400090499216e-09
22 1.23647681875085e-09
23 1.37973957655849e-10
24 8.64119044176524e-11
25 7.94881273897501e-11
26 4.80229695295645e-11
27 7.53378014734756e-12
28 4.43333366049816e-12
29 1.02026522151379e-12
30 8.36384533060246e-13
};
\addlegendentry{$\eta$}
\addplot [thick, color1, dashed]
table {%
0 1
1 0.0142591882104696
2 0.000290151916957558
3 2.16083662357261e-05
4 4.98689894693042e-06
5 4.70515968704204e-06
6 1.0922818762808e-06
7 4.67179954054551e-07
8 1.51464068402887e-07
9 1.77303957247854e-07
10 6.20554933265058e-08
11 2.02888056647836e-08
12 6.08014447268978e-09
13 5.99435297181841e-09
14 2.15022814619033e-09
15 7.55822336917672e-10
16 3.5378735937984e-10
17 1.49573256052504e-10
18 1.06992880868169e-10
19 2.82133502676071e-11
20 5.78225645869476e-12
21 6.30593569710007e-12
22 2.5585361244864e-12
23 5.685980701701e-13
24 4.16003286360201e-13
25 2.3040791442305e-13
26 1.00534473647095e-13
27 2.99967995145931e-14
};
\addlegendentry{$\tMomPlain$}
\end{axis}

\end{tikzpicture}
        }
    }
    \caption{Artery-like test case: results from the construction of the \gls{rom}.}
    \label{fig:artery-likeOffline}
\end{figure}

First, we compute snapshots with the \gls{fom} for $\nTrain = 41$ training samples that are equidistantly distributed. Here, the \gls{fom} involves $\nBasisFOM = 2,194,390$ \glspl{dof} in total. As a result from the \gls{pod}, the distribution of eigenvalues is depicted in \Cref{fig:artery-likeEigenvalues} and we choose $\nBasisVelocityROM= 1,\dots,7$, including the two velocity lifting functions, and $\nBasisPressureROM = 1,\dots,3$. For the \gls{eim}, we set a tolerance of $\SI{1e-12}{}$ and $\SI{1e-13}{}$ and obtain $Q_{\visc} = 30$ and $Q_{\tau} = 27$, respectively. \Cref{fig:artery-likeEIMErrors} shows the maximum interpolation error for the viscosity $\visc$ and the stabilization parameter $\tMomPlain$ during the \gls{eim}. 
\par
\begin{figure}
    \captionsetup[sub]{position=bottom}
    \centering
    \Large
    \subcaptionbox{Velocity.}{
        \resizebox{0.44\textwidth}{!}{
\begin{tikzpicture}

\definecolor{color0}{rgb}{1,0.87843137254902,0.12156862745098}
\definecolor{color1}{rgb}{1,0.627450980392157,0.372549019607843}
\definecolor{color2}{rgb}{1,0.376470588235294,0.623529411764706}
\definecolor{color3}{rgb}{1,0.125490196078431,0.874509803921569}

\begin{axis}[
tick align=outside,
tick pos=left,
title={$\relErrorVelocity$},
x grid style={white!69.0196078431373!black},
xlabel={\(\displaystyle N_p\)},
xmin=1, xmax=3,
xtick style={color=black},
xtick={1,2,3},
xticklabels={\(\displaystyle 1\),\(\displaystyle 2\),\(\displaystyle 3\)},
y grid style={white!69.0196078431373!black},
ylabel={\(\displaystyle N_{\mathbf{u}}\)},
ymin=3, ymax=7,
ytick style={color=black},
ytick={3,4,5,6,7},
yticklabels={
  \(\displaystyle 3\),
  \(\displaystyle 4\),
  \(\displaystyle 5\),
  \(\displaystyle 6\),
  \(\displaystyle 7\)
}
]
\addplot [draw=none, fill=color0]
table{%
x  y
2 4.37650891394412
3 4.35330770412385
3 5
3 6
3 7
2 7
1.99040991259345 7
1.9903127173998 6
1.98760845108052 5
2 4.37650891394412
};
\addplot [draw=none, fill=color1]
table{%
x  y
2 3.92579367910137
3 3.92540471448223
3 4
3 4.35330770412385
2 4.37650891394412
1.98760845108052 5
1.9903127173998 6
1.99040991259345 7
1.86869012295147 7
1.86619036626405 6
1.85013870710813 5
1.80309130655496 4
2 3.92579367910137
};
\addplot [draw=none, fill=color2]
table{%
x  y
2 3.14101234840689
3 3.14081265825156
3 3.92540471448223
2 3.92579367910137
1.80309130655496 4
1.85013870710813 5
1.86619036626405 6
1.86869012295147 7
1 7
1 6
1 5
1 4
1 3.22895886838152
2 3.14101234840689
};
\addplot [draw=none, fill=color3]
table{%
x  y
2 3
3 3
3 3.14081265825156
2 3.14101234840689
1 3.22895886838152
1 3
2 3
};
\path [draw=black, semithick]
(axis cs:3,4.35330770412385)
--(axis cs:2,4.37650891394412)
--(axis cs:1.99313704618645,4.72182414084339);

\path [draw=black, semithick]
(axis cs:1.98836176234379,5.27856400751013)
--(axis cs:1.9903127173998,6)
--(axis cs:1.99040991259345,7);

\path [draw=black, semithick]
(axis cs:3,3.92540471448223)
--(axis cs:2,3.92579367910137)
--(axis cs:1.80309130655496,4)
--(axis cs:1.83713590577805,4.72362338456115);

\path [draw=black, semithick]
(axis cs:1.85460609537104,5.27831317744237)
--(axis cs:1.86619036626405,6)
--(axis cs:1.86869012295147,7);

\path [draw=black, semithick]
(axis cs:3,3.14081265825156)
--(axis cs:2.10379029115697,3.14099162250753);

\path [draw=black, semithick]
(axis cs:1.89626538289669,3.15013544698203)
--(axis cs:1,3.22895886838152);

\draw (axis cs:1.98760845108052,5) node[
  scale=0.8,
  text=black,
  rotate=271.3
]{$10^{-5}$};
\draw (axis cs:1.85013870710813,5) node[
  scale=0.8,
  text=black,
  rotate=85.2
]{$10^{-4}$};
\draw (axis cs:2,3.14101234840689) node[
  scale=0.8,
  text=black,
  rotate=359.1
]{$10^{-3}$};
\end{axis}

\end{tikzpicture}}
    }
    \subcaptionbox{Pressure.}{
        \resizebox{0.44\textwidth}{!}{
\begin{tikzpicture}

\definecolor{color0}{rgb}{1,0.835294117647059,0.164705882352941}
\definecolor{color1}{rgb}{1,0.501960784313725,0.498039215686275}
\definecolor{color2}{rgb}{1,0.164705882352941,0.835294117647059}

\begin{axis}[
tick align=outside,
tick pos=left,
title={$\relErrorPressure$},
x grid style={white!69.0196078431373!black},
xlabel={\(\displaystyle N_p\)},
xmin=1, xmax=3,
xtick style={color=black},
xtick={1,2,3},
xticklabels={\(\displaystyle 1\),\(\displaystyle 2\),\(\displaystyle 3\)},
y grid style={white!69.0196078431373!black},
ylabel={\(\displaystyle N_{\mathbf{u}}\)},
ymin=3, ymax=7,
ytick style={color=black},
ytick={3,4,5,6,7},
yticklabels={
  \(\displaystyle 3\),
  \(\displaystyle 4\),
  \(\displaystyle 5\),
  \(\displaystyle 6\),
  \(\displaystyle 7\)
}
]
\addplot [draw=none, fill=color0]
table{%
x  y
3 3.99841272096808
3 4
3 5
3 6
3 7
2 7
1.99231388686608 7
1.9923086500537 6
1.99186016592534 5
2 4.02906689026844
2.0929081054864 4
3 3.99841272096808
};
\addplot [draw=none, fill=color1]
table{%
x  y
2 3.85666988114318
3 3.84783564712328
3 3.99841272096808
2.0929081054864 4
2 4.02906689026844
1.99186016592534 5
1.9923086500537 6
1.99231388686608 7
1.78528526361933 7
1.78346767600251 6
1.76616198020346 5
1.71107969290604 4
2 3.85666988114318
};
\addplot [draw=none, fill=color2]
table{%
x  y
2 3
3 3
3 3.84783564712328
2 3.85666988114318
1.71107969290604 4
1.76616198020346 5
1.78346767600251 6
1.78528526361933 7
1 7
1 6
1 5
1 4
1 3
2 3
};
\path [draw=black, semithick]
(axis cs:3,3.99841272096808)
--(axis cs:2.0929081054864,4)
--(axis cs:2,4.02906689026844)
--(axis cs:1.99419498249751,4.72149914921187);

\path [draw=black, semithick]
(axis cs:1.99198510066168,5.278571143184)
--(axis cs:1.9923086500537,6)
--(axis cs:1.99231388686608,7);

\path [draw=black, semithick]
(axis cs:3,3.84783564712328)
--(axis cs:2,3.85666988114318)
--(axis cs:1.71107969290604,4)
--(axis cs:1.75098262156661,4.72442395946829);

\path [draw=black, semithick]
(axis cs:1.77097765920204,5.27827133069329)
--(axis cs:1.78346767600251,6)
--(axis cs:1.78528526361933,7);

\draw (axis cs:1.99186016592534,5) node[
  scale=0.8,
  text=black,
  rotate=270.6
]{$10^{-4}$};
\draw (axis cs:1.76616198020346,5) node[
  scale=0.8,
  text=black,
  rotate=84.5
]{$10^{-3}$};
\end{axis}

\end{tikzpicture}}
    }
    \caption{Artery-like test case: maximum relative error of the \gls{rom} over all testing samples.}
    \label{fig:artery-likeRelativeErrorMax}
\end{figure}
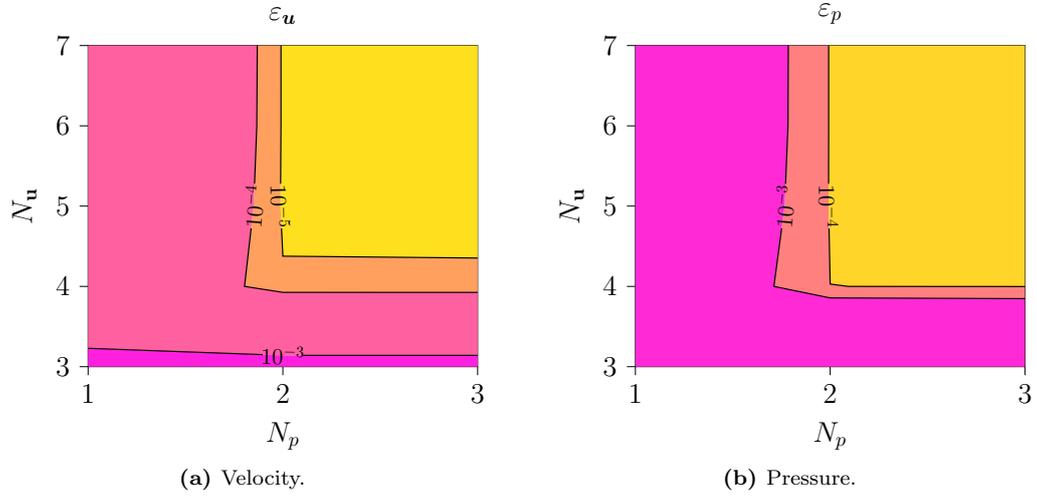

Subsequently, we use $\nTest= 20$ uniformly distributed random samples to carry out the error and performance analysis. The results for the maximal relative errors $\relErrorVelocity$ and $\relErrorPressure$ are presented in \Cref{fig:artery-likeRelativeErrorMax}. The qualitative behavior of the errors is very similar as for the valve-like test case (see \Cref{subsubsec:valve-likeROM}). However, the magnitude of the error is not as small as before. For the velocity, it is smaller than \SI{1e-2}{} and decreases to below \SI{1e-5}{}. In a like manner, the error in the pressure ranges from values smaller than \SI{1e-2}{} to values less than \SI{1e-4}{}.   
\par
\par
\begin{figure}
    \captionsetup[sub]{position=bottom}
    \centering
    \Large
    \subcaptionbox{Average speed up over all testing samples.}{
        \resizebox{0.44\textwidth}{!}{
\begin{tikzpicture}

\definecolor{color0}{rgb}{1,0.952941176470588,0.0470588235294118}
\definecolor{color1}{rgb}{1,0.850980392156863,0.149019607843137}
\definecolor{color2}{rgb}{1,0.749019607843137,0.250980392156863}
\definecolor{color3}{rgb}{1,0.650980392156863,0.349019607843137}
\definecolor{color4}{rgb}{1,0.549019607843137,0.450980392156863}
\definecolor{color5}{rgb}{1,0.450980392156863,0.549019607843137}
\definecolor{color6}{rgb}{1,0.349019607843137,0.650980392156863}
\definecolor{color7}{rgb}{1,0.247058823529412,0.752941176470588}
\definecolor{color8}{rgb}{1,0.149019607843137,0.850980392156863}
\definecolor{color9}{rgb}{1,0.0470588235294118,0.952941176470588}

\begin{axis}[
tick align=outside,
tick pos=left,
title={Speedup - Average},
x grid style={white!69.0196078431373!black},
xlabel={\(\displaystyle N_p\)},
xmin=1, xmax=3,
xtick style={color=black},
xtick={1,2,3},
xticklabels={\(\displaystyle 1\),\(\displaystyle 2\),\(\displaystyle 3\)},
y grid style={white!69.0196078431373!black},
ylabel={\(\displaystyle N_{\mathbf{u}}\)},
ymin=3, ymax=7,
ytick style={color=black},
ytick={3,4,5,6,7},
yticklabels={
  \(\displaystyle 3\),
  \(\displaystyle 4\),
  \(\displaystyle 5\),
  \(\displaystyle 6\),
  \(\displaystyle 7\)
}
]
\addplot [draw=none, fill=color0]
table{%
x  y
2 3.99776029877739
3 3.97186061500088
3 4
3 4.05129247948725
2 4.00401665310977
1.64783776122399 4
2 3.99776029877739
};
\addplot [draw=none, fill=color1]
table{%
x  y
2 3.86632751738316
3 3.82489079580478
3 3.97186061500088
2 3.99776029877739
1.64783776122399 4
2 4.00401665310977
3 4.05129247948725
3 4.31918911036938
2 4.23972661512697
1 4.36447852644752
1 4
1 3.88899847255399
2 3.86632751738316
};
\addplot [draw=none, fill=color1]
table{%
x  y
2 5.95377752187074
3 5.82647002301902
3 6
3 7
2 7
1.13955470147603 7
1.68538959041588 6
2 5.95377752187074
};
\addplot [draw=none, fill=color2]
table{%
x  y
2 3.73489473598894
3 3.67792097660869
3 3.82489079580478
2 3.86632751738316
1 3.88899847255399
1 3.77440465802077
2 3.73489473598894
};
\addplot [draw=none, fill=color2]
table{%
x  y
2 4.23972661512697
3 4.31918911036938
3 4.58708574125151
2 4.47543657714418
1 4.74075248971692
1 4.36447852644752
2 4.23972661512697
};
\addplot [draw=none, fill=color2]
table{%
x  y
2 5.65807486026756
3 5.50125704822069
3 5.82647002301902
2 5.95377752187074
1.68538959041588 6
1.13955470147603 7
1 7
1 6
1 5.51094647432634
2 5.65807486026756
};
\addplot [draw=none, fill=color3]
table{%
x  y
2 3.60346195459471
3 3.5309511574126
3 3.67792097660869
2 3.73489473598894
1 3.77440465802077
1 3.65981084348756
2 3.60346195459471
};
\addplot [draw=none, fill=color3]
table{%
x  y
2 4.47543657714418
3 4.58708574125151
3 4.85498237213365
2 4.71114653916139
1.20242039825908 5
2 5.36237219866437
3 5.17604407342235
3 5.50125704822069
2 5.65807486026756
1 5.51094647432634
1 5
1 4.74075248971692
2 4.47543657714418
};
\addplot [draw=none, fill=color4]
table{%
x  y
2 3.47202917320049
3 3.38398133821651
3 3.5309511574126
2 3.60346195459471
1 3.65981084348756
1 3.54521702895434
2 3.47202917320049
};
\addplot [draw=none, fill=color4]
table{%
x  y
2 4.71114653916139
3 4.85498237213365
3 5
3 5.17604407342235
2 5.36237219866437
1.20242039825908 5
2 4.71114653916139
1.85326064468108 5
2 5.06666953706118
2.3295534889383 5
2 4.94685650117859
1.85326064468108 5
};
\addplot [draw=none, fill=color5]
table{%
x  y
2 3.34059639180626
3 3.23701151902041
3 3.38398133821651
2 3.47202917320049
1 3.54521702895434
1 3.43062321442112
2 3.34059639180626
};
\addplot [draw=none, fill=color5]
table{%
x  y
2 4.94685650117859
2.3295534889383 5
2 5.06666953706118
1.85326064468108 5
2 4.94685650117859
};
\addplot [draw=none, fill=color6]
table{%
x  y
2 3.20916361041203
3 3.09004169982432
3 3.23701151902041
2 3.34059639180626
1 3.43062321442112
1 3.3160293998879
2 3.20916361041203
};
\addplot [draw=none, fill=color7]
table{%
x  y
2 3.07773082901781
2.60424735056556 3
3 3
3 3.09004169982432
2 3.20916361041203
1 3.3160293998879
1 3.20143558535469
2 3.07773082901781
};
\addplot [draw=none, fill=color8]
table{%
x  y
2 3
2.60424735056556 3
2 3.07773082901781
1 3.20143558535469
1 3.08684177082147
1.64970439596316 3
2 3
};
\addplot [draw=none, fill=color9]
table{%
x  y
1 3.08684177082147
1 3
1.64970439596316 3
1 3.08684177082147
};
\path [draw=black, semithick]
(axis cs:3,3.97186061500088)
--(axis cs:2,3.99776029877739)
--(axis cs:1.77726149904818,3.99917688362941);

\path [draw=black, semithick]
(axis cs:1.77726069380053,4.00147615776869)
--(axis cs:2,4.00401665310977)
--(axis cs:3,4.05129247948725);

\path [draw=black, semithick]
(axis cs:3,3.82489079580478)
--(axis cs:2.12940867994302,3.86096524594254);

\path [draw=black, semithick]
(axis cs:1.87058051564284,3.86926158071125)
--(axis cs:1,3.88899847255399);

\path [draw=black, semithick]
(axis cs:1,4.36447852644752)
--(axis cs:1.87071547689406,4.2558551064886);

\path [draw=black, semithick]
(axis cs:2.12936741661236,4.25000647285406)
--(axis cs:3,4.31918911036938);

\path [draw=black, semithick]
(axis cs:3,5.82647002301902)
--(axis cs:2,5.95377752187074)
--(axis cs:1.81462022259509,5.98101346971504);

\path [draw=black, semithick]
(axis cs:1.57849420274952,6.1958383200348)
--(axis cs:1.13955470147603,7);

\path [draw=black, semithick]
(axis cs:3,3.67792097660869)
--(axis cs:2.1293949518265,3.72752261913855);

\path [draw=black, semithick]
(axis cs:1.87058991946085,3.74000771818117)
--(axis cs:1,3.77440465802077);

\path [draw=black, semithick]
(axis cs:1,4.74075248971692)
--(axis cs:1.87120364443596,4.5096082997567);

\path [draw=black, semithick]
(axis cs:2.12931226725914,4.48987418369249)
--(axis cs:3,4.58708574125151);

\path [draw=black, semithick]
(axis cs:3,5.50125704822069)
--(axis cs:2.12920375458042,5.63781341016601);

\path [draw=black, semithick]
(axis cs:1.87076991649158,5.63906144666591)
--(axis cs:1,5.51094647432634);

\path [draw=black, semithick]
(axis cs:3,3.5309511574126)
--(axis cs:2.12937689557672,3.5940807327595);

\path [draw=black, semithick]
(axis cs:1.87060441248816,3.61075325217864)
--(axis cs:1,3.65981084348756);

\path [draw=black, semithick]
(axis cs:3,5.17604407342235)
--(axis cs:2,5.36237219866437)
--(axis cs:1.33002909148532,5.05797771486651);

\path [draw=black, semithick]
(axis cs:1.33068211116719,4.95354841123603)
--(axis cs:2,4.71114653916139)
--(axis cs:3,4.85498237213365);

\path [draw=black, semithick]
(axis cs:3,3.38398133821651)
--(axis cs:2.12935451663035,3.46063978806578);

\path [draw=black, semithick]
(axis cs:1.87062398959549,3.48149792598798)
--(axis cs:1,3.54521702895434);

\path [draw=black, semithick]
(axis cs:1.98152235758919,4.95354841123603)
--(axis cs:2,4.94685650117859)
--(axis cs:2.3295534889383,5)
--(axis cs:2,5.06666953706118)
--(axis cs:1.98086933790731,5.05797771486651);

\path [draw=black, semithick]
(axis cs:3,3.23701151902041)
--(axis cs:2.12932782172,3.32719998584572);

\path [draw=black, semithick]
(axis cs:1.87064864385909,3.35224148340055)
--(axis cs:1,3.43062321442112);

\path [draw=black, semithick]
(axis cs:3,3.09004169982432)
--(axis cs:2.12929681886861,3.19376152631549);

\path [draw=black, semithick]
(axis cs:1.87067836656755,3.2229836688651)
--(axis cs:1,3.3160293998879);

\path [draw=black, semithick]
(axis cs:2.60424735056556,3)
--(axis cs:2.12927570066499,3.06110070664302);

\path [draw=black, semithick]
(axis cs:1.87071314723031,3.09372422763725)
--(axis cs:1,3.20143558535469);

\path [draw=black, semithick]
(axis cs:1.64970439596316,3)
--(axis cs:1.12926390778472,3.06956390303255);

\draw (axis cs:1.64783776122399,4) node[
  scale=0.8,
  text=black,
  rotate=270.1
]{20000};
\draw (axis cs:2,3.86632751738316) node[
  scale=0.8,
  text=black,
  rotate=359.3
]{25000};
\draw (axis cs:2,4.23972661512697) node[
  scale=0.8,
  text=black,
  rotate=359.5
]{25000};
\draw (axis cs:1.68538959041588,6) node[
  scale=0.8,
  text=black,
  rotate=341.3
]{25000};
\draw (axis cs:2,3.73489473598894) node[
  scale=0.8,
  text=black,
  rotate=359.0
]{30000};
\draw (axis cs:2,4.47543657714418) node[
  scale=0.8,
  text=black,
  rotate=358.4
]{30000};
\draw (axis cs:2,5.65807486026756) node[
  scale=0.8,
  text=black,
  rotate=359.9
]{30000};
\draw (axis cs:2,3.60346195459471) node[
  scale=0.8,
  text=black,
  rotate=358.6
]{35000};
\draw (axis cs:1.20242039825908,5) node[
  scale=0.8,
  text=black,
  rotate=271.0
]{35000};
\draw (axis cs:2,3.47202917320049) node[
  scale=0.8,
  text=black,
  rotate=358.3
]{40000};
\draw (axis cs:1.85326064468108,5) node[
  scale=0.8,
  text=black,
  rotate=271.0
]{40000};
\draw (axis cs:2,3.34059639180626) node[
  scale=0.8,
  text=black,
  rotate=357.9
]{45000};
\draw (axis cs:2,3.20916361041203) node[
  scale=0.8,
  text=black,
  rotate=357.6
]{50000};
\draw (axis cs:2,3.07773082901781) node[
  scale=0.8,
  text=black,
  rotate=357.3
]{55000};
\draw (axis cs:1,3.08684177082147) node[
  scale=0.8,
  text=black,
  rotate=357.1
]{60000};
\end{axis}

\end{tikzpicture}}
    }
    \subcaptionbox{Maximum speed up over all testing samples.}{
        \resizebox{0.44\textwidth}{!}{
\begin{tikzpicture}

\definecolor{color0}{rgb}{1,0.917647058823529,0.0823529411764706}
\definecolor{color1}{rgb}{1,0.749019607843137,0.250980392156863}
\definecolor{color2}{rgb}{1,0.584313725490196,0.415686274509804}
\definecolor{color3}{rgb}{1,0.415686274509804,0.584313725490196}
\definecolor{color4}{rgb}{1,0.247058823529412,0.752941176470588}
\definecolor{color5}{rgb}{1,0.0823529411764706,0.917647058823529}

\begin{axis}[
tick align=outside,
tick pos=left,
title={Speedup - Max},
x grid style={white!69.0196078431373!black},
xlabel={\(\displaystyle N_p\)},
xmin=1, xmax=3,
xtick style={color=black},
xtick={1,2,3},
xticklabels={\(\displaystyle 1\),\(\displaystyle 2\),\(\displaystyle 3\)},
y grid style={white!69.0196078431373!black},
ylabel={\(\displaystyle N_{\mathbf{u}}\)},
ymin=3, ymax=7,
ytick style={color=black},
ytick={3,4,5,6,7},
yticklabels={
  \(\displaystyle 3\),
  \(\displaystyle 4\),
  \(\displaystyle 5\),
  \(\displaystyle 6\),
  \(\displaystyle 7\)
}
]
\addplot [draw=none, fill=color0]
table{%
x  y
2 3.89113969366561
3 3.81295761379902
3 4
3 4.33664393194555
2 4.20059943824401
1 4.31046717020061
1 4
1 3.91000519290614
2 3.89113969366561
};
\addplot [draw=none, fill=color0]
table{%
x  y
3 5.8901285826523
3 6
3 7
2 7
1.35534430289395 7
2 6.33923453361608
2.26948248989546 6
3 5.8901285826523
};
\addplot [draw=none, fill=color1]
table{%
x  y
2 3.66515692779255
3 3.56174755756104
3 3.81295761379902
2 3.89113969366561
1 3.91000519290614
1 3.71929144543007
2 3.66515692779255
};
\addplot [draw=none, fill=color1]
table{%
x  y
2 4.20059943824401
3 4.33664393194555
3 4.78877856727545
2 4.61702317811216
1 4.96839799320352
1 4.31046717020061
2 4.20059943824401
};
\addplot [draw=none, fill=color1]
table{%
x  y
2 5.49635319629926
3 5.28342883044444
3 5.8901285826523
2.26948248989546 6
2 6.33923453361608
1.35534430289395 7
1 7
1 6
1 5.12326574274158
2 5.49635319629926
};
\addplot [draw=none, fill=color2]
table{%
x  y
2 3.4391741619195
3 3.31053750132306
3 3.56174755756104
2 3.66515692779255
1 3.71929144543007
1 3.52857769795399
2 3.4391741619195
};
\addplot [draw=none, fill=color2]
table{%
x  y
2 4.61702317811216
3 4.78877856727545
3 5
3 5.28342883044444
2 5.49635319629926
1 5.12326574274158
1 5
1 4.96839799320352
2 4.61702317811216
};
\addplot [draw=none, fill=color3]
table{%
x  y
2 3.21319139604645
3 3.05932744508508
3 3.31053750132306
2 3.4391741619195
1 3.52857769795399
1 3.33786395047791
2 3.21319139604645
};
\addplot [draw=none, fill=color4]
table{%
x  y
2 3
3 3
3 3.05932744508508
2 3.21319139604645
1 3.33786395047791
1 3.14715020300184
1.93165335182826 3
2 3
};
\addplot [draw=none, fill=color5]
table{%
x  y
1 3.14715020300184
1 3
1.93165335182826 3
1 3.14715020300184
};
\path [draw=black, semithick]
(axis cs:3,3.81295761379902)
--(axis cs:2.12936922751048,3.8810253383881);

\path [draw=black, semithick]
(axis cs:1.87057909585802,3.89358128363441)
--(axis cs:1,3.91000519290614);

\path [draw=black, semithick]
(axis cs:1,4.31046717020061)
--(axis cs:1.87068419689434,4.21480707223737);

\path [draw=black, semithick]
(axis cs:2.12925816089073,4.21818429929918)
--(axis cs:3,4.33664393194555);

\path [draw=black, semithick]
(axis cs:3,5.8901285826523)
--(axis cs:2.26948248989546,6)
--(axis cs:2.11717612375158,6.19172889157148);

\path [draw=black, semithick]
(axis cs:1.87909257266006,6.46316340444465)
--(axis cs:1.35534430289395,7);

\path [draw=black, semithick]
(axis cs:3,3.56174755756104)
--(axis cs:2.12932814733133,3.65178318552381);

\path [draw=black, semithick]
(axis cs:1.87060221616005,3.67216181440409)
--(axis cs:1,3.71929144543007);

\path [draw=black, semithick]
(axis cs:1,4.96839799320352)
--(axis cs:1.87167093357662,4.66211478009752);

\path [draw=black, semithick]
(axis cs:2.12915991192253,4.63920708904871)
--(axis cs:3,4.78877856727545);

\path [draw=black, semithick]
(axis cs:3,5.28342883044444)
--(axis cs:2.12901874905936,5.46888196097241);

\path [draw=black, semithick]
(axis cs:1.87180845969371,5.44852654095875)
--(axis cs:1,5.12326574274158);

\path [draw=black, semithick]
(axis cs:3,3.31053750132306)
--(axis cs:2.12927571007098,3.42254456627974);

\path [draw=black, semithick]
(axis cs:1.87064764090286,3.45073872021719)
--(axis cs:1,3.52857769795399);

\path [draw=black, semithick]
(axis cs:3,3.05932744508508)
--(axis cs:2.12921195720541,3.19331033379938);

\path [draw=black, semithick]
(axis cs:1.87071529966112,3.22930964988661)
--(axis cs:1,3.33786395047791);

\path [draw=black, semithick]
(axis cs:1.93165335182826,3)
--(axis cs:1.12920058312662,3.12674358717447);

\draw (axis cs:2,3.89113969366561) node[
  scale=0.8,
  text=black,
  rotate=359.0
]{30000};
\draw (axis cs:2,4.20059943824401) node[
  scale=0.8,
  text=black,
  rotate=0.3
]{30000};
\draw (axis cs:2,6.33923453361608) node[
  scale=0.8,
  text=black,
  rotate=337.0
]{30000};
\draw (axis cs:2,3.66515692779255) node[
  scale=0.8,
  text=black,
  rotate=358.3
]{40000};
\draw (axis cs:2,4.61702317811216) node[
  scale=0.8,
  text=black,
  rotate=358.1
]{40000};
\draw (axis cs:2,5.49635319629926) node[
  scale=0.8,
  text=black,
  rotate=1.7
]{40000};
\draw (axis cs:2,3.4391741619195) node[
  scale=0.8,
  text=black,
  rotate=357.7
]{50000};
\draw (axis cs:2,3.21319139604645) node[
  scale=0.8,
  text=black,
  rotate=357.0
]{60000};
\draw (axis cs:1,3.14715020300184) node[
  scale=0.8,
  text=black,
  rotate=356.6
]{70000};
\end{axis}

\end{tikzpicture}}
    }
    \caption{Artery-like test case: performance results.}
    \label{fig:artery-likePerformance}
\end{figure}
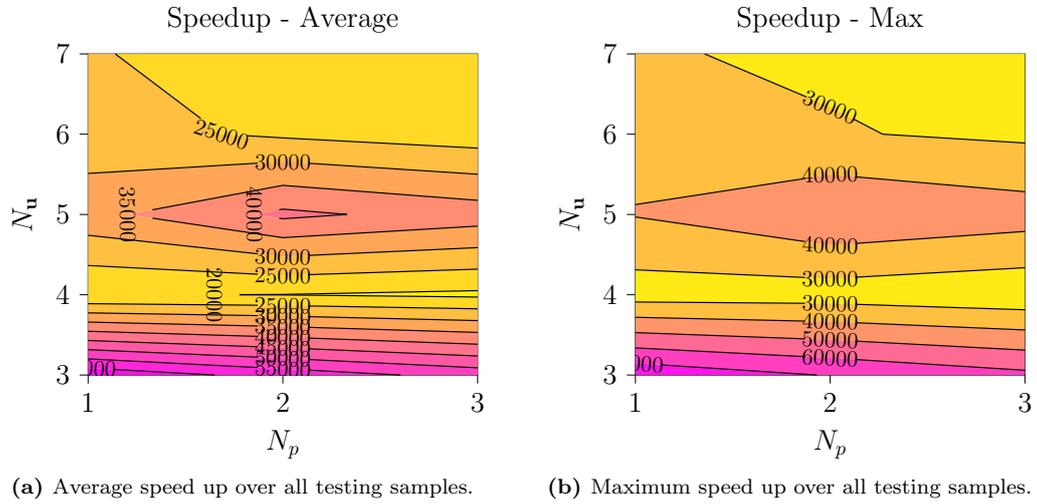

Next, we present the results from the performance analysis. The \gls{rom} has been evaluated on a single core, whereas we have used 240 cores to compute the \gls{fom}. The measured speed up --- the average and the maximum over all the testing samples --- is shown in \Cref{fig:artery-likePerformance}. The results are, to some extent, counter-intuitive. First, one can observe that the speed up is mainly governed by $\nBasisVelocityROM$ and roughly constant for  $\nBasisPressureROM$. Second, the speed up does not change monotonically with increasing $\nBasisVelocityROM$ but has a local maximum at $\nBasisVelocityROM=5$. Analyzing the evaluation of the \gls{rom} in more detail has revealed that this pattern originates from the number of iterations needed to solve the non-linear reduced system. In particular, the value for $\nBasisVelocityROM$ influences the number of non-linear iterations and, in some cases, slows down the solution process. Nevertheless, the order of magnitude of the observed speed up, ranging from 25,000 to 60,000, is still compelling and even much more significant than in the previous test case.
\par
So, also the results for this test case underline the advantages of the proposed approach in general and, specifically, its ability to handle problems defined in a four-dimensional space-time domain. 

\section{Conclusion}
\label{sec:conclusion}
In this work, we presented an \gls{mor} concept for parametric time-dependent problems defined in deforming domains that are allowed to even entail spatial topology changes. There are two main building blocks in this approach. First, we make use of the time-continuous space-time setting. Here, the \gls{csst}-\gls{fem} leads to fixed finite-dimensional subspaces for the entire space-time domain, implicitly accounting for current domain deformations. Working on these fixed subspaces, we make use of a projection-based \gls{mor} technique based on \gls{pod}. In contrast to other \gls{mor} methods applied to deforming domain problems, this particular approach proposed in this paper can be applied in a straightforward manner. Taken together, we argued that it can reduce the computational complexity for the aforementioned class of problems and, at the same time, maintain the desired level of accuracy in the results.
\par
To confirm this claim, we investigated two representative test cases from the fields of engineering and biomedical applications. In particular, we carried out error and performance analyses for both cases. The first test case was composed of a two-dimensional valve-like geometry containing a moving plug. In the resulting three-dimensional space-time domain, we considered the flow of plastics melt, which occurs in many polymer processing techniques. The parametric character of the problem originated from incorporating uncertainty in the viscosity model, i.e., fluctuations in the model parameters. The error analysis showed that we are able to reduce the original model effectively while keeping control over the magnitude of the error in the relevant physical fields. Moreover, a significant speed-up in terms of computational time was proven.
By means of the second test case, we extended our analysis to four-dimensional space-time domains that arise when considering three-dimensional bodies in space. Specifically, we studied blood flow in an artery-like volume that undergoes compression, yielding the domain deformation. The problem was varied via a parametric prescribed inflow velocity. The results of the error and performance analysis confirmed the findings from the previous test case and demonstrates the applicability of the approach, also when interested in spatially three-dimensional problems.
\par
Overall, the approach presented here extends the collection of existing \gls{mor} techniques with a simple but elegant approach for deforming domain problems. Despite the fact that it shows methodologically favorable characteristics and yields plausible results for the presented test cases, practical difficulties may arise. For example, in three-dimensional application cases, the generation of the required four-dimensional meshes can be challenging and limits the applicability as discussed in~\Cref{subsec:methodologicalBackgroundFOM}. Furthermore, the nature of the resulting reduced system may influence and potentially limit the speed-up that can be realized. Especially, the latter observation remains to be elucidated. Nevertheless, the proposed approach is able to make the advantages of \gls{mor} accessible for the class of complex problems that include domain deformations, potentially with spatial topology changes.

\section*{Acknowledgments}
\label{sec:acknowledgments}
MD gratefully acknowledges the funding by dtec.bw \--- Digitalization and Technology Research Center of the Bundeswehr under the project RISK.twin. dtec.bw is funded by the European Union – NextGenerationEU.
FB and GR acknowledge the European Union's Horizon 2020 research and innovation program under the Marie Skłodowska-Curie Actions, grant agreement 872442 (ARIA).
FB also thanks the project ``Reduced order modelling for numerical simulation of partial differential equations'' funded by Università Cattolica del Sacro Cuore, and the INDAM-GNCS project ``Metodi numerici per lo studio di strutture geometriche parametriche complesse'' (CUP\_E53C22001930001, PI Dr. Maria Strazzullo).

\printbibliography

\end{document}